\documentclass[10pt]{article}
\usepackage[margin=0.6in]{geometry}
\usepackage{hyperref}
\usepackage{multicol}
\usepackage{blindtext}
\usepackage{amsmath}
\usepackage{amssymb}
\usepackage{mathrsfs}
\usepackage{dsfont}
\usepackage{framed}
\usepackage{algorithm}
\usepackage{algpseudocode}
\usepackage{authblk}
\usepackage{graphicx}
\usepackage{wrapfig}
\usepackage[table]{xcolor}
\usepackage{wrapfig}
\usepackage{float}
\usepackage{mdframed}
\usepackage{multirow}
\mdfdefinestyle{mddefault}{
rightline=true,innerleftmargin=4pt,innerrightmargin=4pt,innertopmargin=4pt,innerbottommargin=2pt
frametitlerule=false,backgroundcolor=black!5}

\usepackage{framed,color}
\definecolor{shadecolor}{rgb}{0.85,0.85,0.85}
\usepackage{ulem}
\usepackage{url}
\usepackage{cleveref}
\usepackage{hyperref,xcolor}
\hypersetup{
colorlinks,
linkcolor={red!50!black},
citecolor={blue!50!black},
urlcolor={blue!80!black}
}
\usepackage{xr}
\makeatletter
\newcommand*{\addFileDependency}[1]{
\typeout{(#1)}
\@addtofilelist{#1}
\IfFileExists{#1}{}{\typeout{No file #1.}}
}\makeatother

\usepackage{color}
\usepackage{ulem}
\usepackage{cancel}

\usepackage{lineno}

\usepackage{amsthm}
\usepackage{hyperref}
 \newcommand{\ind}{\perp\!\!\!\!\perp} 

\theoremstyle{definition}
\newtheorem{theorem}{Theorem}
\newtheorem{lemma}{Lemma}

\newtheorem{corollary}{Corollary}[theorem]
\newtheorem{result}{Result}
\newtheorem{limit_result}{Result}[result]

\newtheorem{definition}{Definition}
\newtheorem{fact}{Fact}

\newtheorem{equation_env}{Equation}

\RequirePackage{tikz} 
\newcommand{\orcidicon}{\includegraphics[width=0.32cm]{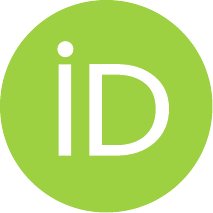}}

\foreach \x in {A, ..., Z}{%
\expandafter\xdef\csname orcid\x\endcsname{\noexpand\href{https://orcid.org/\csname orcidauthor\x\endcsname}{\noexpand\orcidicon}}
}

\begin{document}
\title{Shared-Endpoint Correlations and Hierarchy in Random Flows on Graphs}
\newcommand{\orcidauthorA}{0000-0001-6618-631X} 

\author{Joshua Richland $^1$ and Alexander Strang $^2$  \orcidA{} }

\affil{\footnotesize $^1$Stanford University, Statistics. \\ $^2$University of California, Berkeley, Statistics. \\ Correspondance: alexstrang@berkeley.edu}

\date{\vspace{-1.2em} \today}

\maketitle
\vspace{-30pt}

\section{Importance Statement}

Edge flows describe flux over the edges of a network. Fixing an edge flow is equivalent to fixing a directed, weighted version of an undirected graph. Accordingly, edge flows may be adapted to describe various alternating functions on graphs including preference in comparison networks, advantage in competitive networks, thermodynamic exchanges in biophysical networks, transition odds in discrete-state stochastic processes, or information transfer in neural networks. In these settings, the bulk organization of the flow is of interest. Does it circulate? Does it originate in sources and collect in sinks? Random models for flows allow a user to specify generative assumptions that produce varying degrees of organization. By varying assumptions, the user can study sufficient conditions needed to produce differing degrees of organization. When appropriately minimal, these assumptions may be cast as a null model that should be disproved before assuming more sophisticated generative processes. Thus, characterizing the functional relationship between the expected organization of a flow and a random generating model is important. For example, if the flow on an edge is only determined by its endpoints, and is otherwise independent of the broader topological context of the edge, then the bulk organization of flow is entirely determined by the correlation in flows on edges sharing an endpoint. We compute the shared-endpoint correlation for random edge flows generated by underlying Gaussian processes. We provide both sharp quantitative predictions for specific processes and identify generic trends. These trends share a standard qualitative moral; smoother functions exhibit higher correlation, thus produce more organized flows. 

\section{Abstract}



We analyze the correlation between randomly chosen edge weights on neighboring edges in a directed graph. This shared-endpoint correlation controls the expected organization of randomly drawn edge flows when the flow on each edge is conditionally independent of the flows on other edges given its endpoints. To model different relationships between endpoints and flow, we draw edge weights in two stages. First, assign a random description to the vertices by sampling random attributes at each vertex. Then, sample a Gaussian process (GP) and evaluate it on the pair of endpoints connected by each edge. We model different relationships between endpoint attributes and flow by varying the kernel associated with the GP. We then relate the expected flow structure to the smoothness class containing functions generated by the GP. We compute the exact shared-endpoint correlation for the squared exponential kernel and provide accurate approximations for Mat\'ern kernels. In addition, we provide asymptotics in both smooth and rough limits and isolate three distinct domains distinguished by the regularity of the ensemble of sampled functions. Taken together, these results demonstrate a consistent effect; smoother functions relating attribute to flow produce more organized flows.

\section{Background}

\subsection{Edge Flows}

An edge flow on an undirected graph $\mathcal{G} = (\mathcal{V},\mathcal{E})$ is a real-valued function $f: \mathcal{V} \times \mathcal{V} \rightarrow \mathbb{R}$ on the edges of the graph, where, for each pair of connected vertices $(i,j)$, $f(i,j) = -f(j,i)$ \cite{lim2020hodge}. Edge flows are the discrete analog to vector fields. 

Edge flows have been adopted to describe a plethora of processes on graphs. Most directly, an edge flow may represent the flow of some quantity, such as current, fluid flow, or network traffic (c.f.~\cite{zhanglonghodge}). There are a variety of applications in economics \cite{kichikawa2019community, fujiwara2020hodge,jiang2011statistical, wand2024causal, mike2019combinatorial}, game theory \cite{candogan2011flows, balduzzi2019open,balduzzi2018re,sizemore2013hodgerank,strang2022network}, as well as biophysics and molecular biology \cite{nartallo2024decomposing,strang2020applications}. Edge flows can also describe connectivity patterns in biological neural networks \cite{miura2015hodge,miura2015scaling}, aggregate preference data in pairwise comparison problems and social choice \cite{jiang2011statistical, strang2020applications}, control patterns for multi agent robot routing \cite{kingston2011distributed}, and even the structure of geopolitical sanction networks \cite{hisano2020identifying}. 

Often, it is helpful to decompose an edge flow into simpler components with well-defined properties. Questions about the global structure of the flow can be answered using these components, whose relative sizes summarize global flow structure \cite{jiang2011statistical,strang2022network}. For example, we might ask, does the flow circulate or does it tend to flow ``downhill" from sources to sinks?  If the flow concentrates, where? If it circulates, in what part of the graph, on loops of what size, and how intensely? If the network has open boundaries, is there any component passing through the boundaries? And, if the flow can be decomposed, how much of the flow is represented by each component? 

While interpretations vary by field, these underlying questions remain the same. In transaction networks, circulation may represent arbitrage \cite{jiang2011statistical}, while flow from sources to sinks can be used to separate businesses into hierarchical communities \cite{fujiwara2020hodge} and to target market interventions \cite{wand2024causal}. In a game theoretic setting, the optimal mixed strategies are influenced by circulation in utility \cite{candogan2011flows}, as are ensuing evolutionary dynamics \cite{Xue}. In biophysics and stochastic processes, the work done when a system transitions between states is an edge flow, and the presence or lack of cycles in the edge flow has deep thermodynamic consequences. Flows that do not admit any cyclic component correspond to equilibrium processes, which are time-reversible and reach thermodynamic equilibrium. Flows that admit circulation correspond to non-equilibrium processes, which characterize living systems that continuously transport energy between reservoirs to perform work \cite{schnakenberg1976network,strang2020applications}. In ranking and decision problems, edge flows represent dominance or preference relations, and flow from sources to sinks distinguishes high and low rank \cite{jiang2011statistical,strang2020applications,strang2022network}. Cyclic edge flows cannot be described by rank order, so cycles play an important role in social choice theory where they allow voting paradoxes. The extent of cyclicity in voter preferences determines how vulnerable a system is to strategic deal-making \cite{Lagerspetz}, reordering of choices \cite{Flanagan}, or agenda-setting authority \cite{Morse}. Thus, cycles frustrate fair choice. 

In all of the applications described above, the authors study the global structure of their edge flow via a decomposition into components that either flow downhill, circulate at a chosen scale, or do neither. The associated decomposition is the combinatorial Helmholtz-Hodge Decomposition (HHD) introduced in \cite{jiang2011statistical,lim2020hodge}. The combinatorial HHD adapts the familiar Helmholtz decomposition from vector calculus to suit edge flows on graphs. Like the Helmholtz decomposition, it separates an edge flow into a conservative (or transitive) component that flows from sources to sinks, a circulating component, and a harmonic component which neither circulates nor concentrates at sinks. 

In particular, if the graph $\mathcal{G}$ is closed and finite, then any flow $f$ admits a unique decomposition:
\begin{equation}
    f = f_t + f_c
\end{equation}
where $f_t$ is the transitive component of the flow and $f_c$ is the cyclic component of the flow \cite{strang2020applications}. 

Null distributions that generate flows are needed to evaluate the statistical significance of any observed flow structure. For example, the significance of observed degrees of hierarchy in dominance networks is typically evaluated by comparing against a uniform null \cite{Appleby,de_Vries,Kendall}. While uniform randomization is straightforward, it is plainly implausible in most applied settings. Indeed, uniform randomization produces predominantly cyclic networks for almost any dense network of even moderate size \cite{Gehrlein_d,Gehrlein_e,strang2022network}. For example, the probability that a directed graph with $V = |\mathcal{V}|$ vertices is acyclic, when sampled uniformly among all directed graphs with $V$ vertices, converges to zero with rate at least $V^{-1/3}$ \cite{ralaivaosaona2020probability}. If the edge directions are sampled uniformly and independently, then, on average, a quarter of all triangles form cycles \cite{Shizuka}. Demonstrating significance against a ``strawman" null is necessary, but far from sufficient, if the author aims to establish evidence for a more complicated underlying process. The space of possible random flow models that are not uniform is vast, varied, and deserves more nuanced exploration.

\subsection{Null Models and Shared-Endpoint Correlations}

Following \cite{strang2022network}, we focus on generative models that produce randomly sampled flows with tunable degrees of structure. Note that, since each flow is a directed quantity defined on the edges, sampling a flow is equivalent to sampling a weighted-directed graph on $\mathcal{V}$, with a directed edge for every undirected edge in $\mathcal{E}$. When the edge-set $\mathcal{E}$ is also random, then sampling an edge set and a flow is equivalent to sampling a random, weighted, directed graph.

In particular, we consider random models for both the graph $\mathcal{G}$ and its flow $f$, where, conditional on $\mathcal{G}$, the value of $f$ on a given edge, $(i,j)$, depends only on the identities of the endpoints, and is otherwise independent of their position in the rest of the network. This model can be parameterized by first sampling or fixing the graph of interest, $\mathcal{G}$, then
\begin{enumerate}
    \item choosing an attribute space $\Omega$, where $x \in \Omega$ is a list of attributes used to describe a vertex,
    \item fixing a distribution $\pi_X$ over $\Omega$, and drawing an attribute for each node, $X(i) \sim \pi_X$, i.i.d., 
    \item fixing a function, $f:\Omega \times \Omega \rightarrow \mathbb{R}$ satisfying $f(x,y) = -f(y,x)$, 
    \item then setting the flow on each edge to the function $f$ evaluated at its endpoints, $F(i,j) = f(X(i),X(j))$. 
\end{enumerate}

This model was introduced in \cite{strang2022network}, to study the expected structure of competitive tournaments. There, attributes were traits of individual competitors, and the flow modeled competitive advantage, i.e. performance, between individuals. Accordingly, the name trait-performance was used to describe models satisfying the sampling process (1) - (4). 

Models in this family specify a distribution over flows, $F$, that are minimally structured. The marginal distribution of the flow on all edges are identical, the joint distribution over all pairs of edges are identical, and the flows on two disjoint edges are independent. Indeed, the flow on any subgraph is independent of the flow on any other subgraph if the subgraphs do not share any nodes, and are identical if the subgraphs are topologically equivalent. Nevertheless, these models do fix some structure in the joint distribution, since edges sharing an endpoint are correlated. The topology of the graph $\mathcal{G}$ influences the flow distribution through this coupling. The strength of the coupling, i.e. the correlation in the flow between adjacent edges, is a useful summary statistic that predicts structural characteristics of the overall flow. 

\begin{snugshade}
\begin{theorem} \label{thm: trait performance} \textbf{[Trait-Performance] \cite{strang2022network}}  
Let $\mathcal{G} = (\mathcal{V},\mathcal{E})$ denote a random undirected graph with $V = |\mathcal{V}|$ vertices and $E = |\mathcal{E}|$ edges. Let $L = E - (V - 1)$ denote the dimension of the cycle space of the graph. Let $F$ denote a random flow where $F(i,j) = f(X(i),X(j))$ for $X \sim \pi_X$ drawn i.i.d. Then, the covariance matrix for the flow is of the form:
\begin{equation} \label{eqn: covariance structure}
    \mathbb{C}[F] = \sigma^2(I_{E \times E} + \rho A_{\mathcal{E}}(\mathcal{G}))
\end{equation}
where $\sigma^2 = \text{Var}_{X,Y \sim \pi_X}[f(X,Y)]$ is the variance in the flow on each edge, $\rho = \text{Corr}_{X,Y,W \sim \pi_X}[f(X,Y),f(X,W)]$ is the correlation in the flows on edges sharing an endpoint, $I_{E \times E}$ is the $E \times E$ identity matrix, and $A_{\mathcal{E}}(\mathcal{G})$ is the signed edge-adjacency matrix which is an indicator matrix for connected edge pairs.

\vspace{0.1 in}

\noindent Let $F_t$ and $F_c$ denote the transitive and cyclic components of $F$. Then:
\begin{equation} \label{eqn: trait performance}
    \begin{aligned}
        &\mathbb{E}[\|F_t\|^2] = \sigma^2 \left((\mathbb{E}[V] - 1) + 2 \rho \mathbb{E}[L] \right) \quad \quad \mathbb{E}[\|F_c\|^2] = \sigma^2 (1 - 2 \rho) \mathbb{E}[L] \\
    \end{aligned}
\end{equation}
where $\rho \in [0,1/2]$ and $F_c = 0$ with probability one if $\rho = 1/2$.
\end{theorem}
\end{snugshade}

Equations \eqref{eqn: covariance structure} and \eqref{eqn: trait performance} illustrate that, despite the extreme flexibility available in $\Omega$, $\pi_X$, and $f$, the covariance of the flow and the expected sizes of its components are controlled by only two statistics: the variance in the flow on each edge, which acts as a scaling factor, and the shared-endpoint correlation, $\rho$, which determines the relative sizes of the expected components. More strikingly, only very basic features of the graph influence the expected sizes of the components: the number of edges and vertices. All other topological details are irrelevant. In fact, not even both numbers are needed if we are interested in the relative sizes of the components (c.f.~\cite{wand2024causal}). The relative size of the expected components depends exclusively on the shared-endpoint correlation $\rho$ and the average density of the graph, $\mathbb{E}[E]/\mathbb{E}[V]$. 

The expected degrees of hierarchy and cyclicity, are monotonic in both the shared-endpoint correlation and the density. On average, denser graphs are more cyclic and less transitive. 
Similarly, the larger the correlation, the more transitive the graph, and the smaller the correlation, the less transitive. In the extreme case when $\rho$ is maximized ($\rho = 1/2$), all flows are perfectly transitive. For details and interpretation in a game theory setting, see \cite{strang2022network}.

The shared-endpoint correlation is particularly interesting since the random flow model is consistent upon subdivision of the graph into smaller, disjoint parts. For example, if a user was interested in the expected degree of cyclicity on a particular motif, or on local subsets, then changing the subgraph would change the expected density, but would not change the correlation coefficient. Since the density of a fixed subgraph, or graphlet, is known, the only variable needed to predict the relative cyclicity of flows on that subgraph is the shared-endpoint correlation (see equation \eqref{eqn: trait performance}). In this sense, the shared-endpoint correlation offers a topology-independent characterization of expected flow structure. 

If a user is not interested in the sizes of the HHD components, then there are still good reasons to study the shared-endpoint correlation. For example, given a fixed graph, the covariance matrix for a flow has only two free parameters, $\sigma^2$ and $\rho$. The first is simply a scale parameter, and can be normalized by adjusting the units used to measure the flow. Then up to scale, the form for the covariance matrices depends only on $\rho$. Thus, any downstream quantity that depends on the covariance, independent of scale, will only depend on $\rho$ and the graph topology.

For these reasons, we will focus on the shared-endpoint correlation, $\rho$, as a unit-less, topology-free statistic that summarizes the expected degree of structure in randomly sampled flows. Large $\rho$ suggests organization; flows are strongly coupled on neighboring edges and move predominantly from sources to sinks. Small $\rho$ indicates disorganization. In the extreme case, $\rho = 0$, the flows on all edges are independent and identical, so we recover the standard uniform model. 

The shared-endpoint correlation is determined by the attribute space $\Omega$, attribute distribution $\pi_X$, and the function $f$. In particular, the correlation is strongly influenced by the smoothness of $f$ and breadth of $\pi_X$. If $f$ is linear, then $\rho = 1/2$. When $f$ is nearly linear for most $X$ drawn from $\pi_X$, then $\rho$ is near to $1/2$. If $f$ is differentiable, then it is nearly linear on local neighborhoods, so $\rho$ converges to $1/2$ if the expected errors in linear approximations vanish sufficiently quickly as $\pi_X$ concentrates. Detailed conditions for convergence to $1/2$ and rates are provided in \cite{cebra2023similarity}. In particular:
\begin{snugshade}
\begin{fact}\label{fact 1} \textbf{[Smoothness, Concentration, and Correlation Limits]} \cite{cebra2023similarity}
\vspace{-0.05 in}
\begin{enumerate}
    \item $\rho$ almost surely converges to $1/2$ as $\pi_X$ concentrates if $f$ is differentiable almost everywhere with respect to $\pi_X$,
    \vspace{-0.05 in}
    \item if $f$ is also second differentiable almost everywhere, then $1/2 - \rho$ converges to zero at rate $\mathcal{O}(\kappa^2)$ or faster provided the covariance in $\pi_X$ is $\mathcal{O}(\kappa^2)$,
    \vspace{-0.05 in}
    \item and, if $\Omega \subset \mathbb{R}$, that is, if the attribute-space is one-dimensional, then $1/2 - \rho$ converges to zero at rate $\mathcal{O}(\kappa^4)$ or faster.
\end{enumerate}
\end{fact}
\end{snugshade}

Under reasonable conditions, Fact \ref{fact 1} implies both convergence to perfect transitivity ($F_c = 0$) in expectation, and convergence to transitivity (the directed graph is acyclic) in probability when the $\pi_X$ concentrates \cite{cebra2024almost}. Following \cite{cebra2023similarity} and \cite{cebra2024almost}, we aim to investigate the dependence of $\rho$ on both the similarity of the sampled attributes and on the smoothness of $f$.
The analyses in \cite{cebra2023similarity} and \cite{cebra2024almost} provide guaranteed convergence rates up to order in generic concentration parameters under strict constraints on the regularity of $f$. 
This approach has two major limitations. 

First, it can only provide rates of convergence. It cannot provide specific predictions for $\rho$ since no distribution is introduced over the set of functions satisfying the enforced constraints. Second, it does not provide predictions for intermediate degrees of regularity. Since the conclusions are based on integer degrees of differentiability, there is no method for interpolating between functions that are first but not second, or second but not third differentiable almost everywhere. There is also no way to study functions that are close to, but not, differentiable. To link degrees of differentiability, we should work with ensembles of functions that are almost surely differentiable up to a fractional degree. See \cite{teodoro2019review} for a review of definitions. Then, it may be possible to study the behavior of $\rho$ for functions that are almost but not quite firs or second differentiable.

In this paper, we adopt an alternative tactic. Instead of fixing a set of constraints on the regularity of $f$, we fix a prior distribution over possible functions $f$. The choice of prior may be viewed as either a method for averaging over an ensemble of plausible functions, a generative model, or as an added sampling step. By varying the distribution over $f$, we can prioritize functions of varying degrees of smoothness or regularity without strictly enforcing constraints. 

Let $\mathcal{F}$ denote the randomly sampled function, $\mathcal{X}$ the randomly sampled set of attributes assigned to each vertex, and $\mathcal{G}$ the randomly sampled network. As long as $\mathcal{G}$ is sampled independently of the attributes, and thus, the flow, the trait-performance theorem still holds in expectation over the random function draw:
\begin{equation} \label{eqn: trait performance random f and G}
\mathbb{E}_{\mathcal{F},\mathcal{X},\mathcal{G}} \left[||F||^2 \right] = \sigma^2  \xrightarrow{\text{decompose}}
\left\{
\begin{aligned}
& \mathbb{E}_{\mathcal{F},\mathcal{X},\mathcal{G}} \left[||F_t||^2 \right] = \sigma^2 \left[ (\mathbb{E}_{\mathcal{G}}[V] - 1) + 2 \rho \mathbb{E}_{\mathcal{G}}[L] \right] \\
& \mathbb{E}_{\mathcal{F},\mathcal{X},\mathcal{G}} \left[||F_c||^2 \right] = \sigma^2 \left(1 - 2 \rho \right) \mathbb{E}_{\mathcal{G}}[L] \end{aligned} \right.
\end{equation}
where:
\begin{equation} \label{eqn: sigma and rho for random f}
    \begin{aligned}
        & \sigma^2 = \mathbb{E}_{\mathcal{F}}\left[\mathbb{V}_{\mathcal{X}|\mathcal{F}}[\mathcal{F}(\mathcal{X}_i,\mathcal{X}_j)] \right] = \mathbb{E}_{\mathcal{F}}\left[\mathbb{V}_{X,Y\sim \pi_x|\mathcal{F}}[\mathcal{F}(X,Y)] \right] \\
        & \rho = \frac{1}{\sigma^2} \mathbb{E}_{\mathcal{F}}\left[\mathbb{C}_{\mathcal{X}|\mathcal{F}}[\mathcal{F}(\mathcal{X}_i,\mathcal{X}_j),\mathcal{F}(\mathcal{X}_i,\mathcal{X}_k)] \right] = \frac{1}{\sigma^2} \mathbb{E}_{\mathcal{F}}\left[\mathbb{C}_{X,Y,W\sim \pi_x|\mathcal{F}}[\mathcal{F}(X,Y),\mathcal{F}(X,W)] \right]
    \end{aligned}
\end{equation}
where $\mathbb{V}$ and $\mathbb{C}$ stand for variance and covariance respectively. 

If the attributes $X$ are independent of the flow function $\mathcal{F}$, then the conditional statements can be removed:
\begin{equation} \label{eqn: sigma and rho for random f independent X}
    \begin{aligned}
        & \sigma^2  = \mathbb{E}_{\mathcal{F}}\left[\mathbb{V}_{X,Y\sim \pi_x}[\mathcal{F}(X,Y)] \right] \\
        & \rho = \frac{1}{\sigma^2} \mathbb{E}_{\mathcal{F}}\left[\mathbb{C}_{X,Y,W\sim \pi_x}[\mathcal{F}(X,Y),\mathcal{F}(X,W)] \right]
    \end{aligned}
\end{equation}
and the order of the expectations over function and population can be exchanged, allowing simplification.

In this paper, we will attempt to calculate $\rho$, as defined in equation \eqref{eqn: sigma and rho for random f} for $\mathcal{F}$ drawn from a Gaussian process (GP). In the next section, we will review GP's, and argue that they offer a useful compromise between generality and tractability when constructing random flow models. 

\subsection{Gaussian Processes} \label{sec: GPs}

A Gaussian process (GP) is a stochastic process satisfying the condition that any finite subset of samples from the process are jointly Gaussian distributed \cite{williams2006gaussian}. A GP is parameterized by a mean and covariance function. Formally, $\mathcal{F} \sim \text{GP}(\mu,k)$ where $\mathcal{F}:\Omega \rightarrow \mathbb{R}$ for some domain $\Omega$, $\mu:\Omega \rightarrow \mathbb{R}$ is the mean function, $\mathbb{E}[\mathcal{F}(x)] = \mu(x)$, and $k:\Omega \times \Omega \rightarrow \mathbb{R}$ is the covariance function, $k(x,y) = \mathbb{C}[\mathcal{F}(x),\mathcal{F}(y)]$. Note that the covariance function must be symmetric, $k(x,y) = k(y,x)$ and must produce positive semi-definite matrices for any collection of sample locations. For a collection of locations, $\textbf{X} = \{x_i\}_{i=1}^n$, we write $K_{\textbf{XX}}$ for the $n\times n$ matrix with entries ${K_{\textbf{XX}}}_{ij} = \kappa(x_i,x_j)$.

We adopt Gaussian processes for three reasons:
\begin{enumerate}
\item They are a sufficiently broad class of models. Any linear combination of basis functions with Gaussian coefficients defines a Gaussian process. Famous examples include Brownian motion, Ornstein-Uhlenbeck processes, Gaussian random fields, and autoregressive processes with Gaussian noise \cite{williams2006gaussian}.

In fact, they are representative. When the attributes are drawn independently of $\mathcal{F}$, $\rho$ depends only on the mean and covariance functions of the process, whether or not it is a GP. Thus, the correlation coefficient for $\mathcal{F}$ will equal the correlation coefficient for a GP, $\mathcal{F'}$, with matching mean and covariance functions (see Section \ref{sec: beyond GPs}). From this perspective, we need only fix a mean and covariance function for an otherwise unspecified process, and can adopt the corresponding GP as a representative. This approach can be justified via Occam's razor, as $\mathcal{F'}$ is the maximum entropy process given a known mean and covariance function for a process taking values in $\mathbb{R}$ \cite{cover1999elements}. This approach is further justified since GP's are entirely specified by their mean and covariance functions, so do not require any additional information beyond what is needed to recover $\rho$. 

\item The parameters (mean and covariance) are interpretable. In particular, the smoothness of a GP is controlled by the choice of the covariance function \cite{adler2010geometry}. When stationary, the differentiability of the covariance function at nearby inputs \cite{williams2006gaussian}, and the tail decay rate of its Fourier transform \cite{steinwart2019convergence}, control sample regularity. They control the degree of differentiability of sample draws, the reproducing kernel Hilbert space that almost surely contains draws, and the smaller reproducing kernel Hilbert space that contains posterior estimators \cite{dudley2010sample,williams2006gaussian}. It is straightforward to adopt covariance functions such that draws are almost surely $\nu$-differentiable for any $\nu \in \mathbb{R}^+$. These results offer clear guidance to a modeler interested in sample regularity. Note: while the relation between covariance function and correlation are independent of the choice to use a GP, the implied smoothness properties of draws are not.

\item GP's admit standard reference models. The existence of standard references allows the analyst to prioritize  effort. In this sense, our work is GP specific since we prioritize standard GP models. Accordingly, we develop the main text assuming that $\mathcal{F}$ is a GP, before generalizing in Section \ref{sec: beyond GPs}.
\end{enumerate}


These features strongly recommend GP's for exploratory modeling.  We will encode smoothness assumptions in the choice of covariance function, then study how the shared-endpoint correlation $\rho$ changes as we vary both the length scale of variation and the regularity of the class of sample draws. In contrast to previous work, which relied on local approximations and only provided convergence rates \cite{cebra2023similarity}, we will provide exact correlations.

We will, in most cases, assume two standard symmetries: stationarity and isotropy. A stochastic process is stationary if its joint distribution over any set of samples is invariant to translation of the set. Then, the ensemble is translation-symmetric. A stochastic process is isotropic if the joint distribution is also invariant under unitary transformations of the coordinates used to represent samples. For example, the ensemble should be rotation and reflection-symmetric. 
In this case, all input directions are treated identically. In the following analysis, we will consider anisotropy where it is revealing to do so. All examples are stationary. 

Stationarity and isotropy are often adopted for analytic convenience. In part, our motivations follow. Given stationarity and isotropy, our results can be expressed cleanly in closed form. However, convenience alone is insufficient justification. 

Recall that the shared-endpoint correlation depends on both draws of the sample function $\mathcal{F}$, and draws of a sample population $\mathcal{X}$ representing the set of attributes of nodes in a graph. We need not assume that the underlying process generating $\mathcal{F}$ is stationary or isotropic, or that the joint distribution of $\mathcal{X}$ is isotropic, to employ either assumption. We need only assume that the sampled function and attributes are unrelated. Suppose, for example, that $\mathcal{F}$ is not drawn from a stationary process. Suppose also that $\mathcal{X}$ is drawn from a distribution $\pi_X$ whose centroid is subject to a uniform random shift. Averaging over the possible centroid locations is equivalent to fixing the centroid, and averaging the process $\mathcal{F}$ over all possible translations. The translation-averaged process is stationary. So, if the location of the trait distribution is independent of $\mathcal{F}$, and otherwise uniform, then we can freely center $\pi_X$ at zero, and adopt a stationary process. Isotropy can be justified similarly. If there is no reason to expect $\pi_X$ to align its orientation to the process, or if $\pi_X$ is itself isotropic, then averaging over the possible rotations of $\pi_X$ is equivalent to fixing an isotropic trait distribution and process. 

Neither assumption need hold when the attributes have developed in response to the flow. For example, in a competitive setting, the flow may represent advantage. If agents evolve, then the process will exert selective pressure on the attributes, producing trait distributions that depend on the process, seek its critical points, and align with its local curvature (c.f.~\cite{cebra2024almost}). We save this case for future work. 

If a Gaussian process is stationary, then its covariance function, $k(x,y)$, depends only on the difference in the inputs, $x - y$. In particular, there exists a kernel function $\kappa: \Delta \Omega  \rightarrow \mathbb{R}$ such that $k(x,y) = \kappa(x - y)$ for all $(x,y) \in \Omega \times \Omega$. Since $k(x,y) = k(y,x)$ the kernel function is always even. Most kernel functions are further parameterized by adopting some decaying function of distance $h:\mathbb{R}^+ \rightarrow \mathbb{R}$, and norm $\| \cdot \|$, such that $\kappa(x - y) = h(\|x - y\|)$. The process is stationary and isotropic if and only if the norm is the Euclidean norm $\| \cdot \|_2$.

Popular examples include the exponential and squared exponential (SE) kernels:
\begin{equation}  
    \kappa_{\text{exp}}(x - y;l) = h_{\text{exp}}(\|x - y\|;l) =  \exp \left(- \frac{\|x - y\|}{l} \right), \quad  \kappa_{\text{SE}}(x - y;l) = h_{\text{SE}}(\|x - y\|;l) = \exp \left( -\frac{1}{2} \left(\frac{\|x - y\|}{l}\right)^2\right)
\end{equation}
which are both parameterized by a length scale $l$. The squared exponential kernel produces highly regular functions. Draws from the squared exponential kernel are infinitely differentiable with probability one. In contrast, the exponential kernel produces rough sample functions. Draws from the exponential kernel correspond to draws of an Ornstein-Uhlenbeck (OU) process in one dimension \cite{uhlenbeck1930theory,maller2009ornstein,williams2006gaussian}. Draws of an OU process are nowhere differentiable. 

The exponential and squared exponential are special cases of the Mat\'ern model \cite{stein1999interpolation,matern1960spatial}: 
\begin{equation} \label{eqn: Matern}
\kappa_{\text{Mat}}(x - y;l,\nu) = h_{\text{Mat}}(\|x - y\|;l,\nu) =  \frac{2^{1 - \nu}}{\Gamma(\nu)} \left(\frac{\sqrt{2\nu}\|x - y\|}{l} \right)^{\nu} K_{\nu}\left(\frac{\sqrt{2\nu}\|x - y\|}{l} \right)\end{equation}
where $K_{\nu}$ is the modified Bessel function \cite{abramowitz1965formulas,williams2006gaussian}.

\begin{figure}[t]
    \centering
    \includegraphics[scale = 0.4]{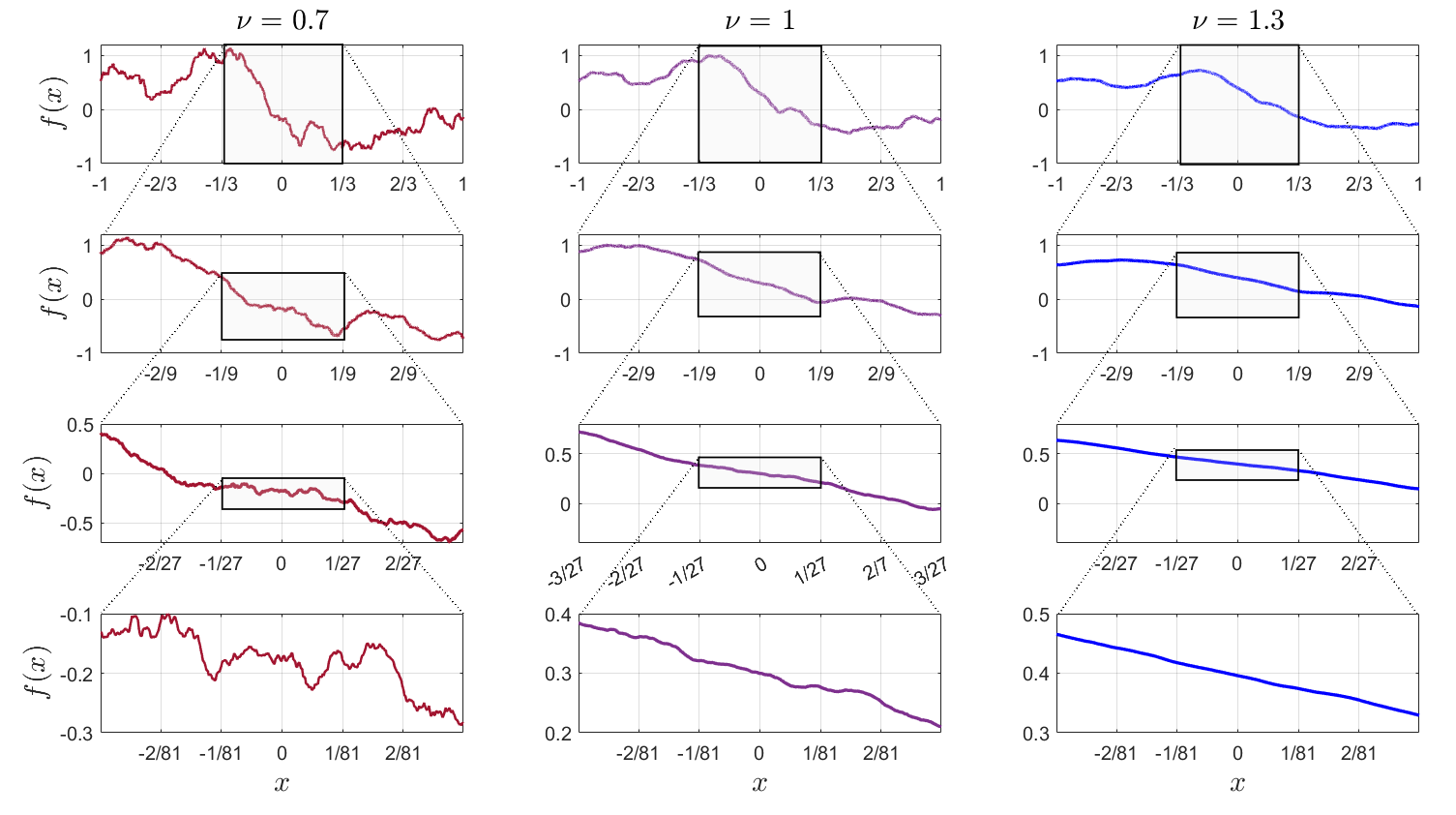}
    \caption{ Sample trajectories drawn from a Mat\'ern process with length scale parameter $l = 1$ and varying regularity parameter $\nu$. All columns use the same random seed. Moving down between rows ``zooms in." Each row scales the $x$ axis of the row above it by a factor of 3. The boxes show the region expanded. All boxes in a given row have the same height.} 
    \label{fig: Matern regularity}
\end{figure}

The Mat\'ern kernel is parameterized by a length scale $l$, and a shape parameter $\nu$ which determines the regularity of draws of the process. When $\nu = \frac{1}{2}$, the process is exponential. When $\nu$ approaches infinity, the process converges to a squared exponential process. Figure \ref{fig: Matern regularity} shows three example trajectories for varying $\nu$.

We will focus on Mat\'ern kernels. To provide intuition, we recall some fundamental relations between the regularity of sample draws of a stationary GP and the Fourier transform of its kernel function. The Fourier transform of the kernel is the power spectral density (PSD). The PSD plays a central role in the analysis of stationary Gaussian processes (c.f.~\cite{bochner2005harmonic,gikhman2004theory,shannon1948mathematical}) since the Fourier transform, in effect, diagonalizes stationary kernels \cite{daw2022overview,karhunen1946spektraltheorie,loeve1955probability}. Loosely, the faster the PSD decays, the less power is assigned to high frequencies, so the more regular the process. In particular, the PSD of the Mat\'ern kernel decays at rate $\mathcal{O}(\|\zeta\|^{-(2 \nu + d)})$ where $\zeta$ is the frequency variable and $d$ is the dimension of the process. Processes whose PSD decays at this rate are almost always differentiable up to order $\nu$, but not beyond:

\begin{snugshade}
\begin{fact} \label{fact 2}
    \textbf{[Regularity and $\nu$]} If $\mathcal{F}$ is drawn from a Mat\'ern process with regularity parameter $\nu$, then $\mathcal{F}$ is $\mu^{th}$ order mean-square differentiable for all $\mu < \nu$, but no $\mu \geq \nu$ \cite{williams2006gaussian}. Sample draws of the process are almost surely contained inside the Sobolev space $W^{m}$ if $m < \nu$ \cite{kanagawa2018gaussian,steinwart2019convergence}, and are almost never contained in $W^{m}$ if $m \geq \nu$ \cite{driscoll1973reproducing}.
\end{fact}
\end{snugshade}

For detailed definitions, see Appendix Section \ref{app: GP regularity}.

Figure \ref{fig: Matern regularity} illustrates this effect as $\nu$ crosses 1. It shows three example trajectories sampled with $\nu = 0.7$, 1, and 1.3, while zooming into increasingly small neighborhoods. When $\nu < 1$ the sampled trajectory remains far from linear, while, for $\nu > 1$, the process converges to a line on small scales. This distinction will be important when studying the shared-endpoint correlation for attributes sampled from a small neighborhood, or for processes with a long length scale $l$. 

\subsection{Outline and Aims}
In summary, we aim to study the shared-endpoint correlation, $\rho$, when $\mathcal{F}$ is drawn from a GP with either a SE or a Mat\'ern kernel. We use the SE case to study the role of the length scale $l$ and the dimension of $\Omega$, for isotropic and anisotropic models (see Section \ref{sec: SE}). We generalize these results to any stationary kernel that can be expressed as a scale mixture of SE kernels. We leverage the mixture result to study $\rho$ when $\mathcal{F}$ is drawn from a GP with a Mat\'ern kernel (see Section \ref{sec: Matern}). We use the Mat\'ern model to isolate the influence of regularity, $\nu$, focusing on the critical transitions when $\nu$ crosses 1 and 2. We then show that the qualitative relationship between shared-endpoint correlation and regularity observed for the Mat\'ern kernel generalizes to arbitrary stationary isotropic kernels (see Section \ref{sec: generalization regularity}). In all cases, we develop asymptotic approximations to $\rho$ in the limit of smooth and rough functions. We use the smoothness limits to confirm the convergence rate analysis presented in \cite{cebra2023similarity}, while extending those results to fractionally differentiable processes. We also provide exact expectations, bounds, and approximations to $\rho$ away from either limit. 

\section{Results} \label{sec: results}

\subsection{Gaussian Process Models for Skew Symmetric Functions}

Suppose that the flow function $\mathcal{F}$ is drawn from a GP over the product space $\Omega \times \Omega$, where $\Omega \subset \mathbb{R}^T$. Here $T$ denotes the number of traits or attributes used to characterize endpoints. 

The flow function $f:\Omega \times \Omega \rightarrow \mathbb{R}$ must be skew symmetric: $f(x,y) = -f(y,x)$ for all $(x,y) \in \Omega \times \Omega$. Enforcing this symmetry on all draws $\mathcal{F}$ from a GP constrains both the mean and covariance functions. In particular, it will rule out most standard covariance functions, thus will require a minorly adapted construction.


\begin{snugshade}
\begin{theorem} \label{thm: skew symmetric GP}
\textbf{[Skew Symmetric GP]} 
Suppose $\mathcal{F} \sim \textnormal{GP}(\mu_f,k_f)$, $\mathcal{F}:\Omega \times \Omega \rightarrow \mathbb{R}$. Then $\mathcal{F}(x,y) = - \mathcal{F}(y,x)$ for all $(x,y) \in \Omega \times \Omega$ if and only if $\mu_f(x,y) = -\mu_f(y,x)$ and $k_f([x,y],[x,y]) = -k_f([x,y],[y,x])$ for all $x,y$.
\end{theorem}
\end{snugshade}

\textbf{Proof Outline:} Necessity is established by assuming that $\mathcal{F}(x,y) = - \mathcal{F}(y,x)$, then computing the moments on each side of the stated equalities. Sufficiency is established by showing that under the stated moment equalities $\mathcal{F}(x,y) + \mathcal{F}(y,x)$ has zero variance and mean zero, thus  $\mathcal{F}(x,y) = - \mathcal{F}(y,x)$ with probability one. See Appendix Section \eqref{app: thm 2}. $\blacksquare$
\vspace{0.1 in}

Theorem \ref{thm: skew symmetric GP} constrains the space of allowed covariance functions. First, notice that the covariance function cannot be strictly positively valued, since if it is strictly positively valued, then $k_f([x,y],[x,y]) > 0$ so $k_f([x,y],[y,x]) < 0$. Indeed, the kernel can only be nonnegative if $k_f([x,y],[x,y]) = 0$ for all $x,y$. This would, in turn, require, $\mathbb{V}[\mathcal{F}(x,y)] = 0$ for all $x$ and $y$, in which case $\mathcal{F}(x,y) = \mu(x,y)$ and there is no randomness in the choice of performance function (See Appendix Section \ref{app: cor 2.1}). Thus:

\begin{snugshade}
\begin{corollary} \label{cor 1: kernel sign}
    \textbf{[Sign of Kernel]} If $\mathcal{F} \sim \textnormal{GP}(\mu_f,k_f)$ for $\mu_f,k_f$ such that $\mathcal{F}(x,y) = -\mathcal{F}(y,x)$, and $\mathcal{F}(x,y) \neq \mu(x,y)$ then $k_f$ must take on both positive and negative values, and $k_f([x,y],[y,x])$ is negative where $\mathcal{F}(x,y)$ is nondeterministic.
\end{corollary}
\end{snugshade}


Theorem \ref{thm: skew symmetric GP} also rules out stationary processes.

\begin{snugshade}
\begin{corollary} \label{cor 2: nonstationarity}
    \textbf{[Nonstationarity] } If $\mathcal{F} \sim \textnormal{GP}(\mu_f,k_f)$ for $\mu_f,k_f$ such that $\mathcal{F}(x,y) = -\mathcal{F}(y,x)$, and $\mathcal{F}(x,y) \neq \mu(x,y)$, then $\mathcal{F}$ is not stationary.
\end{corollary}
\end{snugshade}

\textbf{Proof Outline:} Proceed by contradiction. If stationary, then $\mathbb{V}[\mathcal{F}(x,y)] = k_f([x,y],[x,y]) = \kappa_f(x - x,y - y) = \kappa_f(0,0)$, for all $(x,y) \in \Omega \times \Omega$. So, $\mathbb{V}[\mathcal{F}(x,y)] = \kappa_f(0,0) = \kappa_f(z-z,z-z) = \mathbb{V}[\mathcal{F}(z,z)]$ for any $z$. But, if skew-symmetric, $\mathcal{F}(z,z) = 0$ for any $z \in \Omega$, so $\mathbb{V}[\mathcal{F}(z,z)] = 0$. Then, $\mathbb{V}[\mathcal{F}(x,y)] = \mathbb{V}[\mathcal{F}(z,z)] = 0$ for all $x,y$, so $\mathcal{F}$ is deterministic. $\blacksquare$     
\vspace{0.1 in}

Corollary \ref{cor 2: nonstationarity} presents an interpretative challenge. In Section \ref{sec: GPs}, we argued that it is reasonable to assume stationarity and isotropy if the distribution $\pi_X$ is independent from $\mathcal{F}$. Yet, no skew-symmetric GP is ever stationary.

Recall that $\mathcal{F}$ is defined on the product space $\Omega \times \Omega$ consisting of pairs of attributes. Since $\mathcal{F}$ represents an interaction between two endpoints that is consistent under reversal of the order of the inputs, there is no reason to adopt one description for the endpoint listed first, and a different description for the endpoint listed second. Therefore, we need not consider transformations of $\Omega \times \Omega$ that treat the two inputs differently. Indeed, since we will draw all attributes from $\pi_X$, changes in $\pi_X$ will change both inputs identically. Shifts of $\pi_X$ can be absorbed as translations of $\Omega$. Rotations of $\pi_X$ can be absorbed as rotations of $\Omega$. Therefore, $\mathcal{F}$ need not be invariant to all transformations of $\Omega \times \Omega$, only those which transform the coordinates used for the first and second inputs equivalently. 

Thus, we only need to enforce a more limited form of stationarity in the product space. Namely, the joint distribution of any set of pairs in the product space $\{\mathcal{F}(x_i,y_i)\}$ should be unchanged by any shift of the form $\{\mathcal{F}(x_i + s,y_i + s)\}$. This is weaker than requiring stationarity as a whole, i.e. $\{\mathcal{F}(x_i,y_i)\} \sim \{\mathcal{F}(x_i + s,y_i + t)\}$ for any $(s,t)$. So, while $\mathcal{F}$ cannot be stationary under generic translations in the product space, $\mathcal{F}$ may be stationary with respect to translations in the product space $\Omega \times \Omega$ that correspond to a translation in the original factor space, $\Omega$.


Nevertheless, Corollaries \ref{cor 1: kernel sign} and \ref{cor 2: nonstationarity} eliminate most standard covariance models. To construct a valid covariance model we will adopt an approach that is generic and interpretable.

\begin{snugshade}
\begin{definition} \label{def: difference of utility}
    \textbf{[Difference of Utility Model]:} Let $\mathcal{U}:\Omega \times \Omega \rightarrow \mathbb{R}$ denote a random utility function, $\mathcal{U} \sim \text{GP}(\mu_u,k_u)$ for mean utility function $\mu_u: \Omega \times \Omega \rightarrow \mathbb{R}$ and utility covariance function $k_u: \Omega^4 \rightarrow \mathbb{R}$. Then, let the flow function equal the skew-symmetric component of utility:
    \begin{equation}
        \mathcal{F}(x,y) = \mathcal{U}(x,y) - \mathcal{U}(y,x). 
    \end{equation}
\end{definition}
\end{snugshade}

In a game-theory context, a difference of GP utilities models is an example of a random utility model where an agent with traits $x$ receives utility $\mathcal{U}(x,y)$ when interacting with an agent with traits $y$. Then, $\mathcal{F}(x,y)$ is the difference in utility received by $x$ when interacting with $y$, and $y$ when interacting with $x$. Outside of game theory, it is common to construct a flow as some difference in a weight function on a pair of directed edges between two endpoints. For example, work is equal to a difference in log transition rates on pairs of forward and backward edges in biophysical Markov processes \cite{schnakenberg1976network}.

We adopt the difference of utilities model since it places no restrictions on the GP model for $\mathcal{F}$, while it is, itself, unconstrained, so can be parameterized using traditional kernels. In particular:

\begin{snugshade}
\begin{lemma} \label{lemma 1}
\textbf{[All Skew-Symmetric Functions are a Difference in Utilities]} 
A function $\mathcal{F} \sim \text{GP}$ is skew-symmetric if and only if there exists a utility function $\mathcal{U} \sim \text{GP}$ such that $\mathcal{F}$ can be expressed as a difference in utilities. 
\end{lemma}
\end{snugshade}

\textbf{Proof:} Suppose that $\mathcal{U} \sim \text{GP}$ with some mean and covariance. Then, $\mathcal{F}(x,y)$ is also a GP, since linear transformations of Gaussian processes are also Gaussian processes. It is skew-symmetric by construction. Given $\mathcal{F}(x,y)$ it suffices to find a utility function. Let $\mathcal{U}(x,y) = \frac{1}{2} \mathcal{F}(x,y)$. $\blacksquare$ 

The mean and covariance functions for $\mathcal{F}$ follow directly from the mean and covariance functions for $\mathcal{U}$:
\begin{equation} \label{eqn: utility to performance}
    \begin{aligned}
        & \mu_f(x,y) = \mu_u(x,y) - \mu_u(y,x) \\
        & k_f([x,y],[v,w]) = k_u([x,y],[v,w]) - k_u([x,y],[w,v]) - k_u([y,x],[v,w]) + k_u([y,x],[w,v])
    \end{aligned}
\end{equation}

By adopting a difference in utilities model, we can choose any $\mu_u,k_u$ for $\mathcal{U}$ while ensuring that $\mu_f,k_f$ satisfy the constraints of Theorem \ref{thm: skew symmetric GP}. In particular, we can use any standard stationary covariance model for $\mathcal{U}$. In principle, this places no restrictions on the GP model for $\mathcal{F}$, while allowing a more standard parameterization. In practice, it acts as a soft assumption, since we do not test every possible covariance for $\mathcal{U}$. We are willing to accept this tacit restriction in scope as many flows are, in practice, defined as a skew-symmetric difference of some other unconstrained function.

In all that follows, we will also assume that $\mathbb{E}[\mathcal{F}(x,y)] = 0$ for all $x$ and $y$. That is, a priori, we have no reason to expect a particular flow direction between $x$ and $y$. We adopt this assumption since incorporating a nonzero mean only leads to trivial changes in the theory. In particular, we could model $\mathcal{F} = \mu_f + \epsilon$ where $\epsilon$ represents the randomness in flow, and $\mu_f$ is the expected performance. This case is treated in \cite{strang2022network}. Since the transitive/cyclic decomposition is performed by orthogonal projection, adding a nonzero $\mu_f$ only changes the trait-performance result in equation \eqref{eqn: trait performance} by adding the norms of the projections of $\mu_f$ onto the corresponding components. Since this is a trivial update, we will focus on the $\mu_f = 0$ case. This requires the added assumption that the mean utility function is symmetric $\mu_u(x,y) = \mu_u(y,x)$. 

In the following analysis, we will gradually introduce further assumptions, producing ever simpler, more specific, formulas for $\rho$. When we reach a closed form for the SE kernel, we will work outwards to generic scale mixtures. We use this generalization to study the Mat\'ern kernel and the relationship between correlation and regularity. Table \ref{table: Results + Assumptions} outlines the sequence of assumptions and links to the corresponding results.
\begin{table}[] 
\centering
\begin{minipage}[t]{0.45\textwidth}
\centering
\small
\begin{tabular}{ll}
\hline
Number & Assumption \\
\hline
\hline
1 & $\mathcal{X}$ are drawn $\ind$ of the flow       \\
2 & $\mathcal{F} \sim \text{GP}$ invariant under agent label orderings \\ 3 & $\mathcal{F}$ is a difference in GP utilities $\mathcal{U}$ \\ 4 & $\mathcal{U}$ is stationary \\
5 & $\mathcal{U}$ is isotropic\\ 6 & $\mathcal{X}$ are multivariate Gaussian \\  
7 & $k_u$ is squared-exponential (SE) \\
8 & $k_u$ is a mixture of SE kernels 
\\ 
9 & $k_u$ is Mat\'ern 
\\ 
10 & $\mathcal{U}$ is $\nu^{th}$ order M.S.~differentiable 
\end{tabular}
\label{table: Assumption Enumeration}
\end{minipage}
\hspace{0.03\textwidth} 
\begin{minipage}[t]{0.45\textwidth}
\centering
\small
\begin{tabular}{ll||ll}
\hline
Result & Assumptions Used & Result & Assumptions Used\\
\hline
\hline
\ref{result 1} & 1, 2 & \ref{result 8: smoothness} & 1, 2, 3, 4, 6, 7  \\
\ref{result 2} & 1, 2, 3 & \ref{result 8: roughness} & 1, 2, 3, 4, 6, 7\\ 
\ref{result 3} & 1, 2, 3, 4 & \ref{result 9} & 1, 2, 3, 4, 8 \\ 
\ref{result 4} & 1, 2, 3, 4, 5 & \ref{result 10} & 1, 2, 3, 4, 5, 8\\
\ref{result 5} & 1, 2, 3, 4, 6 & \ref{result 10: smoothness} & 1, 2, 3, 4, 5, 8\\ 
\ref{result 6} & 1, 2, 3, 4, 5, 6  & \ref{result 11} & 1, 2, 3, 4, 5, 9\\  
\ref{result 7} & 1, 2, 3, 4, 5, 6, 7  & \ref{result 11: roughness} & 1, 2, 3, 4, 5, 9\\
\ref{result 7: smoothness} & 1, 2, 3, 4, 5, 6, 7  & \ref{result 11: smoothness} & 1, 2, 3, 4, 5, 9\\ 
\ref{result 7: roughness} & 1, 2, 3, 4, 5, 6, 7 & \ref{result 12} & 1, 2, 3, 4, 5, 9 \\ 
\ref{result 8} & 1, 2, 3, 4, 6, 7 & Thm \ref{thm: regularity/correlation in smoothness} & 1, 2, 3, 4, 5, 6, 8, 10\\ 
\end{tabular}
\end{minipage}
\caption{Assumptions sequence and results key. \textbf{Left:} the assumptions introduced. Recall the notation: $\mathcal{X}$ = attributes, $\mathcal{F}$ = flow function, $\mathcal{U}$ = utility function, and $k_u$ = utility kernel. \textbf{Right:} assumptions adopted per result.}
\label{table: Results + Assumptions}
\end{table}

\subsection{Correlation in Gaussian Process Flow Models}

Assume (1) and (2). Then equation \eqref{eqn: sigma and rho for random f independent X} applies, $\mathcal{F}$ is drawn from a GP, and is mean zero. Reducing:
\begin{snugshade}
\begin{result} \label{result 1}
\textbf{[$\sigma^2, \rho$ under (\hyperref[table: Results + Assumptions]{1 - 2})]} 
When the attributes are independent of $\mathcal{F}$, the flow is mean zero, and drawn from a GP:
\begin{equation} \label{eqn: sigma and rho f gp mean zero}
    \begin{aligned}
         & \sigma^2 = \mathbb{E}_{X,Y}[k_{F}([X,Y],[X,Y])] \\
        & \rho 
        = \frac{1}{\sigma^2} \mathbb{E}_{X,Y,W}[k_f([X,Y],[X,W])] = \frac{\mathbb{E}_{X,Y,W}[k_f([X,Y],[X,W])]}{\mathbb{E}_{X,Y}[k_{F}([X,Y],[X,Y])]}.    \end{aligned}
\end{equation}
\end{result}
\end{snugshade}

So, under assumptions (1) and (2), the shared-endpoint correlation is a ratio of expectations of the covariance function $k_f$. Notice that equation \eqref{eqn: sigma and rho f gp mean zero} applies without enforcing skew-symmetry.

To enforce skew symmetry, adopt the difference in utilities model. Then, the covariance function for $\mathcal{F}$ is determined by the covariance function for the utilities $\mathcal{U}$. By equation \eqref{eqn: utility to performance}:
$$
\begin{aligned}
k_f([X,Y],[X,Y]) = (k_{u}([X,Y],[X,Y]) + k_u([Y,X],[Y,X])) - (k_u([X,Y],[Y,X]) + k_u([Y,X],[X,Y])).
\end{aligned}
$$

Since $k$ is a covariance function, it must be symmetric under an exchange of its inputs. Therefore, the last two terms are identical. Since $X$ and $Y$ are independent and identically distributed, they are interchangeable, so $k_{u}([X,Y],[X,Y])$ and $k_u([Y,X],[Y,X])$ are identically distributed random variables. Therefore:
\begin{equation}
    \sigma^2  = 2 \left( \mathbb{E}_{X,Y}[k_u([X,Y],[X,Y])] - \mathbb{E}_{X,Y}[k_u([X,Y],[Y,X])] \right).
\end{equation}


The correlation coefficient follows similarly:
$$
\sigma^2 \rho = \mathbb{E}_{X,Y,W}\left[k_u([X,Y],[X,W]) + k_u([Y,X],[W,X]) \right] - \mathbb{E}_{X,Y,W}\left[ k_u([X,Y],[W,X]) + k_u([Y,X],[X,W])  \right].
$$


Group the terms that share a common orientation. First, since $k_u$ is a covariance function, $k_u([Y,X],[X,W]) = k_u([X,W], [Y,X])$. Then, since $Y$ and $W$ are identical and independent, $k_u([X,W], [Y,X])$ and $k_u([X,Y],[W,X])$ are identically distributed random variables. Consequently, they share the same expectation. Therefore:
\begin{equation}
    \sigma^2 \rho = \mathbb{E}_{X,Y,W}\left[k_u([X,Y],[X,W]) + k_u([Y,X],[W,X]) \right] - 2\mathbb{E}_{X,Y,W}\left[ k_u([X,Y],[W,X])  \right]
\end{equation}

Collapsing the first term requires an additional symmetry assumption. Namely, that $k_u$ cannot distinguish between the first and second entries of its bracketed input. We assume, without loss of generality, that $k_u([x,y],[w,v]) = k_u([y,x],[v,w])$ as, if not, any such difference between the two terms would be averaged away in equation \eqref{eqn: utility to performance}. Then:
\begin{equation}
    \sigma^2 \rho = 2 \left( \mathbb{E}_{X,Y,W}\left[k_u([X,Y],[X,W]) \right] - \mathbb{E}_{X,Y,W}\left[ k_u([X,Y],[W,X]) \right] \right)
\end{equation}

Therefore:
\begin{snugshade}
\begin{result} \label{result 2}
\textbf{[$\sigma^2, \rho$ under (\hyperref[table: Results + Assumptions]{1 - 3})]} 
When the attributes are independent of $\mathcal{F}$, the flow is mean zero, and drawn from a difference in GP utilities:
\begin{equation} \label{eqn: sigma and rho utility difference}
\begin{aligned}
    & \sigma^2 = 2 \left( \mathbb{E}_{X,Y}[k_u([X,Y],[X,Y])] - \mathbb{E}_{X,Y}[k_u([X,Y],[Y,X])] \right) \\
    & \rho = \frac{\mathbb{E}_{X,Y,W}\left[k_u([X,Y],[X,W]) \right] - \mathbb{E}_{X,Y,W}\left[ k_u([X,Y],[W,X]) \right]}{\mathbb{E}_{X,Y}[k_u([X,Y],[X,Y])] - \mathbb{E}_{X,Y}[k_u([X,Y],[Y,X])]}.
\end{aligned}
\end{equation}
\end{result}
\end{snugshade}

Now, suppose that $\mathcal{U}$ is stationary. Then:
\begin{equation}
    k_u([x,y],[v,w]) = \kappa([x - v,y - w])
\end{equation}
for some utility kernel $k$. So, applying equation \eqref{eqn: sigma and rho utility difference}:
\begin{snugshade}
\begin{result} \label{result 3}
\textbf{[$\sigma^2, \rho$ under (\hyperref[table: Results + Assumptions]{1 - 4})]}  When the attributes are independent of $\mathcal{F}$, the flow is mean zero, and drawn from a difference of stationary GP utilities:
\begin{equation} \label{eqn: sigma and rho stationary}
\begin{aligned}
    & \sigma^2 = 2 \left( \mathbb{E}_{X,Y}[\kappa([0,0])] - \mathbb{E}_{X,Y}[\kappa([X-Y,Y-X])] \right) \\
    & \rho = \frac{\mathbb{E}_{X,Y,W}\left[\kappa([0,Y-W]) \right] - \mathbb{E}_{X,Y,W}\left[ \kappa([X-W,Y-X]) \right]}{\mathbb{E}_{X,Y}[\kappa([0,0])] - \mathbb{E}_{X,Y}[\kappa([X-Y,Y-X])]}.
\end{aligned}
\end{equation}
\end{result}
\end{snugshade}

If $k$ is isotropic, then it is symmetric under all unitary transformations of the input, so:
\begin{equation} \label{eqn: isotropic}
k(\tau) = \sigma_u^2 h(\|\tau\|_2)
\end{equation}
for some $h:\mathbb{R}^+ \rightarrow \mathbb{R}$ and $h(0) = 1$. Then, simplifying:
%
%
%
%
\begin{snugshade} 
\begin{result} \label{result 4}
\textbf{[$\sigma^2,\rho$ under (\hyperref[table: Results + Assumptions]{1 - 5})]}  When the attributes are drawn independently of $\mathcal{F}$, the flow is mean zero, and is drawn from a difference of stationary, isotropic, GP utilities:
\begin{equation} \label{eqn: sigma and rho isotropic l2}
\begin{aligned}
    & \sigma^2 = 2 \sigma_u^2 \left( 1 - \mathbb{E}_{X,Y}[h(\sqrt{2}\|X-Y\|_2)] \right) \\
    & \rho = \frac{\mathbb{E}_{X,Y}\left[h(\|X - Y\|_2) \right] - \mathbb{E}_{X,Y,W}\left[ h(\|[X-W,X-Y] \|_2) \right]}{1 - \mathbb{E}_{X,Y}[h(\sqrt{2}\|X-Y\|_2)]}.
\end{aligned}
\end{equation}
\end{result}
\end{snugshade}

Note that, the form of isotropy enforced by equation \eqref{eqn: isotropic} is tighter than necessary, as equation \eqref{eqn: isotropic} enforces isotropy over all unitary transformations of the product space, $\Omega \times \Omega$, not over all unitary transformations of $\Omega$. This is a stricter symmetry than necessary. 
We will relax it when studying anisotropic processes in Section \ref{sec: anisotropic}.

\subsection{Gaussian Traits}

Our following results largely assume Gaussian traits. 
This is, by far, the most restrictive assumption enforced. We adopt a Gaussian model since we are primarily interested in the role of the kernel and since we aim for closed formulas for $\rho$ as a function of the kernel. No such formulas would exist without specifying a distribution for the attributes, and, in absence of a more principled choice, the Gaussian is both plausibly generic and is tractable, allowing closed forms, or nearly closed forms, for the kernels of interest. For sample results that relax the Gaussian traits assumption see Section \ref{sec: mixture models}.

Suppose that $\mathcal{X} = \{X_i\}$ where $X_i \sim \mathcal{N}(\bar{x},\Sigma_x)$. Without loss of generality, we can set $\bar{x} = 0$ since the flow model is stationary on $\Omega$. Then $Z^{(1)} = X - Y \sim \mathcal{N}(0,2 \Sigma_x)$. Similarly, $[Z^{(1)},Z^{(2)}] = [X - Y, X - W]$ is normally distributed with:
\begin{equation} \label{eqn: sigma z}
[Z^{(1)},Z^{(2)}] \sim \mathcal{N}(0,\Sigma_z) \text{ where } 
\Sigma_z = \left[\begin{array}{cc} 2 \Sigma_x & \Sigma_x \\ \Sigma_x & 2 \Sigma_x \end{array} \right] = (I + \textbf{1} \textbf{1}^\intercal) \otimes \Sigma_x
\end{equation}
and where $\textbf{1}$ is the vector of all ones and $\otimes$ is the Kronecker product. Then:

%
\begin{snugshade} \begin{result} \label{result 5}
\textbf{[$\sigma^2,\rho$ under (\hyperref[table: Results + Assumptions]{1 - 4, 6})]}  When the attributes are drawn independently of $\mathcal{F}$, the flow is mean zero, is drawn from a difference of stationary, GP utilities, and the attributes are normally distributed with covariance $\Sigma_x$:
\begin{equation} \label{eqn: sigma and rho isotropic Gaussian}
\begin{aligned}
    & \sigma^2 = 2 \sigma_u^2 \left( 1 - \mathbb{E}[h(\sqrt{2}\|Z^{(1)}\|_2)] \right) \\
    & \rho = \frac{\mathbb{E}\left[h(\|Z^{(1)}\|_2) \right] - \mathbb{E}\left[ h(\|[Z^{(1)}, Z^{(2)}]\|_2) \right]}{1 - \mathbb{E}[h(\sqrt{2}\|Z^{(1)}\|_2)]}
\end{aligned}
\end{equation}
where $[Z^{(1)},Z^{(2)}]$ is normal, mean zero, with covariance $\Sigma_z = (I + \textbf{1} \textbf{1}^{\intercal}) \otimes \Sigma_x$.
\end{result}
\end{snugshade}

We will use equation \eqref{eqn: sigma and rho isotropic Gaussian} to compute $\rho$ when $h$ is a squared exponential or Mat\'ern function, and will generalize the equation to allow anisotropic $k_u$. 

Suppose $X$ are also drawn from an isotropic distribution. Then $\Sigma_x = \sigma^2 I$. Then let $S^{(1)}$ and $S^{(2)}$ be independent, $\chi^2(T)$ random variables. Then $\|Z^{(1)}\|_2 \sim \sigma_x \sqrt{2 S^{(1)}}$. To express $\|[Z^{(1)},Z^{(2)}]\|$, diagonalize $\Sigma_z = (I + \textbf{1} \textbf{1}^{\intercal}) \otimes I_{T \times T}$. The covariance $\Sigma_z$ has eigenvalues $\Lambda(\Sigma_z) = \{3 (\times T), 1 (\times T)\}$ where $\times T$ denotes multiplicity. Since $\| \cdot \|_2$ is unitarily invariant, the specific eigenvectors do not matter. Then, instead of evaluating the expectation over $[Z^{(1)}, Z^{(2)}]$ we can rotate the coordinate system to align with the eigenvectors. In this basis, we need the norm of a random vector with independent Gaussian entries, $T$ of which have standard deviation 3, and $T$ of which have standard deviation 1. Therefore $\|[Z^{(1)},Z^{(2)}]\|^2 \sim S^{(1)} + 3 S^{(2)}$. 

%
\begin{snugshade}
\begin{result} \label{result 6}
\textbf{[$\sigma^2,\rho$ under (\hyperref[table: Results + Assumptions]{1 - 6})]}  When the attributes are drawn independently $\mathcal{F}$, the flow is mean zero, is drawn from a difference of stationary, isotropic, GP utilities, and the attributes are normally distributed with covariance $\Sigma_x = \sigma_x^2 I$:
\begin{equation} \label{eqn: sigma and rho iso iso Gaussian}
\begin{aligned}
    & \sigma^2 = 2 \sigma_u^2 \left( 1 - \mathbb{E}[h(\sigma_x \sqrt{4 S^{(1)} })] \right) \\
    & \rho = \frac{\mathbb{E}\left[h(\sigma_x \sqrt{S^{(1)}}) \right] - \mathbb{E}\left[ h(\sigma_x \sqrt{S^{(1)} + 3 S^{(2)}}) \right]}{1 - \mathbb{E}[h(\sigma_x \sqrt{4 S^{(1)} })]}
\end{aligned}
\end{equation}
where $S^{(1)}, S^{(2)}$ are independent standard $\chi^2(T)$ random variables.
\end{result}
\end{snugshade}

This form is convenient since it eliminates the need to perform integrals over a $T$ dimensional space when $\Omega \subseteq \mathbb{R}^T$. It allows easy numerical evaluation for arbitrarily large $T$. We use equation \eqref{eqn: sigma and rho iso iso Gaussian} for numerical validation, and equation \eqref{eqn: sigma and rho isotropic Gaussian} for analysis, since it generalizes more directly in the anisotropic case, and is simpler to close.

\subsection{Squared Exponential Kernels and Roughness} \label{sec: SE}

We turn now to specific kernels. Suppose that the utility function is drawn with a SE kernel. Then:
\begin{equation}
    h(\tau;\Sigma_u) = \exp\left(- \frac{1}{2} \|\phi \|_{\Sigma_u^{-1}}^2 \right) \text{ where } \| \tau \|_{\Sigma_{u}^{-1}} = \tau^{\intercal} \Sigma_u^{-1} \tau.
\end{equation}

In the isotropic case, $\Sigma_u = l^2 I$ for some length scale parameter $l$. In general, there is no reason to consider generic $\Sigma_u$ and $\Sigma_x$ since all the expectations required are preserved under a coordinate change that either whitens the trait distribution ($\Sigma_x \propto I$) or isotropizes the process ($\Sigma_u \propto I)$. So, we could, without loss of generality, fix $\Sigma_x = I$ or $\Sigma_u = I$. We will, instead, show that all results reduce to functions of the matrix $R = \Sigma_u^{-1} \Sigma_x$ which measures the relative length scales of variation in sampled attributes, and of variation in the process. This is natural since the length scale of the process and variation in the traits only have concrete meaning relative to each other. 

\begin{snugshade}
\begin{definition} \label{def: roughness coefficients}
        \textbf{[Roughness Coefficients]} Any similarity transformation of $R = \Sigma_u^{-1} \Sigma_x$ is a \textit{roughness matrix}. The eigenvalues, $\{r_j\}_{j=1}^T = \{ \lambda_j(R)^{1/2} \}_{j=1}^T$ are the \textit{roughness coefficients}. When both process and distribution are isotropic $r_j = r = \sigma_x/l$, where $l$ is the characteristic length scale for variation in the process. 
\end{definition}
\end{snugshade}

The roughness coefficients are dimensionless. They are the principal standard deviations in the attributes in the coordinate system that isotropizes the process. 
Changes in $\Sigma_x$ and $\Sigma_u$ that preserve the roughness coefficients have no impact on $\sigma^2$ and $\rho$. This reduction reflects the duality of similarity and smoothness. The apparent smoothness of a process depends on the units used to measure length. Changing $\Sigma_x$ amounts to changing those units.

The shared-endpoint correlation $\rho$, is completely determined by the roughness coefficients and shape parameters of the kernel. In the squared exponential case, $\rho$ depends on the roughness coefficients alone. Accordingly, we use the SE kernel to isolate and illustrate the relationship between roughness, as measured by $\{r_j\}_{j=1}^T$, and shared-endpoint correlation $\rho$.

\subsubsection{Isotropic Models}

Since the variance and shared-endpoint correlation in flow depend only on the roughness $R$, there is no reason to consider mixed isotropy cases where the process is isotropic, and the distribution is not. Instead, we should require either isotropy in the process given a coordinate system that whitens the trait distribution, or, isotropy in the trait distribution in a coordinate system that isotropizes the process. Both cases reduce to enforcing $R = r^2 I$ for some roughness coefficient $r$. 

\begin{snugshade}
\begin{result} \label{result 7}
    \textbf{[$\rho$ under (\hyperref[table: Results + Assumptions]{1 - 7})]}  When the attributes are drawn independently of $\mathcal{F}$, the flow is mean zero, is drawn from a difference of stationary, isotropic, GP utilities with SE kernel with length scales $\Sigma_u$, and the attributes are normally distributed with covariance $\Sigma_x = \sigma_x^2 I$:
    \begin{equation} \label{eqn: rho isotropic SE gaussian traits}
\rho(r) = \frac{\left(1 + 2 r^2 \right)^{-T/2} - (1 + r^2)^{-T/2} (1 + 3 r^2)^{-T/2}}{1 - (1 + 4 r^2)^{-T/2}}
\end{equation}
provided $R = \Sigma_u^{-1} \Sigma_x = r^2 I$ for some scalar $r$. Then $\{r_j\}_j = r$ for all $j$, so $r$ is the roughness coefficient.   
\end{result}
\end{snugshade}

\textbf{Calculation Outline:}  Equation \eqref{eqn: rho isotropic SE gaussian traits} follows from the observation that the product of two squared-exponential functions is a squared-exponential function, and, all squared-exponential functions can be integrated in closed form by finding the normalizing factor for the matching Gaussian density. All of the expectations required are products of squared-exponentials, so can be recovered by completing the square, then normalizing. For details, see Appendix Section \ref{app: rho SE}. $\blacksquare$

\begin{figure}[t]
    \centering
    \includegraphics[width = \textwidth]{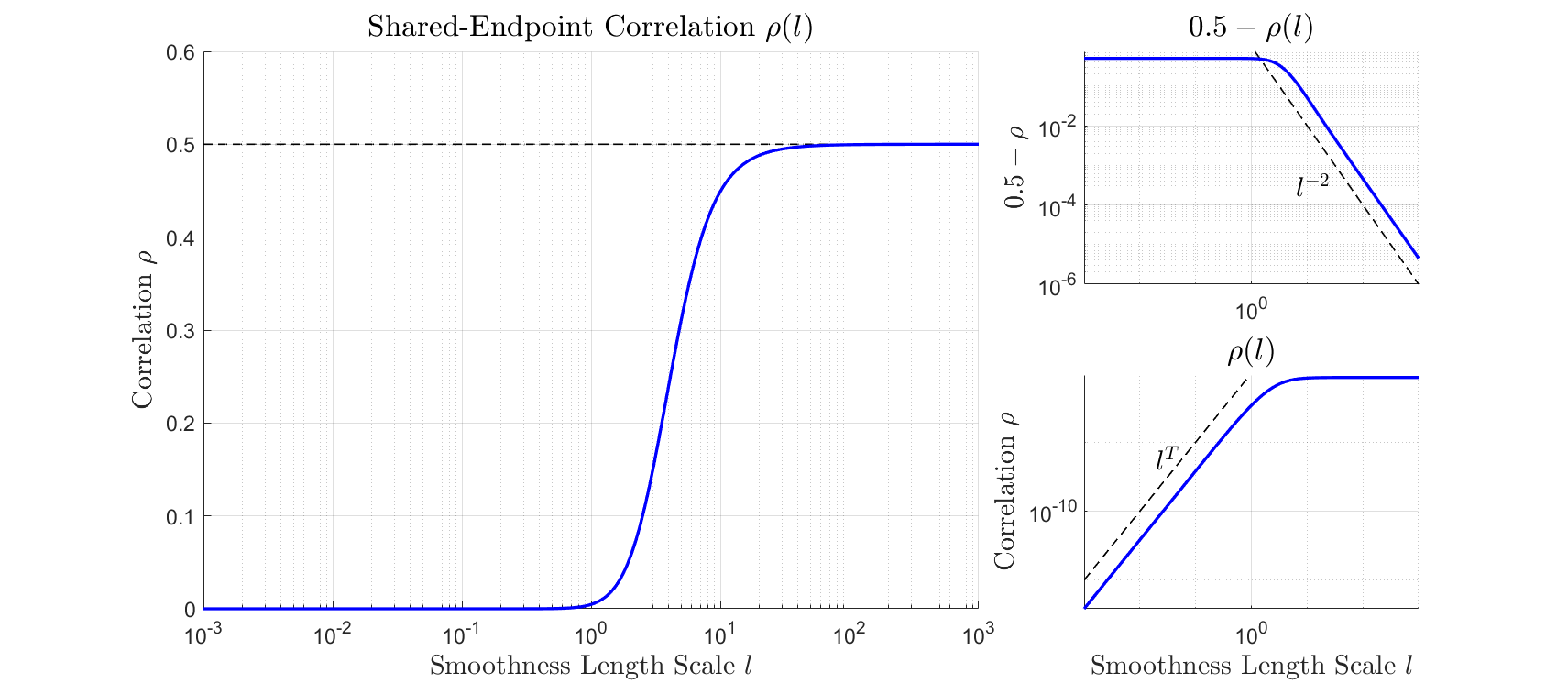}
    \caption{Shared-endpoint correlation $\rho$ as a function of the smoothness length scale $l \propto r^{-1}$ for $T = 5$, $\Sigma_x$ drawn randomly as the product of two Gaussian random matrices. The left-most plot illustrates the sigmoid behavior of $\rho$ as a function of $l$, and that $\rho$ converges to $0.5$ for sufficiently smooth functions. The right panels illustrate the predicted convergence rate to $0.5$, $\mathcal{O}(l^{-2}) = \mathcal{O}(r^{2})$, and to 0, $\mathcal{O}(l^T) = \mathcal{O}(r^{-T})$. Both axes in the right column are log scaled.}
    \label{fig: Correlation isotropic}
\end{figure}

Equation \eqref{eqn: rho isotropic SE gaussian traits} defines a monotonically decreasing sigmoid function in $r$. It approaches $1/2$ as $r$ approaches zero, and approaches 0 as $r$ approaches infinity. The sigmoid is illustrated in Figure \ref{fig: Correlation isotropic}. Since the limit as $r$ approaches zero corresponds to a concentration limit ($\sigma_x \rightarrow 0$) when the kernel is fixed, and draws from a SE process are almost surely infinitely differentiable \cite{williams2006gaussian}, the fact that $\rho(r)$ converges to $1/2$ as $r \rightarrow 0$ is consistent with the analysis presented in Fact \ref{fact 1} and in \cite{cebra2023similarity}. Unlike the local analysis developed in \cite{cebra2023similarity} equation \eqref{eqn: rho isotropic SE gaussian traits} is global, exact, and non-asymptotic. It is less general, but allows detailed analysis of the asymptotics for small $r$  beyond convergence rates alone. It also allows asymptotic analysis for large $r$, and characterizes the transition region, where the correlation is neither near $1/2$ nor 0. 

The asymptotic behavior of $\rho(r)$ is summarized in Results \ref{result 7: smoothness} and \ref{result 7: roughness}.

\begin{snugshade}
\begin{limit_result} \label{result 7: smoothness}
\textbf{[Smoothness Limit of $\rho$ under (\hyperref[table: Results + Assumptions]{1 - 7})]}  Under the same assumptions as stated in \ref{result 7}, in the limit as $r \rightarrow 0$:
\begin{equation} \label{eqn: isotropic SE smoothness}
    \rho(r) \simeq \frac{1}{2} \left[1 - \frac{(T-1)}{2} r^2 + \mathcal{O}(r^4) \right] \xrightarrow{r \rightarrow 0} \frac{1}{2}.
\end{equation}
\end{limit_result}
\end{snugshade}

\begin{snugshade}
\begin{limit_result} \label{result 7: roughness}
\textbf{[Roughness Limit of $\rho$ under (\hyperref[table: Results + Assumptions]{1 - 7})]}  Under the same assumptions as stated in \ref{result 7}, in the limit as $r \rightarrow \infty$:
\begin{equation}  \label{eqn: rho roughness limit isotropic SE}
    \rho(r) \simeq (\sqrt{2} r)^{-T} + \mathcal{O}(r^{-2T}) \xrightarrow{r \rightarrow \infty} 0.
\end{equation}

\end{limit_result}
\end{snugshade}

For detailed calculations see Appendix Section \ref{app: smoothness and roughness iso SE}. 

Result \ref{result 7: smoothness} is consistent with the analysis developed in \cite{cebra2023similarity}. In particular, $\rho(r)$ converges to $1/2$ at order $\mathcal{O}(r^2) = \mathcal{O}(\sigma_x^2)$. Note the dependency on dimension. In particular, if $T = 1$, then $\rho(r) = 1/2$ to order $r^4$ as predicted in Fact \ref{fact 1} and in \cite{cebra2023similarity}. 

Result \ref{result 7: roughness} contrasts Result \ref{result 7: smoothness}. Notice that the convergence rate to $1/2$ as $r$ vanishes is, to order of convergence, independent of the dimension $T$ if $T > 1$. In contrast, the rate at which $\rho(r)$ approaches 0 as $r$ diverges depends on the order $T$. In both cases, the larger $T$, the faster $r$ decreases, however, the convergence to $\rho = 0$ is much more sensitive to the dimension $T$, and, for large $T$, occurs much faster than $\rho(r)$ diverges from $1/2$. We will see that this result is a special case of the more general asymptotic behavior when using an anisotropic model. Deviations from the smoothness limit ($r = 0$) are controlled by the sum of the roughness coefficients. Deviations from the roughness limit ($r \rightarrow \infty$) are controlled by the product of the roughness coefficients. 

The asymptotic behavior of $\rho(r)$ can be used to build simpler expressions that approximate $\rho(r)$. The simplest two-point Pad\'e approximation \cite{brezinski1996extrapolation} consistent with both limits is the rational function:
\begin{equation}
    \rho \simeq \frac{1}{2} \left(1 + \frac{(T-1)}{2} r^2 + \frac{1}{2} (\sqrt{2} r)^{T} \right)^{-1}
\end{equation}
where the approximation has exact asymptotics for both $r$ near zero and near infinity. We will use approximations of this form to approximate $\rho$ for more complicated kernels (Mat\'ern) when we cannot close the necessary expectations exactly, but can compute asymptotics. The accuracy of the approximation is verified in Figure \ref{fig: nu validation}. While the approximation is not exact for intermediate $r$, matching the asymptotic behavior near $r = 0$ and $r \rightarrow \infty$ produces approximations that extrapolate well to intermediate $r$.

\subsubsection{Anisotropic Models} \label{sec: anisotropic}

Suppose, now, that $R$ is not proportional to an identity matrix. Then (see Appendix Section \ref{app: smoothness and roughness aniso SE}):
\begin{snugshade}
\begin{result} \label{result 8}
    \textbf{[$\rho$ under (\hyperref[table: Results + Assumptions]{1 - 4, 6 - 7})]}  When the attributes are drawn independently of $\mathcal{F}$, the flow is mean zero, is drawn from a difference of stationary, GP utilities with SE kernel with length scales $\Sigma_u$, and the attributes are normally distributed with covariance $\Sigma_x$:
    \begin{equation} \label{eqn: rho anisotropic}
    \rho(R) = \frac{\det\left(I + 2R \right)^{-1/2} - \det\left(I + R  \right)^{-1/2} \det\left(I + 3 R  \right)^{-1/2}}{1 - \det\left(I + 4 R \right)^{-1/2}}
    \end{equation}
    where $R$ is any matrix similar to $\Sigma_u^{-1} \Sigma_x$. If $r_j = \sqrt{\lambda_j(R)}$ are the roughness values, then:
    \begin{equation} \label{eqn: rho anisotropic roughness values}
        \rho(r) = \frac{\prod_{j=1}^T (1 + 2 r_j^2)^{-1/2} - \prod_{j=1}^T (1 +  r_j^2)^{-1/2} (1 + 3 r_j^2)^{-1/2}}{1 - \prod_{j=1}^T (1 + 4 r_j^2)^{-1/2}}.
    \end{equation}
\end{result}
\end{snugshade}

Notice that, as in the isotropic case, the actual form for $\rho$ depends only on the list of roughness values, so can be expressed as a function of roughness rather than the original matrices $\Sigma_u$ and $\Sigma_x$. We will use this observation in Section \ref{sec: mixture models} to derive a general form for $\rho$ given a kernel and trait distribution that can be expressed as Gaussian mixtures, and thus, a joint distribution over the roughness values.

Like equation \eqref{eqn: rho isotropic SE gaussian traits}, equation \eqref{eqn: rho anisotropic roughness values} defines a sigmoid function along each roughness value, and, along any ray pointing from $r = 0$ into $\mathbb{R}^{T+}$. As before, when $r$ approaches zero, $\rho(r)$ approaches $1/2$, and as $\|r\|$ approaches infinity $\rho$ approaches zero. The exact asymptotics follow.

First, all the terms in equation \eqref{eqn: rho anisotropic roughness values} are soft-thresholded versions of $k r^2$ for some $k \in {1,2,3,4}$, bounded below by:
\begin{equation}
    (1 + k r^2) \geq \min\{k r^2, 1\}.
\end{equation}

When $r_j \ll 1$, the corresponding terms in equation \eqref{eqn: rho anisotropic roughness values} are all approximately 1. When $r_j^2 \gg 1$, then all terms are approximately $k r^2$. This suggests an initial approximation. If the roughness values span many orders of magnitude, and $T$ is small relative to the orders of magnitude spanned, then the range of possible values between $\max_j\{r_j\}$ and $\min_j\{r_j\}$ may be sparsely sampled, so there may not be any roughness values such that $k r_j^2$ is on the order of 1 for $k = 1,2,3,4$. Then, we may be able to approximate the shared-competitor correlation by separating the roughness values into a rough set and smooth set:
\begin{equation}
\begin{aligned}
    & \textbf{Rough Set: } \mathcal{R} = \{j \mid r_j^2 \gg 1\}\\
    & \textbf{Smooth Set: } \mathcal{S} = \{j \mid r_j^2 \ll \tfrac{1}{4}\}\\
\end{aligned}
\end{equation}

If $\mathcal{R} \cup \mathcal{S}$ covers all $j \in [1,T]$, then the set of roughness values can be partitioned into the rough and smooth sets. Let $ |\mathcal{S}|$ and $|\mathcal{R}|$ denote the cardinality of these sets. Then $|\mathcal{S}| + |\mathcal{R}| = T$. 

Now, assuming such a partition exists, it is reasonable to adopt the approximation:
\begin{equation}
    (1 + k r_j^2) \simeq \begin{cases} & 1 \text{ if } j \in \mathcal{S} \\
    & k r_j^2 \text{ if } j \in \mathcal{R}
    \end{cases}
\end{equation}

Under this approximation, the shared-competitor correlation simplifies dramatically since the smooth set drops out of the calculation, and the powers can be distributed.

The smoothness limit follows by direct expansion of equation \eqref{eqn: rho anisotropic roughness values} (see Appendix Section \ref{app: smoothness and roughness aniso SE}).  Suppose that all the roughness values are small. Then $\mathcal{R} = \emptyset$, and, in the limit as the largest roughness value vanishes:

\begin{snugshade}
\begin{limit_result} \label{result 8: smoothness}
    \textbf{[Smoothness Limit of $\rho$, under (\hyperref[table: Results + Assumptions]{1 - 4, 6 - 7})]}  Under the same assumptions as stated in \ref{result 8}, in the limit as $r \rightarrow 0$:
    \begin{equation} \label{eqn: half less rho smooth}
        \frac{1}{2} - \rho(R) \simeq \frac{1}{4} \left( \text{tr}(R) - \frac{\text{tr}(R^2)}{\text{tr}(R)}  \right) + \mathcal{O}(R^2) = \frac{1}{4}((T - 1) \text{av}(r^2) - \text{dis}(r^2)) + \mathcal{O}(r^4)
    \end{equation}
    where $T$ is the dimension of the trait space, $\text{av}$ denotes an arithmetic average and $\text{dis}$ denotes the index of dispersion (ratio of variance to average).
\end{limit_result}
\end{snugshade}

If $T = 1$ then there is no variance in the roughness values, so the dispersion is zero. Then, $\rho(R) = \frac{1}{2} - \mathcal{O}(r^4)$. This result mirrors analogous conclusions in \cite{cebra2023similarity}, where it was shown that cyclicity is, to lowest order, associated with an off-diagonal block of the Hessian of $\mathcal{F}$, and this block vanishes if $T = 1$. Alternately, when $T$ is large, cyclicity is essentially proportional to the average roughness value, and for fixed average roughness, grows proportional to $T$. Dispersion in the roughness values only weakly reduces this tendency since the dispersion-dependent term does not scale with the trait space dimension.

Unlike the smoothness limit, which only applies if $\mathcal{R} = \emptyset$, the roughness limit does not require that $\mathcal{S} = \emptyset$. We will show that the roughness limit requires only one large $r_j$. Let $\bar{r}_{\mathcal{R}}$ denote the geometric average of the roughness values in the rough set. Then, provided we can partition the roughness coefficients into a rough and smooth set:

\begin{snugshade}
\begin{limit_result} \label{result 8: roughness}
    \textbf{[Roughness Limit of $\rho$, under (\hyperref[table: Results + Assumptions]{1 - 4, 6 - 7})]}  Under the same assumptions as stated in \ref{result 8}, in the limit as $r \rightarrow \infty$:
    \begin{equation} \label{eqn: rho rough smooth approx}
    \rho(R) \approx 2^{-|\mathcal{R}|/2} \bar{r}_{\mathcal{R}}^{-|\mathcal{R}|} \left[1 - \left((3/2)^{-|\mathcal{R}|/2} - 4^{-|\mathcal{R}|/2} \right) \bar{r}_{\mathcal{R}}^{-|\mathcal{R}|}\right] = \mathcal{O}(\bar{r}_{\mathcal{R}}^{-|\mathcal{R}|})\xrightarrow{\bar{r}_{\mathcal{R}} \rightarrow \infty} 0.
\end{equation}
\end{limit_result}
\end{snugshade}

Result \ref{result 8: roughness} generalizes Result \ref{result 7: roughness}, which assumes that $|\mathcal{R}| = T$, and $\bar{r}_{\mathcal{R}} = r$ since $r_j = r$ for all $j$. Result \ref{result 8: roughness} implies that roughness in any direction is sufficient to decorrelate the flow on edges sharing an endpoint. If any non-empty subset of $|\mathcal{R}|$ roughness values are diverging, and $\bar{r}_{\mathcal{R}}$ is the geometric average of the rough values, then the shared-competitor correlation converges to zero at rate $\mathcal{O}(\bar{r}_{\mathcal{R}}^{-|\mathcal{R}|})$. 

\begin{wrapfigure}[30]{r}{0.45\textwidth}
    \centering
    \includegraphics[trim = 15 8 12 8,clip,width = 0.4\textwidth]{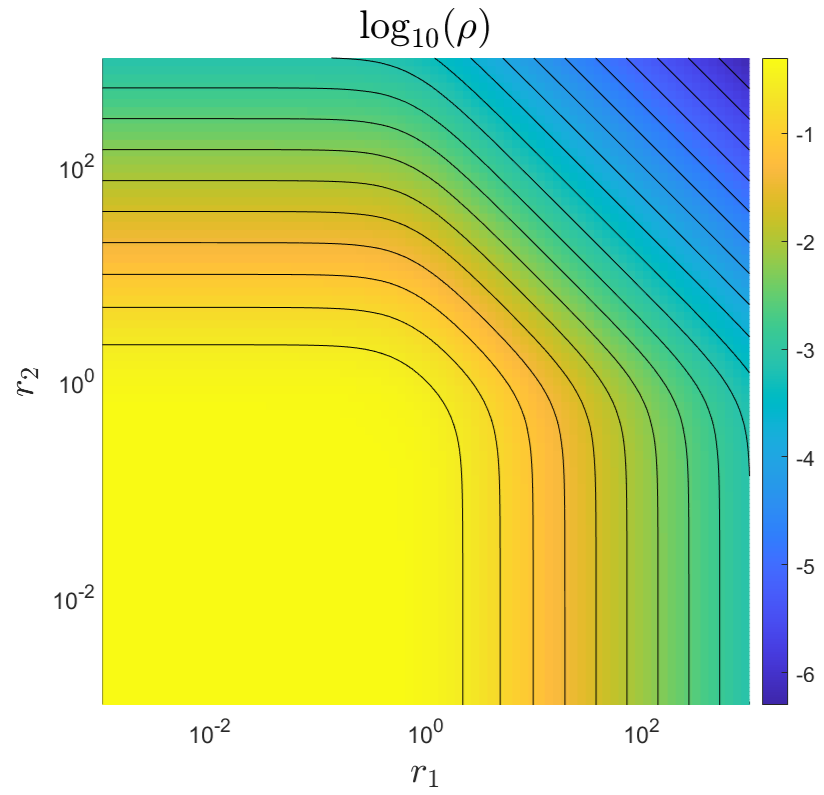}
    \caption{Shared-endpoint correlation $\rho$ as a function of roughness coefficients $r_1$ and $r_2$. Yellow (small $r$) corresponds to correlations near $1/2$. Blue (large $r$) correspond to correlations near zero. Only one large roughness value is needed to trigger a transition from order (high correlation) to disorder (low correlation). The change from convex to log-linear level curves for large $r$ reflects the shift in convergence rate from a rate that depends on the arithmetic average of the roughness values in the smooth limit, to their geometric average in the rough limit. The sigmoidal curve in Figure \ref{fig: Correlation isotropic} corresponds to diagonal where $r_2 = r_1 = l^{-1}$.}
    \label{fig: Correlation anisotropic}
\end{wrapfigure}

Figure \ref{fig: Correlation anisotropic} shows the correlation coefficient $\rho$ as a function of the first two roughness values for $T = 2$, and illustrates the predicted convergence behavior to $1/2$ when smooth, and $0$ when any roughness value is large.


Equation \eqref{eqn: rho anisotropic roughness values} provides the generic relationship between roughness and the shared-endpoint correlation when the kernel is squared exponential, and the attribute distribution is Gaussian. In the next section we will show how to generalize this result to kernels and attribute distributions that can be expressed as scale-mixtures of Gaussian distributions.

\subsection{Mixture Models} \label{sec: mixture models}


To generalize our results, we turn now to scale mixtures. A scale mixture distribution is formed by averaging a parametric distribution over a range of values of its scale parameter. In particular, a scale mixture of Gaussian distributions is the distribution for a random variable $X$ that is conditionally Gaussian, conditioned on the covariance $V = v$, i.e.~$X|V = v \sim \mathcal{N}(0,v)$, where $V$ is itself drawn from a distribution of possible covariances \cite{andrews1974scale}. The distribution over the scale parameter, $V$, is the mixing distribution.  

Scale mixtures of Gaussian distributions are a reasonably flexible family of even, unimodal distributions widely employed in hierarchical models for efficient, robust inference \cite{abanto2010robust,calvetti2019hierachical,calvetti2020sparse,gneiting1997normal}. If $X$ is a univariate random variable, with density $f_X$ that is symmetric about zero, then $X$ may be expressed as a scale mixture of normals if and only if $f_X(\sqrt{y})$ is completely monotone; $(-\frac{d}{dy})^k f_X(\sqrt{y}) \geq 0$ for all $k \in \mathbb{Z}$ and all $y \geq 0$ \cite{andrews1974scale}. All logistic distributions, Laplace distributions, and Student's T distributions are scale mixtures of normals. Laplace distributions correspond to exponential mixing distributions. Student's T distributions correspond to inverse gamma mixing distributions \cite{andrews1974scale}. Thus, by adopting a mixture model, we can easily express symmetric unimodal distributions with dramatically different tail behaviors (Gaussian, sub-Gaussian, exponential, or, power-law) \cite{wainwright2019high}.

Suppose that both $\pi_X$ and the kernel function, $\kappa_u$ can be expressed as scale-mixtures of normal densities. Note that, the mixture used to express $\kappa_u$ need not be normalized. However, since $\rho$ is invariant to the overall variance, changing the normalization of $\kappa_u$ has no impact on $\rho$. Also note that we are free to use a joint mixture model, where $\Sigma_x$ and $\Sigma_u$ are coupled. In either case, $\Sigma_x$ and $\Sigma_u$ are both random matrices, drawn from either a pair of mixing distributions, or from a joint mixing distribution. Then, reversing the order of expectation inside each term defining $\rho$ gives:
\begin{equation} \label{eqn: mixture unsimplified}
    \rho = \frac{\mathbb{E}_{\Sigma_x, \Sigma_u}\left[\mathbb{E}_{Z^{(1)} \sim \mathcal{N}(0,\Sigma_x)}\left[h_{\text{SE}}(\|Z^{(1)}\|_2;\Sigma_u) \right] \right] - \mathbb{E}_{\Sigma_x, \Sigma_u}\left[\mathbb{E}_{[Z^{(1)},Z^{(2)}] \sim \mathcal{N}(0,(I + \textbf{1} \textbf{1}^{\intercal}) \otimes \Sigma_x)}\left[ h_{\text{SE}}(\|[Z^{(1)}, Z^{(2)}]\|_2;\Sigma_u) \right] \right]}{1 - \mathbb{E}_{\Sigma_x, \Sigma_u}\left[\mathbb{E}_{Z^{(1)} \sim \mathcal{N}(0,\Sigma_x)}[h_{\text{SE}}(\sqrt{2}\|Z^{(1)}\|_2;\Sigma_u)] \right]}.
\end{equation}

The inner expectations in equation \eqref{eqn: mixture unsimplified} are each of the form simplified by equation \eqref{eqn: rho anisotropic}. Therefore:

\begin{snugshade}
\begin{result} \label{result 9}
\textbf{[$\rho$ under (\hyperref[table: Results + Assumptions]{1 - 4, 8})]}  When the attributes are drawn independently of $\mathcal{F}$, the flow is mean zero, and is drawn from a difference of stationary, GP utilities with a mixture of SE kernels:
\begin{equation} \label{eqn: mixture simplified}
    \rho = \frac{\mathbb{E}_{r \sim \pi_r}\left[ \prod_{j=1}^{T} (1 + 2 r_j^2)^{-1/2} \right] - \mathbb{E}_{r \sim \pi_r}\left[ \prod_{j=1}^{T} (1 + r_j^2)^{-1/2} (1 + 3 r_j^2)^{-1/2}  \right]}{1 - \mathbb{E}_{r \sim \pi_r}\left[ (1 + 4 r_j^2)^{-1/2}  \right]}.
\end{equation}
where $\pi_r$ is the marginal distribution for the roughness coefficients, $r_j = \sqrt{\lambda_j(R)}$ where $R = \Sigma_u^{-1} \Sigma_x$ (or any matrix similar to $\Sigma_u^{-1} \Sigma_x$). 
\end{result}
\end{snugshade}

\begin{snugshade}
\begin{result} \label{result 10}
\textbf{[$\rho$ under (\hyperref[table: Results + Assumptions]{1 - 5, 8})]}  When the attributes are drawn independently of $\mathcal{F}$, the flow is mean zero, and is drawn from a difference of stationary, isotropic, GP utilities with a mixture of SE kernels:
\begin{equation} \label{eqn: isotropic mixture simplified}
    \rho = \frac{\mathbb{E}_{r \sim \pi_r}\left[ (1 + 2 r^2)^{-T/2} \right] - \mathbb{E}_{r \sim \pi_r}\left[ (1 + r^2)^{-T/2} (1 + 3 r^2)^{-T/2}  \right]}{1 - \mathbb{E}_{r \sim \pi_r}\left[ (1 + 4 r^2)^{-T/2}  \right]}.
\end{equation}
\end{result}
\end{snugshade}

Equation \eqref{eqn: mixture simplified} and \eqref{eqn: isotropic mixture simplified} provide compelling general forms. First, all of the expectations run over the same distribution; the distribution of the roughness coefficients. When isotropic, these expectations reduce to one-dimensional integrals, no matter the initial dimension $T$. This fact is numerically convenient since the necessary expectations can be performed with one-dimensional quadrature. We will see that this special case is convenient when calculating $\rho$ for standard, non-SE, kernel models. More strongly, adopting an isotropic mixture extends our SE result to all isotropic kernels \cite{williams2006gaussian}:
\begin{snugshade}
\begin{fact} \label{fact 3}
    \textbf{[Mixtures Model Isotropic Kernels] \cite{stein1999interpolation}} All isotropic kernel functions $\kappa_u$ in $\mathbb{R}^T$ may be represented as a scale-mixture of SE kernels.
\end{fact}
\end{snugshade}

Second, the expectations can be reduced, generically, to expectations over the $T$ roughness-coefficients, rather than the $T(T+1)$ free entries of $\Sigma_x$ and $\Sigma_u$. This reduces the number of different models an analyst need consider, since any two joint distributions over $\Sigma_x$ and $\Sigma_u$ that produce the same marginal distribution over the roughness coefficients will share the same correlation $\rho$.


To illustrate the isotropic case, suppose that $\Sigma_x$ and $\Sigma_u$ are drawn independently, each proportional to the identity with $\Sigma_x = \sigma_x^2 I$, $\Sigma_u = \sigma_u^2 I$, for random variances $\sigma_x^2$ and $\sigma_u^2$. Then, Table \ref{table: mixture key} provides a reference key that relates trait distribution and kernel, to the choice of mixture distribution for the variances, and the ensuing distribution of the roughness coefficients squared $(\sigma_x/\sigma_u)^2$:

\begin{table}[] 
\centering
\begin{tabular}{ll|ll|l}
\hline
              Trait Dist: $\pi_X$  & Utility Kernel: $\kappa_u$     &  Trait Mix: $\sigma_x^2$ & Kernel Mix: $\sigma_u^2$ & Roughness: $r^2 = (\sigma_x/\sigma_u)^2$ \\ \hline \hline
Laplace & Squared Exp. & Exponential & $-$ & Exponential          \\
Gaussian & Rational Quad. & $-$ & $\text{Inv. Gamma}(1,\cdot)$ & Exponential \\ \hline
Gaussian & Rational Quad. & $-$ & Inv. Gamma & Gamma \\ \hline
Student's T & Squared Exp. & Inv. Gamma & $-$ & Inv. Gamma \\
Gaussian & Exponential & $-$ & Exponential & $\text{Inv. Gamma}(1,\cdot)$        \\ \hline
Laplace & Exponential & Exponential & Exponential & Half-Cauchy \\ \hline 
Student's T & Rational Quad. & Inv. Gamma & Inv. Gamma & Beta Prime \\ \hline

\end{tabular}
\caption{Mixture key. The first two columns provide the family for $\pi_X$ and $\kappa_u$. The second pair of columns provide corresponding mixture distributions for $\sigma_x^2$ and $\sigma_u^2$. The final column provides the corresponding distribution for $r^2 = (\sigma_x/\sigma_u)^2$. Rows of the table are organized according to the distribution over $r^2$, since all models with the same distribution over $r^2$ share the same correlation $\rho$. Empty cells correspond to no mixture, i.e.~mixture against a delta distribution.  } 
\label{table: mixture key}
\end{table}

As an example, consider the case when $r^2$ is drawn from an exponential distribution. This occurs when the trait variance is drawn from an exponential distribution and the kernel variance is a fixed constant. Equivalently, when the trait-distribution is a Laplace distribution, and the kernel is a squared exponential kernel \cite{andrews1974scale}. The roughness squared is also an exponential random variable when the trait variance is fixed, and the kernel variance is drawn from an inverse Gamma distribution with shape parameter $\alpha = 1$, i.e.~when $\sigma_u^{-2}$ is exponentially distributed. The corresponding trait-distribution is normal, and kernel is rational quadratic with shape parameter $\alpha = 1$ (see \cite{williams2006gaussian}). 

As in the squared exponential case, the correlation coefficient $\rho$ admits an interpretable expansion in the smoothness limit. The limit in the mixture case can be expressed as a correction to the limit in the SE case (see equation \eqref{eqn: isotropic SE smoothness}), with roughness set to the standard deviation in $r \sim \pi_r$. In particular, mixing decreases the correlation by a term proportional to the dispersion in roughness values over the mixture distribution. For example, 

\begin{snugshade}
\begin{limit_result} \label{result 10: smoothness}
    \textbf{[Smoothness Limit of $\rho$ under (\hyperref[table: Results + Assumptions]{1 - 5, 8})]}  Under the same assumptions as stated in \ref{result 10}, in the limit as $r \rightarrow 0$:
    \begin{equation} \label{eqn: smoothness limit of iso mixture}
        \frac{1}{2} - \rho \approx \frac{(T-1)}{4} \frac{\mathbb{E}_{r \sim \pi_r}\left[ r^4 \right]}{ 
 \mathbb{E}_{r \sim \pi_r}\left[r^2 \right]} =  \frac{(T-1)}{4}\left(\hat{r}^2 + \frac{ \mathbb{V}_{r \sim \pi_r}[r^2]}{\hat{r}^2 } \right) 
    \end{equation}
    where $\hat{r}^2 = \mathbb{E}_{r \sim \pi_r}[r^2]$, and $\mathbb{V}$ denotes variance.
\end{limit_result}
\end{snugshade} 

All details, including the anisotropic case, are provided in Appendix Section \ref{app:  smoothness limit mixture}.


In Section \ref{sec: Matern} we will show that equation \eqref{eqn: isotropic mixture simplified} also applies to scale-mixture representations of the power-spectral density of a stationary isotropic process. We will use this result to calculate the asymptotic behavior of $\rho$ when $\kappa_u$ is Mat\'ern. This corresponds to evaluating the expectations in \eqref{eqn: isotropic mixture simplified} with respect to an inverse gamma mixture. In Section \ref{sec: generalization regularity} we use the mixture approach in frequency space to show that all stationary, isotropic processes with the same degree of differentiability share the same convergence rates in the smoothness limit. 

We turn, now, to an important example of a mixture model,  Mat\'ern processes.

\subsection{Mat\'ern Kernels and Regularity} \label{sec: Matern}

We introduce the Mat\'ern process to study the influence of regularity. The Mat\'ern kernel includes an additional shape parameter, $\nu$, which controls the regularity of sample draws from the corresponding process (see Section \ref{sec: GPs}). We aim to isolate its impact on the shared-endpoint correlation. For simplicity, we restrict our attention to the case when the kernel is isotropic, so $\sigma^2$ and $\rho$ take the forms shown in equation \eqref{eqn: sigma and rho isotropic l2}.

We will also use the Mat\'ern kernel to motivate a last transformation of the integrals defining $\sigma^2$ and $\rho$. All expectations are an inner product between a probability density and a kernel. These inner products can be expressed in different function bases. For example, these can be converted into inner products with respect to a frequency variable via a Fourier transform. This replaces the kernel with its power spectral density (PSD), and the distribution with its characteristic function. This change of basis is useful if the the new inner product is simpler to compute than the old one.

The Mat\'ern kernel is simpler in frequency space. So, replace the Mat\'ern kernel with its PSD and the Gaussian density with its characteristic function. Then, recognizing that the PSD of the Mat\'ern kernel is a scale mixture of Gaussians allows simplification via the mixture approach demonstrated in equation \eqref{eqn: mixture simplified}. 

First, use Parseval's identity to represent an expectation with an inner-product in frequency space:
\begin{equation} \label{eqn: expectation in frequency space}
    \mathbb{E}_{Z \sim f_Z}\left[ g(Z) \right] = \langle \hat{g}, \hat{f}_Z
 \rangle = \int_{\zeta \in \mathbb{R}^T} \hat{g}(\zeta) \hat{f}_Z(\zeta) d \zeta.
 \end{equation}
Recall that the characteristic function of a Gaussian distribution is itself Gaussian in form, with inverted covariance:
\begin{equation}
    \hat{f}_{Z \sim \mathcal{N}(0,\Sigma_z)}(\zeta) = \exp\left( - \frac{(2 \pi)^2}{2} \zeta^{\intercal} \Sigma_z \zeta \right) = \exp\left( - \frac{1}{2} \|\zeta \|_{4 \pi^2 \Sigma_z}^2 \right).
\end{equation}
Then, a Gaussian expectation can be written in frequency space as an expectation against a different, scaled Gaussian:
\begin{equation} \label{eqn: gaussian expectation in frequency space}
    \mathbb{E}_{Z \sim \mathcal{N}(0,\Sigma_z)}[g(Z)] = \frac{1}{(2 \pi)^{T/2}| \text{det}(\Sigma_z)|^{1/2}} \mathbb{E}_{\zeta \sim \mathcal{N}(0, (2 \pi)^{-2} \Sigma_z^{-1})}[\hat{g}(\zeta)].
\end{equation}

We will also exploit the Mat\'ern kernel's \eqref{eqn: Matern} simple power spectral density to compute the shared-endpoint correlation. The form of the PSD depends on the dimension of the ambient space, $d$, kernel length scale, $l$, and regularity, $\nu$:
    \begin{equation} \label{eqn: Matern kernel in frequency space}
        \kappa_{d}(\tau; \nu, l) \leftrightarrow \hat{\kappa}_{d}(\zeta; \nu, l) = \frac{\Gamma(\nu+\tfrac{d}{2})}{\Gamma(\nu)} 2^d \pi^{d/2} \left(\frac{l^2}{2\nu}\right)^{d/2} \left(1 + \frac{l^{2} \|2 \pi \zeta\|^2}{2\nu} \right)^{-\left(\nu+\tfrac{d}{2}\right)}.
    \end{equation}

In order to compute the shared-endpoint correlation, we require the following expectations:
\begin{equation} \label{eqn: necessary expectations pre frequency space}
    \begin{aligned}
        & \textbf{a\hspace{0.035 in})}\indent  \mathbb{E}_{Z \sim \mathcal{N}(0,2 \sigma_x^2 I)}\left[k_{2T}([Z,0]; \nu_{2T}, l_{2T}) \right] = \mathbb{E}_{Z \sim \mathcal{N}(0,2 \sigma_x^2 I)}\left[k_{T}(Z; \nu_{T}, l_{T}) \right] \\
        & \textbf{a')}\indent  \mathbb{E}_{Z \sim \mathcal{N}(0,2 \sigma_x^2 I)}\left[k_{2T}([\sqrt{2}Z,0]; \nu_{2T}, l_{2T}) \right] = \mathbb{E}_{Z \sim \mathcal{N}(0,4 \sigma_x^2 I)}\left[k_{2T}([Z,0]; \nu_{2T}, l_{2T}) \right]  = \mathbb{E}_{Z \sim \mathcal{N}(0,4 \sigma_x^2 I)}\left[k_{T}(Z; \nu_{T}, l_{T}) \right] \\
        & \textbf{b\hspace{0.04 in})}\indent \mathbb{E}_{Z^{(1)},Z^{(2)}}\left[ k_{2T}([Z^{(1)},Z^{(2)}]; \nu_{2T}, l_{2T}) \right] 
    \end{aligned}
\end{equation}
where $[Z^{(1)},Z^{(2)}]$ is normal, mean zero, with covariance $\Sigma_z = (I + \textbf{1} \textbf{1}^{\intercal}) \otimes \Sigma_x$ (see equation \eqref{eqn: sigma and rho isotropic Gaussian}).

Expectations \textbf{a)} and \textbf{a')} are taken on a Mat\'ern kernel defined in $2T$ dimensions, but only the first $T$ inputs to the kernel vary. Since the Mat\'ern kernel is a univariate function of the magnitude of its argument, the Mat\'ern kernel in $d$ and $d'$ dimensions are related for the correct change of parameters. Therefore, we can replace the expectations \textbf{a)} and \textbf{a')} with an average over a $T$ dimensional space, provided we transform kernels correctly. 
The Mat\'ern kernel is parameterized so that this exchange is trivial - all parameters remain the same. See Appendix Section \ref{app: dimensionality reduction} for details:
    \begin{equation} \label{eqn: kernel dimension relationship}
        \kappa_{2T}([Z,0]; \nu, l) = \kappa_{T}(Z; \nu, l) 
    \end{equation}

Since the expectations are taken over Gaussian random variables, equation \eqref{eqn: gaussian expectation in frequency space} applies. Then, in frequency space: 
\begin{equation} \label{eqn: necessary expectations in frequency space}
    \begin{aligned}
        & \textbf{a, a') }\indent \mathbb{E}_{Z \sim \mathcal{N}(0,n \sigma_x^2 I)}\left[\kappa_T(Z;\nu,l) \right] = ((2 \pi) n \sigma_x^2)^{-\tfrac{T}{2}}  \mathbb{E}_{\zeta \sim \mathcal{N}(0, (2\pi)^{-2} (n \sigma_x^{2})^{-1} I)}[\hat{\kappa}_T(\zeta;\nu,l)] \\
        & \textbf{b)}\indent \hspace{0.26in}   \mathbb{E}_{Z^{(1)},Z^{(2)}}\left[ \kappa_u([Z^{(1)},Z^{(2)}]) \right] = ((2 \pi) \sqrt{3} \sigma_x^2)^{-T} \mathbb{E}_{\zeta \sim \mathcal{N}(0,(2\pi)^{-2} \sigma_x^{-2} S)}[\hat{\kappa}_{2T}(\zeta;\nu,l)]  
    \end{aligned}
\end{equation}
where:
$$
S  = \frac{1}{3} \left[\begin{array}{cc} 2 & -1 \\ -1 & 2 \end{array} \right] \otimes I_{T \times T} = \left( \left[\begin{array}{cc} 2 & 1 \\ 1 & 2 \end{array} \right] \otimes I_{T \times T} \right)^{-1} = \left[\begin{array}{cc} 2 I_{T \times T} & I_{T \times T} \\ I_{T \times T} & 2 I_{T \times T} \end{array} \right]^{-1}
$$
is the inverse covariance matrix for $[Z^{(1)},Z^{(2)}]$, $n = 2$ for \textbf{a)} an 4 for \textbf{a')}, and where the factor of $3$ inside the scaling comes from the determinant of the covariance of $[Z^{(1)},Z^{(2)}]$, as seen in equation \eqref{eqn: sigma z}.

The PSD of the Mat\'ern kernel, equation \eqref{eqn: Matern kernel in frequency space}, is proportional to a Student's T density. The Student's T density can be expressed as a scale mixture of normals, with inverse gamma variances (see Table \ref{table: mixture key}). Thus, the mixture method developed in Section \ref{sec: anisotropic} (see equation \eqref{eqn: isotropic mixture simplified}) can be used to rewrite the $T$ and $2T$ dimensional integrals in equation \eqref{eqn: necessary expectations in frequency space} with one-dimensional expectations over a scale parameter that is inverse gamma distributed. 
\begin{equation} \label{eqn: expectations in scale mixture form}
    \begin{aligned}
        & \textbf{a, a') } \indent \mathbb{E}_{V \sim \text{ Gamma}^{-1}(\nu,\nu)}\left[ (1 + n r^2 V)^{-\frac{T}{2}} \right]  \\
        & \textbf{b)} \indent \hspace{0.26in} \mathbb{E}_{V \sim \text{Gamma}^{-1}(\nu,\nu)}\left[ (1 + r^2 V)^{-\frac{T}{2}} (1 + 3 r^2 V)^{-\frac{T}{2}}\right]\\
    \end{aligned}
\end{equation}
where $n = 2$ for \textbf{a)} and $n = 4$ for \textbf{a')}.

Combining \textbf{a)}, \textbf{a')}, and \textbf{b)}:
\begin{snugshade}
\begin{result} \label{result 11}
\textbf{[$\rho$ under (\hyperref[table: Results + Assumptions]{1 - 5, 9})]}  When the attributes are drawn independently of $\mathcal{F}$, the flow is mean zero, and is drawn from a difference of stationary, isotropic, GP utilities with a Mat\'ern kernels:
\begin{equation} \label{eqn: Matern mixture simplified}
    \rho(r,\nu) = \frac{\mathbb{E}_{V \sim \text{ Gamma}^{-1}(\nu,\nu)}\left[ (1 + 2 r^2 V)^{-\frac{T}{2}} \right] - \mathbb{E}_{V \sim \text{ Gamma}^{-1}(\nu,\nu)}\left[ (1 + r^2 V)^{-\frac{T}{2}} (1 + 3 r^2 V)^{-\frac{T}{2}}\right]}{1 - \mathbb{E}_{V \sim \text{ Gamma}^{-1}(\nu,\nu)}\left[ (1 + 4 r^2 V)^{-\frac{T}{2}} \right]}.
\end{equation}
where $\text{Gamma}^{-1}(\alpha,\beta)$ is the inverse gamma distribution with parameters $\alpha$ and $\beta$. 
\end{result}
\end{snugshade}

Equation \eqref{eqn: Matern mixture simplified} offers a simple form for $\rho$ as a function of $r$, $\nu$, and the trait dimension $T$. All expectations are evaluated against the same  one-dimensional density. The parameters play distinct roles: $\nu$ determines the shape of distribution that the scale parameter $V$ is drawn from, $r$ determines the rate of change of the argument in $V$, and $T$ determines the shape of the argument. 
As $\nu$ goes to infinity, the inverse gamma distribution concentrates about $V = 1$, in which case $\rho(r,\nu)$ converges to $\rho(r)$ for the isotropic squared exponential model (see equation \eqref{eqn: rho isotropic SE gaussian traits}). The correlation for less regular ensembles, i.e.~for finite $\nu$, differs from the correlation in the SE case by the mixture over $V \neq 1$.

Unlike the squared exponential case, the integrals defined in equations \eqref{eqn: necessary expectations in frequency space} and \eqref{eqn: expectations in scale mixture form} cannot be solved in a closed form. So, we proceed via approximation. In the squared exponential case, we saw that the limiting behavior of $\rho$ in the smoothness and roughness limits could be used to produce reasonable approximations to $\rho$ for intermediate $r$. We will attempt the same analysis for finite $\nu$. First, we compute the correlation $\rho$, in the roughness and smoothness limits (i.e. as $r$ tends to infinity and to 0). Then, we extrapolate to intermediate $r$ via a Pad\'e approximation that merges the limits. 

\subsubsection{Roughness Limits for the Mat\'ern Kernel}

In frequency space, the roughness limits are easier to compute than the smoothness limits since roughness implies a concentrated PSD (small inverse covariance). To compute the roughness limit, we exploit isotropy, as in equation \eqref{eqn: sigma and rho iso iso Gaussian}, to re-express the expectations \eqref{eqn: necessary expectations in frequency space} in terms of a $\chi^2$ random variable. Next we Taylor expand the argument in $r^{-1}$ as $r \rightarrow \infty$, then drop lower order terms. This approximation is asymptotically valid for large $r$ since the $\chi^2$ distribution has raw moments of all orders, and since the distribution concentrates relative to the argument as $r$ approaches infinity. 

As before, let $n$ represent $2$ or $4$ for $\textbf{a)}$ and $\textbf{a')}$ respectively. Then:
\begin{equation} \label{eqn: star and star prime simplification}
    \begin{aligned}
        \mathbb{E}_{Z \sim \mathcal{N}(0,n \sigma_x^2 I)}\left[k_{2T}([Z,0]; \nu, l) \right] \simeq
        2^{\tfrac{T}{2}}\frac{\Gamma(\nu)}{\Gamma(\nu - \tfrac{T}{2})} n ^{-\tfrac{T}{2}} \Tilde{r}^{-T} \left( 1 - \frac{T}{2}\frac{\nu}{\nu-\tfrac{T}{2}}n^{-1}\Tilde{r}^{-2} + \mathcal{O}(\Tilde{r}^{-4}) \right).
    \end{aligned}
\end{equation}
where $\Tilde{r} = \sqrt{|2\nu-T|} \cdot r$ is a scaled version of the roughness coefficient, $r$ (see Appendix Section \ref{app: roughness limit calculations}).  

For \textbf{b)}, we perform the same approximation:
\begin{equation} \label{eqn: star star simplification}
    \begin{aligned}
       \mathbb{E}_{Z^{(1)},Z^{(2)}\sim \mathcal{N}(0, \sigma_x^2 S^{-1})}\left[k_{2T}([Z^{(1)},Z^{(2)}]; \nu, l) \right]       \simeq 2^T \frac{\Gamma(\nu+\tfrac{T}{2})}{\Gamma(\nu-\tfrac{T}{2})} 3^{\tfrac{-T}{2}}\Tilde{r}^{-2T} \left(1- \frac{4}{3}\frac{T}{2}\frac{\nu+\tfrac{T}{2}}{\nu-\tfrac{T}{2}}\Tilde{r}^{-2} + \mathcal{O}(\Tilde{r}^{-4})\right)
    \end{aligned}
\end{equation}

Finally, we return to the shared-endpoint correlation in equation \eqref{eqn: sigma and rho isotropic Gaussian} and re-substitute $r$ for instances of $\Tilde{r}$ (see Appendix Section \ref{app: roughness limit calculations}). In the limit as $r \rightarrow \infty$:
\begin{snugshade}
\begin{limit_result} \label{result 11: roughness}
\textbf{[Roughness Limit of $\rho$ under (\hyperref[table: Results + Assumptions]{1 - 5, 9})]}  Under the same assumptions of \ref{result 11}, in the limit as $r \rightarrow \infty$:
\begin{equation} \label{eqn: roughness limit}
    \rho(r) \simeq \frac{\Gamma(\nu)}{\Gamma(\nu-\tfrac{T}{2})(\nu-\tfrac{T}{2})^{\tfrac{T}{2}}}(\sqrt{2}r)^{-T}+ \mathcal{O}\left(r^{-(T + 2)} \right)
\end{equation}
\end{limit_result}
\end{snugshade}

As $\nu \rightarrow \infty$, the Mat\'ern kernel converges to a SE kernel. The $\nu$ dependent term in equation \eqref{eqn: roughness limit} approaches 1 as $\nu$ goes to infinity so, in the roughness limit, equation \eqref{eqn: roughness limit} reduces to $(\sqrt{2}r)^{-T}$, the roughness limit for the SE kernel \eqref{eqn: rho roughness limit isotropic SE}. Otherwise, the roughness limit is largely independent of $\nu$. In contrast, the smoothness limit will depend critically on $\nu$. 

\subsubsection{Smoothness Limits for the Mat\'ern Kernel}

The smoothness limit is considerably more involved in the joint limit when both $\nu$ and $r$ are small. The smaller $\nu$, the less regular the ensemble of functions, so small $\nu$ and small $r$ correspond to functions that are simultaneously smooth and irregular. Formally, small $\nu$ presents a problem since, when $\nu$ is small, naive expansion methods fail. 

Consider, for example, the scale-mixture form presented in equation \eqref{eqn: Matern mixture simplified}. When $r$ is small, it is tempting to Taylor expand the argument in $r$, then expand into a linear combination of the raw moments of the inverse gamma distribution. For examples see Appendix Section \ref{app: naive expansion}. This approach is invalid since the inverse gamma distribution only has finite moments up to degree $\nu$. The integral and the limit implicit in the sequence of partial sums, cannot be exchanged. 

Naive expansion fails since the tails of inverse gamma decay according to a power law of the form $v^{-(\nu + 1)}$. Similarly, the argument of the expectations decays according to a power law of the form $v^{-T/2}$ for large $v$.   
Slow tail decay changes the asymptotic behavior of $\rho(r)$ when $r$ is small. In all previous work, $\rho(r)$ approached $1/2$ from below at rate $r^2$ when $r$ is small. This behavior holds whenever the sample functions are almost surely second differentiable \cite{cebra2023similarity}. However, when $\nu \leq 2$, sample functions are almost never second differentiable. In this case, the analysis in \cite{cebra2023similarity} does not apply, so $\rho(r)$ need not differ from $1/2$ by a decrement of $\mathcal{O}(r^2)$. We will see that, when $\nu \in [1,2)$, $\rho(r)$ approaches $1/2$ as $r$ vanishes, but does so slower than $\mathcal{O}(r^2)$. At the critical value, $\nu = 1$, convergence is very slow. When $\nu < 1$, the sample draws do not admit accurate linear approximations, so $\rho(r)$ converges to a value smaller than $1/2$. 

To capture the asymptotics for small $r$, it is important to retain separate approximations to both a bulk term, that accounts for most of the mass of the distribution, and a tail term.  Expanding the bulk produces a power series in $r^2$ and a term of intermediate order, $r^{2 \nu}$ associated with the boundary between the tail and bulk terms. Expanding the tail produces terms of matching intermediate power. The smaller $\nu$, the smaller the fractional power of the intermediate order terms. For $\nu < 2$, the fractional power is less than relevant terms in the power series expansion of the bulk, so the tail term dominates some, or all, of the relevant terms in the bulk, thereby changing the limiting behavior.

To simplify the limiting analysis, we will first adopt an approximation to equation \eqref{eqn: Matern mixture simplified} that is asymptotically accurate as $r \rightarrow 0$ provided terms of order $r^4$ can be dropped. The approximation provides a lower bound on the exact value of $\rho$. It provides simpler results that capture the qualitative relation between regularity and correlation. 

The argument of the expectation \textbf{b)} differs in form from the arguments of expectations \textbf{a)} and \textbf{a')}. Both \textbf{a)} and \textbf{a')} have arguments of the form $(1 + n r^2 V)^{-T/2}$ while \textbf{b)} has an argument of the form $(1 + 4 r^2 V + 3 r^4 V^2)^{-T/2}$. When $r$ is small, $r^2 v$ dominates $r^4 v^2$ for all fixed $v$, indeed for all $v$ diverging slower than $r^{-2}$. Therefore, in the smoothness limit, it is tempting to drop the $3 r^4 V^2$ term, leaving $(1 + 4 r^2 V)^{-T/2}$. This produces a lower bound on $\rho$, since $3 r^4 V^2 > 0$. Then:
\begin{equation} \label{eqn: rho lower bound}
    \rho(r,\nu) \gtrsim \tilde{\rho}(r,\nu) = \frac{\mathbb{E}_{V \sim \text{ Gamma}^{-1}(\nu,\nu)}\left[ (1 + 2 r^2 V)^{-\frac{T}{2}} \right] - \mathbb{E}_{V \sim \text{Gamma}^{-1}(\nu,\nu)}\left[ (1 + 4 r^2 V)^{-\frac{T}{2}}\right]}{1 - \mathbb{E}_{V \sim \text{ Gamma}^{-1}(\nu,\nu)}\left[ (1 + 4 r^2 V)^{-\frac{T}{2}} \right]}
\end{equation}

Define the integral functions:
\begin{equation} \label{eqn: integral functions}
\begin{aligned}
    I(r,\nu) = \mathbb{E}_{V \sim \text{ Gamma}^{-1}(\nu,\nu)}\left[ (1 + 2 r^2 V)^{-\frac{T}{2}} \right], \quad J(r,\nu) = 1 - I(r,\nu).
\end{aligned}
\end{equation}

Then:
\begin{equation} \label{eqn: rho in integral functions}
    \rho(r,\nu) \gtrsim \tilde{\rho}(r,\nu) = \frac{J(\sqrt{2} r, \nu) - J(r, \nu)}{J(\sqrt{2} r, \nu)} = 1 - \frac{J(r,\nu)}{J(\sqrt{2} r, \nu)}.
\end{equation}
    
Equation \eqref{eqn: rho in integral functions} is easier to analyze, since the asymptotics of $\tilde{\rho}(r,\nu)$ are fully determined by the asymptotics of a single integral function, $J(r,\nu)$. Specifically, we only need to derive accurate approximations to $J(r,\nu)$ in the limit as $r \rightarrow 0$.

To analyze the asymptotic behavior of the integral function, $J(r,\nu)$, we split the integral over $v > 0$ into two pieces, a bulk term: $v \in (0,v_*(r))$, and a tail term: $v \geq v_*(r)$. The bulk and tail terms are:
\begin{equation}
J(r,\nu) = J_{\text{b}}(r,\nu) + J_{\text{t}}(r,\nu) \text{ where }
\begin{cases}
    & J_{\text{b}}(r,\nu) = \left(1 - \mathbb{E}[(1 + r^2 V)^{-T/2} \mid V < v_*(r)] \right) \text{Pr}(V < v_*(r)) \\
    & J_{\text{t}}(r,\nu) = \left(1 - \mathbb{E}[(1 + r^2 V)^{-T/2} \mid V \geq v_*(r)] \right) \text{Pr}(V \geq v_*(r)) \\
\end{cases}
\end{equation}
for $V \sim \text{Gamma}^{-1}(\nu,\nu)$.

Next, we develop separate approximations for the bulk and the tail that are each asymptotically accurate as $r$ vanishes. We use separate approximations for each based on whether the argument of the expectation, or the distribution of the random variable in the expectation, converge to an asymptotic expansion in the smoothness limit. On the bulk, we Taylor expand the argument, $(1 + r^2 V)^{-T/2}$ in $r$. On the tail, we adopt a Laurent expansion of the inverse Gamma distribution that approximates the inverse Gamma tail with a Pareto tail. We will choose $v_*(r)$ so that both approximations converge in the smoothness limit. We justify each approximation below.

\begin{enumerate}
    \item \textbf{Bulk: } The bulk term is an expectation over the conditional distribution of $V$ given $V < v_*(r)$. Unlike the naive expansion approaches discussed before, an expansion of the bulk into a linear combination of moments via a Taylor expansion of the argument, is asymptotically exact for two reasons. First, the argument is analytic at zero. Therefore, the Taylor expansion converges for all $v$ in the domain of integration provided $r^2 v_*(r)$ vanishes as $r \rightarrow 0$. In practice we will choose $v_*(r) = r^{-(2 - \epsilon)}$ for $\epsilon > 0$.  Second, the integral runs over a finite domain. Then, each \textit{conditional} moment remains finite, so the integral and the infinite sum implied by the Taylor series commute. The conditional moments can be evaluated exactly. Thus, setting $v_*(r) = r^{-(2 - \epsilon)}$, Taylor expanding, then exchanging the order of the expectation and the sum provides exact, analytic, asymptotics for the bulk in the smoothness limit. 

    \item \textbf{Tail: } To approximate the tail term we take the opposite approach. Instead of developing a Taylor series approximation to the argument of the expectation in small $r^2 v$, we develop a Laurent series approximation to the inverse gamma density in large $v$. This is valid since the density does not depend on $r$, so may be approximated assuming large $v$ when $v_*(r)$ diverges. To ensure $v_*(r)$ diverges we pick $\epsilon$ small.  
    Then, the tail of the mixture distribution approaches the tail of a Pareto distribution.
\end{enumerate}

After expanding the bulk and tail terms separately, we combine the approximations, then take $\epsilon$ to zero from above. For detailed analysis, see Appendix Section \ref{app: smoothess limit Matern}.

This approach yields the expansion:
\begin{equation} \label{eqn: expanded tail and bulk}
\begin{aligned}
    & \textbf{bulk: } \frac{1}{\Gamma(\nu)} \left[ \frac{T}{2} (\nu r^2 ) \left(\Gamma(\nu - 1) - \frac{(\nu r^{2 - \epsilon})^{\nu - 1}}{\nu - 1} \right) - \frac{T(T+2)}{8} (\nu r^2)^2 \Gamma(\nu - 2) + \mathcal{O}(r^6)\right] \\
    & \textbf{tail: } \frac{\nu^{\nu - 1}}{\Gamma(\nu)} r^{(2 - \epsilon) \nu} \left[1 - r^{2 \nu + \epsilon(T - \nu)}\left(1 - \frac{T}{2} \frac{\nu}{\nu - 1} r^{\epsilon} + \frac{T(T+2)}{8} \frac{\nu}{\nu - 2} r^{2 \epsilon} + \mathcal{O}(r^{3 \epsilon}) \right) \right]\\
\end{aligned}
\end{equation}

Notice that, for arbitrarily small $\epsilon$ the bulk includes a term of order $r^2$, of order $r^{2(\nu - 1) + 2} = r^{2 \nu}$ and of order $r^4$. The middle term is of intermediate order. It's order depends on the shape parameter $\nu$. The smaller $\nu$, the slower it approaches zero. It matches the lowest order term in the tail, which is also of order $r^{2 \nu}$ for $\epsilon$ sufficiently small. Together, these terms account for the slow tail decay ignored by naive expansions. 

If $\nu > 2$, then $2 \nu > 4$, so the second and fourth order terms dominate. If $\nu = 2$, then the middle term is of order $r^4$. If $\nu \in (1,2)$, then it is of intermediate order between $r^2$ and $r^4$. If $\nu = 1$, then it is order $r^2$, and if $\nu < 1$ then it is of intermediate order between $r^0$ and $r^2$. Thus, the dominating term, and the largest correction to the dominating term, depend on $\nu$. For $\nu > 2$, we can drop the intermediate term since it is dominated by the second and fourth order terms. In this case, the expansion of the bulk term matches the fourth order naive expansion. When $\nu < 2$ the intermediate term supplants the fourth order term. When $\nu < 1$ it supplants the second-order term. This produces distinct asymptotic regimes separated by the critical values $\nu = 1$, and $\nu = 2$, where the sample functions drawn from the corresponding ensemble gain almost sure differentiability, and almost sure second differentiability.

Taking $\epsilon$ to zero from above, and keeping only the terms in the bulk and tail that approach zero slowly enough to contribute to the limiting behavior of $\rho$ leaves:
\begin{equation} \label{eqn: expanded tail and bulk by nu}
    \begin{aligned}
        & \textbf{bulk: } \begin{cases}
            \frac{1}{\Gamma(\nu)} \frac{T}{2}(\nu r^2) \left[ \Gamma(\nu - 1) - \frac{(\nu r^2)^{\nu - 1}}{\nu - 1}\right] = \mathcal{O}(r^2) + \mathcal{O}(r^{2 \nu}) & \text{ if } \nu \in (0,2) \setminus 1 \\
            T r^2 (\log(r^{-1}) - \gamma) = \mathcal{O}(r^2 \log(r^{-1})) + \mathcal{O}(r^2) & \text{ if } \nu = 1 \\
            \frac{1}{\Gamma(\nu)} \left[\frac{T}{2} \Gamma(\nu - 1) (\nu r^2) - \frac{T(T+2)}{8} \Gamma(\nu - 2) (\nu r^2)^2 \right] = \mathcal{O}(r^2) + \mathcal{O}(r^4) & \text{ if } \nu > 2
        \end{cases} 
        \\
        \\
        & \textbf{tail: } \hspace{2.2mm}\begin{cases}\frac{\nu^{\nu - 1}}{\Gamma(\nu)} r^{2 \nu} + \mathcal{O}(r^{4 \nu}) = \mathcal{O}(r^{2 \nu}) & \text{ if } \nu \leq 2 \\
        \text{negligible} & \text{ if } \nu > 2
        \end{cases}
    \end{aligned}
\end{equation}

Finally, substituting each expansion into equation \eqref{eqn: rho in integral functions} gives separate asymptotic approximations to $\tilde{\rho}(r,\nu)$ in the smoothness limit:
\begin{snugshade}
\begin{limit_result} \label{result 11: smoothness}
    \textbf{[Smoothness Limit of $\rho$ under (\hyperref[table: Results + Assumptions]{1 - 5, 9})]}  Under the same assumptions of \ref{result 11}, in the limit as $r \rightarrow 0$:
    \begin{equation} \label{eqn: rho matern smoothness limits}
        \rho(r,\nu) > \tilde{\rho}(r,\nu) \simeq \begin{cases}
            \left(1 - (\frac{1}{2})^{\nu} \right) \left[ 1 - \frac{2^{1 - \nu} - 1}{2^{\nu} - 1} \Gamma(\nu) (2 \nu)^{1- \nu} r^{2(1 - \nu)} \right] = (1 - 2^{-\nu}) - \mathcal{O}(r^{2(1 - \nu)}) & \text{ if } \nu \in (0,1)\\
            \frac{1}{2} \left[1 - \frac{1}{2}| \log_2(r)|^{-1} \right] \hspace{1.375 in} = \frac{1}{2} - \mathcal{O}(|\log_2(r)|^{-1}) & \text{ if } \nu = 1 \\
            \frac{1}{2} \left[1 - \frac{(2 \nu)^{\nu - 1}}{\Gamma(\nu)} \left[2^{\nu - 1} - 1 \right] \frac{\nu(T-2) + 2}{T} r^{2(\nu - 1)} \right] \hspace{0.05 in} = \frac{1}{2} - \mathcal{O}(r^{2(\nu - 1)}) & \text{ if } \nu \in (1, 2) \\
        \end{cases}
    \end{equation}
    where the lower bound is asymptotically accurate for $\nu \in [1, 2)$. If $\nu > 2$, then the tail term is dominated by the lowest order terms in the bulk, so:
    \begin{equation} \label{eqn: rho matern smoothness limits large nu}
        \rho(r) \simeq \frac{1}{2} \left[ 1 - \frac{T-1}{2} \frac{\nu}{\nu - 2} r^2 \right] = \frac{1}{2} - \mathcal{O}(r^2) \text{ if } \nu > 2.
    \end{equation}
\end{limit_result}
\end{snugshade}

Result \ref{result 11: smoothness} identifies three distinct regimes for the smoothness limit, separated by critical values of $\nu$. These are $(0,1)$, $(1,2)$, and $(2,\infty)$ with $\nu = 1$ and $\nu = 2$ as critical shape parameters where the smoothness limit changes behavior. 

In particular, for $\nu \geq 1$, $\tilde{\rho}(r) \leq \rho(r)$ approaches $1/2$ as $r \rightarrow 0$. Thus, as long as functions sampled from the Gaussian process admit fractional derivatives of order arbitrarily close to, or greater than one, then, $\rho(r)$ saturates at the shared-endpoint correlation for linear functions in the smoothness limit. Previous work \cite{cebra2023similarity} had shown the same limiting behavior in an equivalent limit, provided the sample functions were second differentiable. Here, we see that second differentiability is not required, as linear approximations may dominate on small spatial scales even when the sample functions are not second differentiable (recall Figure \ref{fig: Matern regularity}). 

That said, the rate of convergence to $1/2$ does depend on the shape parameter $\nu$, and the second differentiability of the sample draws. If $\nu > 2$, then sample functions are second differentiable, and the shared-endpoint correlation approaches its maximum value from below at rate $\mathcal{O}(r^2)$ \cite{cebra2023similarity}. However, if $\nu \in [1,2)$ then $\rho$ approaches $1/2$ more slowly. For $\nu \in (1,2)$, it approaches at rate $\mathcal{O}(r^{2(\nu - 1)})$. This rate starts at $\mathcal{O}(r^2)$ when $\nu - 2$, crossing $\mathcal{O}(r)$ for $\nu = 1.5$ and approaches $r^0$ as $\nu$ approaches 1. At $\nu = 1$, convergence no longer occurs at a power law rate, and is, instead, logarithmic. Thus, while first differentiability is not needed to for local linear approximations to dominate, when the sample functions are fractionally differentiable for $\nu \in (1,2)$, $\rho(r)$ approaches $1/2$ very slowly.

For $\nu < 1$, sample functions are not approximately linear on small scales (see Figure \ref{fig: Matern regularity}), so $\rho(r)$ no longer approaches $1/2$. For the exact limit, repeat the same asymptotic analysis developed here, only using $\rho(r)$ instead of $\tilde{\rho}(r)$. We will confirm numerically that, when $\nu < 1$, $\rho(r)$ does not approach $1/2$ in the smoothness limit. Instead, it saturates at a value between $1 - 2^{-\nu}$ and 1, and decreases monotonically as $\nu$ approaches 0. The shared-endpoint correlation approaches this limit from below at rate $\mathcal{O}(r^{2(1 - \nu)})$ so approaches slowest for $\nu$ near 1, and fastest for $\nu$ near zero. As $\nu$ approaches zero, convergence occurs at a rate arbitrarily close to $\mathcal{O}(r^2)$. 

\subsubsection{Intermediate Roughness}

To conclude our analysis, we use a two-point Pad\'e approximation to merge the smoothness and roughness limits of $\rho(r)$. This produces an approximation to $\rho(r)$ that is asymptotically accurate to $\tilde{\rho}(r)$ as $r$ approaches 0, and to $\rho(r)$ as $r$ infinity. As in the smoothness limits, we use the simpler lower bound on $\rho(r)$, $\tilde{\rho}(r)$, since it has the same qualitative behavior, admits a simpler form, and nearly exact calculation of $\rho(r)$ is possible through numerical integration of \eqref{eqn: Matern mixture simplified} or \eqref{eqn: sigma and rho iso iso Gaussian}. 

To marry the smooth and rough limits, we adopt an approximant of the form:
\begin{equation}
    \rho(r;\nu,T) \approx C(\nu,T) (1 + S(r;\nu,T) + R(r;\nu,T))^{-1}
\end{equation}
where $C(\nu,T)$ is a constant chosen to match $\lim_{r \rightarrow 0} \tilde{\rho}(r)$, $S(r;\nu,T)$ is chosen such that $C(\nu,T) (1 + S(r;\nu,T))^{-1}$ has the same asymptotic approximation in the smoothness limit as $\tilde{\rho}$ (see Result \ref{result 11: smoothness}), and $R(r;\nu,T)$ is chosen so that $C(\nu,T) R(r;\nu,T)^{-1}$ has the same asymptotic approximation as $\rho(r;\nu,T)$ in the roughness limit. These terms are chosen so that $S$ dominates $R$ when $r$ is sufficiently small, and $R$ dominates $S$ and 1 when $r$ is sufficiently large. We extend the approximation developed for $\nu \in (1,2)$ to $\nu = 2$ since $\rho(r;\nu,T)$ should be a continuous function of $\nu$. Similar extensions are not possible from above, since the approximation developed above is ill-behaved as $\nu$ approaches 2. Similarly, both approximations developed for $\nu$ about one, fail to continue to accurate approximations at $\nu = 1$, since convergence in the smoothness limit occurs logarithmically, i.e. slower than any power law, when $\nu = 1$, but only when $\nu = 1$.

\begin{snugshade}
\begin{result} \label{result 12}
    \textbf{[Pad\'e Approximations for $\rho$ under (\hyperref[table: Results + Assumptions]{1 - 5, 9})]}  When the attributes are drawn independently of $\mathcal{F}$, the flow is mean zero, and is drawn from a difference of stationary, isotropic, GP utilities with a Mat\'ern kernel:
    \begin{equation} \label{eqn: pade approx to rho}
        \rho(r;\nu,T) \approx \begin{cases}
            \left(1 - (\frac{1}{2})^{\nu} \right) \left(1 + \frac{2^{1 - \nu} - 1}{2^{\nu} - 1} \Gamma(\nu) (2 \nu)^{1- \nu} r^{2(1 - \nu)} + \tfrac{2^{\nu + T/2}}{2^{\nu} + 1} G(\nu,T) r^T \right)^{-1} & \text{ if } \nu \in (0,1) \\
            \frac{1}{2} \left(1 + \frac{1}{2}| \log_2(r)|^{-1} + 2^{\tfrac{T}{2}-1} G(\nu,T) r^T \right)^{-1} & \text{ if } \nu = 1 \\
            \frac{1}{2} \left(1 + \frac{(2 \nu)^{\nu - 1}}{\Gamma(\nu)} \left[2^{\nu - 1} - 1 \right] \frac{\nu(T-2) + 2}{T} r^{2(\nu - 1)} + 2^{\tfrac{T}{2}-1}G(\nu,T) r^T \right)^{-1} & \text{ if } \nu \in (1,2] \\
            \frac{1}{2} \left(1 + \tfrac{T-1}{2} \frac{\nu}{\nu - 2} r^2 + 2^{\tfrac{T}{2}-1} G(\nu,T) r^T \right)^{-1} & \text{ if } \nu > 2 \\
        \end{cases}
    \end{equation}
    where:
    \begin{equation}
        G(\nu,T) = \tfrac{\Gamma(\nu)}{\Gamma\left(\nu - \tfrac{T}{2} \right)\left(\nu - \tfrac{T}{2}\right)^{T/2}}.
    \end{equation}
    The Pad\'e approximations differ in their smoothness limits, but converge in the roughness limit. The final approximation converges to the Pad\'e approximation to $\rho(r)$ for a SE kernel in the limit as $\nu \rightarrow \infty$. 
\end{result}
\end{snugshade}

To test the accuracy of these approximations, we evaluated both integral expressions for $\rho$ (equations \eqref{eqn: Matern mixture simplified} and \eqref{eqn: sigma and rho iso iso Gaussian}), compared to the simplified lower bound $\tilde{\rho}$ (see equation \eqref{eqn: rho in integral functions}), then compared to two-point Pad\'e approximations \eqref{eqn: pade approx to rho}. We tested each approach for dimensions $T = 2$ to $5$, for $\nu$ ranging from near 0 to 5, and for $r$ ranging from $10^{-5}$ to $10^{5}$. To validate the integrals, we also approximated $\rho$ by sampling random utility functions from the generative model.

In all cases, the two integral approaches agreed closely, and matched the Monte Carlo approach. The lower bound, $\tilde{\rho}$ provides an extremely accurate approximation for all $\nu \geq 1$. For $\nu < 1$, $\tilde{\rho}$ underestimates $\rho$, especially in the smoothness limit. Regardless, the two show the same qualitative dependence on $r$, $\nu$, and $T$. Moreover, the degree of underestimation is small for $T$ much larger than 5, especially at $\nu = 0.5$, the most commonly used choice beneath 1 \cite{williams2006gaussian}. 

\begin{figure}[t] 
    \centering
    \includegraphics[trim = 0 60 0 0, clip, width = \textwidth]{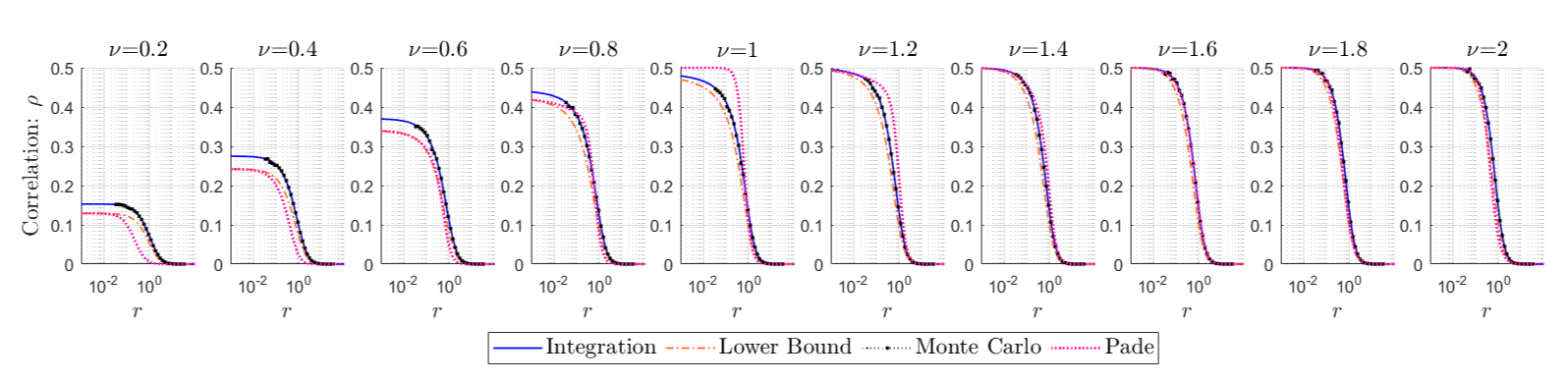}
    \includegraphics[trim = 0 25 0 0, clip, width = \textwidth]{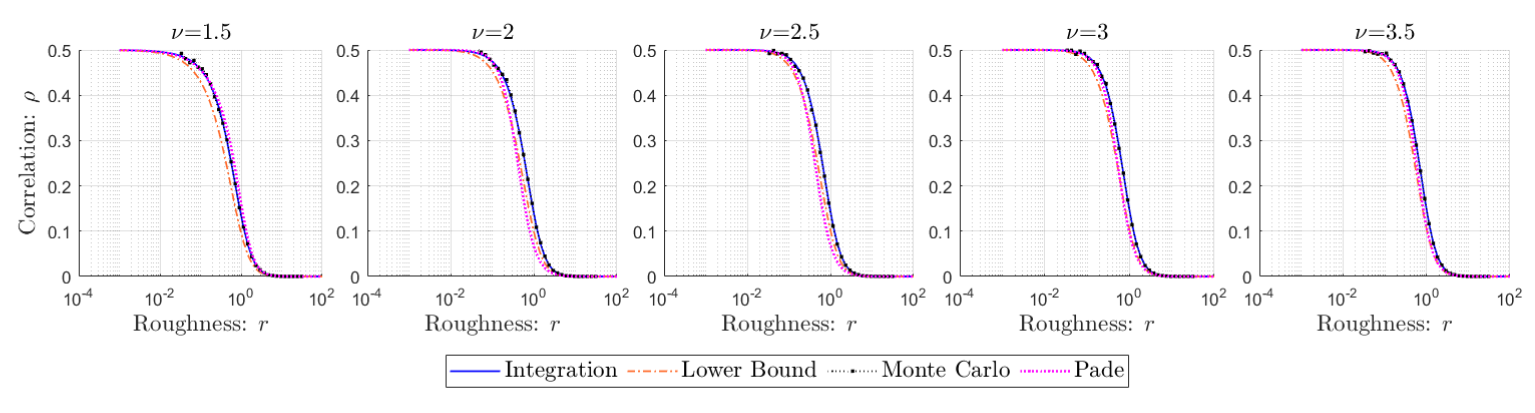}
    \caption{Shared-endpoint correlation $\rho$ as a function of the roughness $r$ given an isotropic Mat\'ern utility and isotropic Gaussian traits. Here $r$ varies from $10^{-4}$ to $10^{2}$ and $T = 3$. The shape parameter $\nu$ varies from 0.2 to 3.5. The solid blue line is integration using the $\chi$-squared density. The dashed orange line is the scale mixture lower bound $\tilde{\rho}$ (see equation \eqref{eqn: rho in integral functions}) to equation \eqref{eqn: Matern mixture simplified}, which is exact. The black dots are Monte Carlo estimates generated by sampling utility functions and trait vectors. The dotted magenta line is the Pad\'e approximation \eqref{eqn: pade approx to rho}. In the case when $\nu = 1$ we extend the solutions for $\nu < 1$ and $\nu \in (1,2)$, which yield the same estimate for $\nu = 1$.  }
    \label{fig: nu validation}
\end{figure}

To validate the Pad\'e approximations, we first tested  $\nu \in [1.5, 3.5]$. At $\nu = 1.5$, the sampled functions are almost surely differentiable, but not second differentiable. For $\nu > 2$, they are almost surely second differentiable. For $\nu > 3$, they closely resemble draws from the matching SE model. Figure \ref{fig: nu validation} compares the Pad\'e approximation in equation \eqref{eqn: pade approx to rho} to numerical integration, Monte Carlo simulation, and numerical integration of the lower bound $\tilde{\rho}$. 

In general, the Pad\'e approximation extrapolates well to intermediate $r$ for $\nu > 1.5$, closely matching both direct numerical integration and Monte Carlo simulation. The approximation improves as $\nu$ increases, and is very accurate for $\nu > 2.5$. Importantly, it remains accurate for $\nu = 2$, where the Pad\'e approximation was extended by continuity. For $\nu$ greater than, but close to 2, the approximation underestimates the true correlation for intermediate $r$. For $\nu$ between 1.5 and 2, and for $T \geq 4$, the approximation overestimates the correlation. We will see that this is a result of an extrapolation failure near $\nu = 1$, where accurate asymptotic approximations fail to produce an accurate Pad\'e approximation for intermediate $r$. In all cases, the Pad\'e approximation is qualitatively accurate; it is sigmoidal, monotonically decreasing, matches the correct limits, and decays over the correct range of $r$ at roughly the correct rate.

Next, we tested $\nu \in (0.2,2)$. The top row of Figure \ref{fig: nu validation} illustrates the results.
%
%
As before, exact integration closely matches Monte Carlo simulation. For $\nu < 1.5$, the Pad\'e approximation is considerably worse, with larger extrapolation errors away from the smoothness limit. In particular, extrapolation is quite poor for $\nu \in (0.8,1.2)$, with dramatic errors for $\nu = 1$. These extrapolation errors correspond to neglected higher order terms in the smoothness limit that are negligible when $r$ is sufficiently small, but contribute for intermediate $r$. In this case, we have opted to extend the neighboring solutions to highlight this effect. The asymptotic behavior in the smoothness limit is distinctly different when $\nu = 1$ than for any $\nu$ near 1. At $\nu = 1$, $\rho$ still approaches $1/2$, but does so logarithmically, thus more slowly, than when converging according to any power law. Using the correct (logarithmic) asymptotics \eqref{eqn: pade approx to rho} provides a more accurate approximation for small $r$, but produces large extrapolation errors as $r$ approaches 2.  Here we observe the dying gasp of the local linear model. The sampled functions are almost linear on small neighborhoods, but the errors in any local linear approximation decay very slowly as the neighborhood shrinks.

For $\nu < 1$, $\rho(r;\nu,T)$ no longer approaches $1/2$ since there is no scale at which the sample functions are approximately linear. For example, $\lim_{r \rightarrow 0} \rho(r) \approx 0.15$ when $\nu = 0.2$. In general, the less regular the sampled functions, the less smoothness implies correlation. In this range, $\tilde{\rho}$ underestimates $\rho$ by a small margin. Importantly, for processes  with approximately exponential kernels ($\nu \approx 0.5$), the Pad\'e approximation accurately predicts $\tilde{\rho}(r;\nu,T)$. 

\subsection{Generalization}

\subsubsection{Regularity and Correlation} \label{sec: generalization regularity}

The bulk-tail approach used to compute the smoothness limit of $\rho$ given a Mat\'ern model can be generalized to establish the convergence rates for $\rho$ given only the mean-square differentiability of the Gaussian process.

Suppose, that $\mathcal{F}$ is stationary, isotropic, and mean-square differentiable to any order $\mu < \nu$, but for no order $\mu \geq \nu$. Then every fractional derivative of the kernel up to order $2 \nu$ must be finite, while every fractional derivative of the kernel beyond order $2 \nu$ must be infinite, when the derivatives are evaluated at zero. So, after taking a Fourier transform, the PSD must admit all raw moments up to order $2 \nu$ but none beyond. If the PSD is expressed as a scale mixture, then the existence of a set of moments bounds the tail-decay rate of the mixture distribution in terms of a Pareto distribution.

%
%
%
In this case, the bulk-tail analysis developed for handling the tails of the Mat\'ern kernel can be generalized (compare equations \eqref{eqn: smoothness limit generalizes} and \eqref{eqn: rho matern smoothness limits}). 

\begin{snugshade}
\begin{theorem} \label{thm: regularity/correlation in smoothness}
    \textbf{[Regularity and Correlation in the Smoothness Limit under (\hyperref[table: Results + Assumptions]{1 - 6, 8, 10})]}  When the attributes are drawn independently $\mathcal{F}$, the flow is mean zero, is drawn from a difference of stationary, isotropic, GP utilities with a mixture of SE kernels, the attributes are normally distributed with covariance $\Sigma_x = \sigma_x^2 I$, and admits all fractional mean-square partial derivatives  of orders $\mu < \nu$, and none of order $\mu \geq \nu$, then as  $\lim_{r \rightarrow 0} \rho(r) = \rho_s \leq \tfrac{1}{2}$:
    \begin{equation} \label{eqn: smoothness limit generalizes}
        \rho_s - \rho(r) = \begin{cases} \mathcal{O}(r^{2(1 - \nu}) & \text{ if } \nu < 1 \\
        \mathcal{O}(|\log(r)|^{-1}) & \text{ if } \nu = 1 \\
        \mathcal{O}(r^{2(\nu - 1)}) & \text{ if } \nu \in (1,2) \\
        \mathcal{O}(r^2) & \text{ if } \nu \geq 2
        \end{cases} 
    \end{equation}
    and where $\rho_s = \tfrac{1}{2}$ if $\nu > 1$.
\end{theorem}
\end{snugshade}

\textbf{Proof Outline:} The proof follows largely as in the Mat\'ern case, without tracking specific integration constants. First, expand $\hat{\kappa}$ as a scale mixture of normals to recover equation \eqref{eqn: Matern mixture simplified} with a generic mixture distribution $\pi_v$. This is valid as all stationary isotropic kernels may be expressed as a scale mixture of Gaussians \cite{stein1999interpolation}. Let $\pi_V$ denote the mixture distribution over the variance. Then:
\begin{equation} \label{eqn: generalized frequency mixture}
        \rho(r) = \frac{\mathbb{E}_{V \sim \pi_V}\left[ (1 + 2 r^2 V)^{-\frac{T}{2}} \right] - \mathbb{E}_{V \sim \pi_V}\left[ (1 + r^2 V)^{-\frac{T}{2}} (1 + 3 r^2 V)^{-\frac{T}{2}}\right]}{1 - \mathbb{E}_{V \sim \pi_V} \left[ (1 + 4 r^2 V)^{-\frac{T}{2}} \right]}.
\end{equation}

Next, notice that any raw moment of order $\mu$ of the PSD is proportional to a raw moment of $\pi_V$ of order $\mu/2$. Therefore, the mixture distribution must admit all fractional moments up to order $\nu$, but not beyond. This is only possible if $\pi_V(v)$ decays at rate $v^{-(\nu + 1)}$ for $v$ sufficiently large. Associated tail bounds provide enough control over the tail decay rate to replicate the bulk-tail results derived for the Mat\'ern model up to order of convergence.

In particular, divide each expectation over the variance $V$ into a bulk and a tail at $v_*(r) = r^{-(2 - \epsilon)}$, enforce the tail-bounds for $v \geq v_{**}$ for some $v_{**}$ sufficiently large, then take the limits as $r$ and $\epsilon$ go to zero in order. Then:
\begin{equation}
    \mathbb{E}_{V \sim \pi_V}\left[1 - ( 1 + r^2 V + c r^4 V^2)^{-T/2} \right] = \begin{cases}
        \mathcal{O}(r^2) + \mathcal{O}(r^{2\nu}) + \mathcal{O}(r^4) & \text{ if } \nu \neq 1 \\
        \mathcal{O}(r^2 \log(r)) + \mathcal{O}(r^2) & \text{ if } \nu = 1
    \end{cases} 
\end{equation}
where $c$ is some constant greater than or equal to zero. For $c = 0$ we can use the same tail analysis as developed for $J$ under the Mat\'ern model. This provides a lower bound. For $c > 0$ we derive an upper bound on the tail by replacing $1 - (1 + \hdots)^{-T/2}$ with 1. Both approaches provide tail terms of matching order as functions of $r$.

Substitution into equation \eqref{eqn: generalized frequency mixture} recovers the convergence rates stated in equation \eqref{eqn: smoothness limit generalizes}. For $\nu > 1$, the limit is recovered by comparing the ratio of the slowest decaying terms in the numerator and denominato. For $\nu > 1$, this ratio always reduces to $1/2$, since the lowest order terms correspond to a local linearization of the sample function.

For details see Appendix Section \ref{app: generalized smoothness limit} $\blacksquare$

\subsubsection{Beyond Gaussian Processes} \label{sec: beyond GPs}

Our work appears limited to GP models. 
However, our results are more general than they appear since $\rho$ only depends on the covariance function. All processes with well-defined pairwise marginals admit a covariance function. Therefore:
\begin{snugshade} 
\begin{theorem} \label{thm: rho for equivalent GP}
    \textbf{[$\rho$ for Equivalent GP's]} Consider a mean-zero stochastic process $\mathcal{F}$ on $\Omega \times \Omega$, with covariance function $k([x,y],[x',y']) = \mathbb{C}[\mathbb{F}(x,y),\mathbb{F}(x',y')]$ satisfying $\mathcal{F}(x,y) = - \mathcal{F}(y,x)$ with probability one. Let $\mathcal{F}'$ denote the corresponding mean-zero Gaussian process with covariance function $k(\cdot,\cdot)$. Then, provided the end-point attributes are drawn independently from $\mathcal{F}$, $\rho$ for $\mathcal{F}$ will equal $\rho$ for $\mathcal{F}'$.
\end{theorem}
\end{snugshade}

Thus, our results extend to all stochastic processes with matching covariance functions. We've emphasized the story using Gaussian processes since the example covariance functions are most popular, and are a sufficient parameterization, in the GP setting. Theorem \ref{thm: rho for equivalent GP} is useful since it allows an analyst to predict $\rho$ for a stochastic process by recalling the result for the equivalent GP.  Note that Theorem \ref{thm: rho for equivalent GP} does not extend Theorem \ref{thm: regularity/correlation in smoothness} to arbitrary stochastic processes since the relations between sample regularity and covariance function used to prove Theorem \ref{thm: regularity/correlation in smoothness} may not generalize.

\section{Discussion} \label{sec: discussion}

The shared-endpoint correlation provides insight into the degree of structure in randomly sampled flows within a network. The size of $\rho$ predicts the bulk organization of flows: values close to 0 indicate disorder and values close to $\frac{1}{2}$ imply order. In this paper, we have examined the relationship between the shared-endpoint correlation, and the parameters of a trait-performance model \cite{strang2022network}, with performance sampled from a Gaussian process. Table \ref{table: results summary} summarizes the qualitative relations observed. Here, we discuss potential generalizations of increasing breadth and avenues for future work. 

\begin{table}[] 
\centering
\begin{tabular}{l|l|l|l}
\hline
               & Roughness  & Dimension  &  Regularity \\ \hline \hline
 Symbol: & $r$ & $T$ & $\nu$ \\
Meaning: & relative length scales $\sigma_x/l$ & \# of attributes to specify flow & regularity of function producing flow \\ \hline \hline

Trend: & $r \nearrow$ then  $\rho \searrow$ & $T \nearrow$ then $\rho \searrow$ & $\nu \nearrow$ then $\rho \nearrow$ \\
 & rougher is less correlated & higher dim. is less correlated & more regular is more correlated \\ \hline \hline
 If smooth: & \colorbox{gray!20}{$\lim_{r \rightarrow 0} \rho = \rho_s = \tfrac{1}{2}$} * & \colorbox{gray!20}{$\rho_s - \rho = \mathcal{O}(r^2)$ if $T > 1$} ** & * \hspace{0.02 in} $\rho_s \leq \tfrac{1}{2}$ with equality iff $\nu \geq 1$  \\
 (all $r$ small) & & \colorbox{gray!20}{\hspace{0.37 in} $= \mathcal{O}(r^4)$ if $T = 1$} & ** $\rho_s - \rho = \mathcal{O}(r^{\min(2,2(\nu - 1))})$ if $\nu \neq 1$ \\
 & & & \hspace{0.7 in} $\mathcal{O}(|\log(r)|^{-1})$ \hspace{0.12 in} if $\nu = 1$ \\ 
 If anisotropic & use arithmetic average & & \\ \hline \hline
 If rough: & $\lim_{r \rightarrow \infty} \rho = 0$ & $\rho = \mathcal{O}(r^T)$ & only constants depend on $\nu$ \\
 (any $r$ large) & & & \\
  & & & \\
 If anisotropic: & use geom. av. of large $r$ & $\rho = \mathcal{O}(\bar{r}^{|\mathcal{R}|})$  & \\ \hline
\end{tabular}
\caption{Results Summary. The top row specifies the key modeling parameters. Then next row identifies how each modeling parameter influences the shared endpoint correlation $\rho$. All trends are monotonic, and each trend assumes only one parameter varies. The middle row summarizes the behavior of $\rho$ in the smoothness limit. Results boxed in gray are restatements of results established in \cite{cebra2023similarity}. All smoothness results are either stated assuming isotropy, or using the arithmetic average of the roughness values. All roughness values must be small for the smoothness limit to apply. The bottom row summarizes the behavior of $\rho$ in the roughness limit, when \textit{any} roughness value is large. Results are stated for the isotropic and anisotropic case. When anisotropic, $\bar{r}$ stands for the geometric average over the set of large $r$, $\mathcal{R}$.}
\label{table: results summary}
\end{table}

This paper reports precise and suggestive results for specific, important kernels as well as for related families that can be expressed as a mixture or whose power spectral densities decay at equivalent rates. The example kernels were chosen to compromise tractability and relevance. Future work could address other families of kernels, or, anisotropic models where the trait distribution is aligned to the sampled function $\mathcal{F}$.  The latter are relevant when the endpoint attributes evolve in response to $\mathcal{F}$ (c.f.~\cite{cebra2024almost}). The former may be relevant for different stochastic processes.

This paper restricts its attention to random flow models where the flow distribution on an edge depends only on the character of its endpoints, and is otherwise independent of the edge's topological context. As a result, the flows on disjoint subgraphs are independent. In particular, two edges have a correlated flow if and only if the edges share an endpoint. Edges that are farther removed in the edge graph do not share any information \cite{strang2022network}.

In some applications, the character of an edge's endpoints \textit{and} it's position within the network influences the flow on that edge. For instance, if neighboring nodes have correlated traits, then the edge flow on an edge will depend on properties of nodes it does not touch. There are many natural settings where the existence of a connection between nodes is more likely when the nodes are similar. In these cases, the flows on a pair of edges may remain correlated even among distant edges. Then, a correlative analysis would need distinct correlations for edge pairs at distinct distances, or for every pair of edges. In this sense, the trait-performance model studied here is appealingly simple. It is fully characterized by one correlation. We hope to extend our analysis by considering similar models with structured edge flow covariances that are controlled by the distance between edges in the edge graph.


This work, like its immediate predecessors, \cite{cebra2023similarity,cebra2024almost,strang2020applications,strang2022network}, describes a specific null model, and characterizes the expected features of flows sampled under that null. 
It is important that future work moves beyond characterization of expected flow behavior and develops hypothesis tests based on those expectations. 

\bibliographystyle{siam}
\bibliography{Refs.bib}

\pagebreak
\section{Appendices}

\subsection{Supplementary Discussion}

\subsubsection{Regularity of Gaussian Processes} \label{app: GP regularity}

We will focus on Mat\'ern kernels. To better understand their regularity, we recall some foundational results that relate the regularity of sample draws of a stationary GP to its kernel function.

\begin{enumerate}
    \item \textbf{Mean-Square (MS) Differentiability:} The differentiability of draws from a GP is well characterized \cite{adler2010geometry}. Derivatives are limits of difference operators, difference operators are linear combinations of samples, and linear combinations of Gaussian random variables are Gaussian. Therefore, derivatives of a Gaussian process are themselves a Gaussian process provided the variances and expectations of the sequence of differences remain finite and converge to values that do not depend on the choice of differencing scheme. The mean square limit of a sequence of Gaussian random variables $\{g_i\}_{i=1}^{\infty}$ exists and takes the value $g_*$ if $\mathbb{E}[|g_i - g_*|^2] \rightarrow 0$ as $i$ goes to infinity. Then, the mean square directional derivatives are defined $\partial_{\delta x} \mathcal{F}(x_i) = \lim_{h \rightarrow 0} \frac{\mathcal{F}(x + h \delta x) - \mathcal{F}(x)}{h}$ where the limit is meant in the mean squared sense \cite{williams2006gaussian}. 

    If $\mathcal{F}$ is drawn from a stationary GP then the $\mu^{th}$ directional derivative along any direction $\delta x$ exists and is finite at all $x$ if the $2\mu^{th}$ directional derivative of the kernel also exists and is finite at zero. That is, $\partial^{\mu}_{\delta x} \mathcal{F}(x_i)$ exists if $\partial^{2\mu}_{\delta x} \kappa(x)|_{x = 0}$ exists \cite{williams2006gaussian}. Here we allow $\mu$ fractional, where fractional derivatives are defined by fractional moments of the Fourier transform of the kernel. 
    
    In particular, the derivatives of any kernel about zero are related to the central moments of the Fourier transform of the kernel. The kernel of a stationary process, and its Fourier transform (the power spectral density) are related by Bochner's theorem \cite{bochner2005harmonic,gikhman2004theory}, and the Karhunen-Lo\`eve (KL) expansion \cite{daw2022overview,karhunen1946spektraltheorie,loeve1955probability}. Bochner's theorem states that every valid kernel can be expressed as the Fourier transform of a nonnegtive measure. When that measure admits a density, the KL expansion expresses the GP as the Fourier transform of a GP in frequency space, with coefficients drawn independently, with variance at frequency $\omega$ equal to the value of the density there. Accordingly, the Fourier transform of the kernel is the power spectral density. The tail decay rate of the power spectral density determines the relative power contributed by high-frequency components as the frequency diverges, so slower tail decay rates produce rougher processes with faster local variation. Indeed, the $2\mu^{th}$ directional derivative of the kernel at zero is a $2\mu^{th}$ central moment of the power spectral density, so $\mathcal{F}$ is $\mu^{th}$ order MS differentiable if and only if the power spectral density admits moments up to order $2k$. 
    
    The Mat\'ern kernels have a power spectral density whose tails decay according to a power law, with exponent determined by $\nu$ and the dimension of the input space. Thus, for finite $\nu$, draws from a GP with a Mat\'ern kernel are not infinitely MS differentiable, and only admit derivatives up to an order controlled by, and monotonically increasing in $\nu$. Indeed, draws from a GP with a Mat\'ern kernel are $\mu$ times MS differentiable if and only if $\mu < \nu$ \cite{williams2006gaussian}. 

    \item \textbf{Sobolev Spaces and Reproducing Kernel Hilbert Spaces (RKHS)} The regularity of draws from a GP can also be studied by asking for the probability that a draw from a GP will lie in a particular smoothness class. Together, Driscoll's zero-one law \cite{driscoll1973reproducing}, and Steinwart's Theorem \cite{steinwart2019convergence} give a clear answer when the smoothness classes are Sobolev spaces.

    A real-valued function $f: \Omega \subseteq \mathbb{R}^d \rightarrow \mathbb{R}$ is a member of the Sobolev space $W^{m}$ if $\sum_{n = 0}^m \sum_{|\alpha| = n} \int_{x \in \Omega} (\partial^{\alpha}_x f(x))^2 dx$ is finite. Here $\alpha$ denotes a multi-index $\alpha = [\alpha_1,\alpha_2,\hdots,\alpha_d] \in \mathbb{Z}^d$, and $\partial^{\alpha}_x = \prod_{j=1}^d \partial_{x_j}^{\alpha_j}$. The Sobolev space $W^{m}$ contains all functions whose $n^{th}$ order partial derivatives are square-integrable for all $n \leq m$. 

    If $\mathcal{F}$ is drawn from a GP with mean function in $W^{m}$ and kernel $\kappa$, what is the probability that $\mathcal{F} \in W^m$? 

    The answer depends on the class of functions that can be expressed as infinite linear combinations of the kernel. Recall that, given a kernel $\kappa$, the reproducing kernel Hilbert space (RKHS) of $\kappa$, $\mathcal{H}_{\kappa}$ is the completion (closure) of all finite linear combinations of the kernel $\kappa(\cdot - x)$ for some finite set of centers $x$. In other words, $\mathcal{H}_{\kappa}$ contains all functions that can be built as the limit of the convolution of the kernel with a finite set of Dirac-$\delta$ functions \cite{berlinet2011reproducing}.

    Driscoll's zero-one law states that, the probability $f \in \mathcal{H}_{\kappa}$ is zero \cite{driscoll1973reproducing}.\footnote{In contrast, maximum a posteriori, or mean estimators, to the sample path of a GP conditioned on finitely many observations are in $\mathcal{H}_{\kappa}$.} However, \cite{steinwart2019convergence} if $\mathcal{H}_{\kappa} \subseteq W^{k}$ then the probability that $\mathcal{F} \in W^{k - \frac{d}{2} - \epsilon}$ is one for all $\epsilon > 0$. In particular, the Mat\'ern kernels with regularity parameter $\nu$ have $\mathcal{H}_{\kappa} = W^{\nu + d/2}$, so, if $\mathcal{F}$ is drawn from a GP with a Mat\'ern kernel and regularity $\nu$, then $\mathcal{F}$ is contained inside $W^{k}$ for any $m < \nu$ \cite{kanagawa2018gaussian}.  
    
\end{enumerate}

Therefore, while $l$ determines the characteristic length scale over which $\mathcal{F}$ varies, $\nu$ determines its regularity:

\setcounter{fact}{1}

\begin{snugshade}
\begin{fact} \label{app fact: Fact 2}
    \textbf{[Regularity and $\nu$]} If $\mathcal{F}$ is drawn from a Mat\'ern process with regularity parameter $\nu$, then $\mathcal{F}$ is $\mu^{th}$ order mean-square differentiable for all $\mu < \nu$, but no $\mu \geq \nu$ \cite{williams2006gaussian}. Sample draws of the process are almost surely contained inside the Sobolev space $W^{m}$ if $m < \nu$ \cite{kanagawa2018gaussian,steinwart2019convergence}, and are almost never contained in $W^{m}$ if $m \geq \nu$ \cite{driscoll1973reproducing}.
\end{fact}
\end{snugshade}

\subsection{Supplementary Proofs}

\subsubsection{Proof of Theorem 2:} \label{app: thm 2}

\setcounter{theorem}{1}

\begin{snugshade}
\begin{theorem} \label{app thm: thm 2}
\textbf{[Skew Symmetric GP]} 
Suppose $\mathcal{F} \sim \textnormal{GP}(\mu_f,k_f)$, $\mathcal{F}:\Omega \times \Omega \rightarrow \mathbb{R}$. Then $\mathcal{F}(x,y) = - \mathcal{F}(y,x)$ for all $(x,y) \in \Omega \times \Omega$ if and only if $\mu_f(x,y) = -\mu_f(y,x)$ and $k_f([x,y],[x,y]) = -k_f([x,y],[y,x])$ for all $x,y$.
\end{theorem}
\end{snugshade}

\noindent \textbf{Proof: } First, suppose that $\mathcal{F}$ is skew symmetric. Then $\mathcal{F}(x,y) = -\mathcal{F}(y,x)$. It follows that $\mu_f(x,y) = \mathbb{E}_{\mathcal{F}}[\mathcal{F}(x,y)] = -\mathbb{E}_{\mathcal{F}}[\mathcal{F}(y,x)] = - \mu_f(y,x)$. It also follows that $\kappa_f([x,y],[x,y]) = \mathbb{V}_{\mathcal{F}}[\mathcal{F}(x,y)] = \mathcal{V}_{\mathcal{F}}[- \mathcal{F}(y,x)] = \mathcal{V}_{\mathcal{F}}[ \mathcal{F}(y,x)] = \kappa_f([y,x],[y,x])$. Moreover, $\kappa_f([x,y],[y,x]) = \mathbb{C}[\mathcal{F}(x,y),\mathcal{F}(y,x)] = -\mathbb{C}_{\mathcal{F}}[\mathcal{F}(x,y),\mathcal{F}(x,y)] = -\mathbb{V}_{\mathcal{F}}[\mathcal{F}(x,y)] = -\kappa_f([x,y],[x,y])$. The remaining equalities follow from the fact that $\kappa$ is a covariance so $\kappa(t,s) = \kappa(s,t)$ for any $s$ and $t$.

Suppose instead that $\mu_f(x,y) = -\mu_f(y,x)$ and $\kappa_f([x,y],[x,y])  = -\kappa_f([x,y],[y,x])$ for all $x,y$. Since $\mu_f(x,y) = -\mu_f(y,x)$, $\mathbb{E}_{\mathcal{F}}[\mathcal{F}(x,y)] = - \mathbb{E}_{\mathcal{F}}[\mathcal{F}(y,x)]$. Then, by the linearity of expectation, $\mathbb{E}_{\mathcal{F}}[\mathcal{F}(x,y) + \mathcal{F}(y,x)] = 0$. 

Now, any random variable with zero variance equals its expectation. So, if we can show that $\mathbb{V}_{\mathcal{F}}[\mathcal{F}(x,y) + \mathcal{F}(y,x)] = 0$ then $\mathcal{F}(x,y) + \mathcal{F}(y,x) = \mathbb{E}_{\mathcal{F}}[\mathcal{F}(x,y) + \mathcal{F}(y,x)] = 0$ so $\mathcal{F}(x,y) = -\mathcal{F}(y,x)$ for all $x,y$ and all $\mathcal{F}$ drawn from the GP. Therefore, it is enough to show that $\mathbb{V}_{\mathcal{F}}[\mathcal{F}(x,y) + \mathcal{F}(y,x)] = 0$. Expanding:
$$
\begin{aligned}
\mathbb{V}_{\mathcal{F}}[\mathcal{F}(x,y) + \mathcal{F}(y,x)] & = \mathbb{V}_{\mathcal{F}}[\mathcal{F}(x,y)] + 2 \mathbb{C}_{\mathcal{F}}[\mathcal{F}(x,y),\mathcal{F}(y,x)] + \mathbb{V}_{\mathcal{F}}[\mathcal{F}(y,x)] \\
& = \kappa_f([x,y],[x,y]) + \kappa_f([y,x],[y,x]) - 2 \kappa_f([x,y],[y,x]) = 0
\end{aligned}
$$
where the last equality follows since we assumed that $\kappa_f([x,y],[y,x]) = - \kappa_f([x,y],[x,y])$. Then, $\mathbb{V}_{\mathcal{F}}[\mathcal{F}(x,y) + \mathcal{F}(y,x)] = \kappa_f([y,x],[y,x]) + \kappa_f([x,y],[y,x])$. These also cancel each other out since $\kappa$ is a covariance function, thus $\kappa(t,s) = \kappa(s,t)$ for all $s$ and $t$, so $\kappa_f([x,y],[y,x]) = \kappa_f([y,x],[x,y]) = -\kappa_f([y,x],[y,x])$ by our original assumption after relabeling $x$, $y$, and $y$, $x$. $\blacksquare$

\vspace{0.1 in}

\subsubsection{Proof of Theorem 3} \label{app: generalized smoothness limit}

\begin{snugshade}
\begin{theorem} \label{app thm: regularity/correlation in smoothness}
    \textbf{[Regularity and Correlation in the Smoothness Limit under (\hyperref[table: Results + Assumptions]{1 - 6, 8, 10})]}  When the attributes are drawn independently $\mathcal{F}$, the flow is mean zero, is drawn from a difference of stationary, isotropic, GP utilities with a mixture of SE kernels, the attributes are normally distributed with covariance $\Sigma_x = \sigma_x^2 I$, and admits all fractional mean-square partial derivatives  of orders $\mu < \nu$, and none of order $\mu \geq \nu$, then as  $\lim_{r \rightarrow 0} \rho(r) = \rho_s \leq \tfrac{1}{2}$:
    \begin{equation} \label{appeqn: smoothness limit generalizes}
        \rho_s - \rho(r) = \begin{cases} \mathcal{O}(r^{2(1 - \nu}) & \text{ if } \nu < 1 \\
        \mathcal{O}(|\log(r)|^{-1}) & \text{ if } \nu = 1 \\
        \mathcal{O}(r^{2(\nu - 1)}) & \text{ if } \nu \in (1,2) \\
        \mathcal{O}(r^2) & \text{ if } \nu \geq 2
        \end{cases} 
    \end{equation}
    and where $\rho_s = \tfrac{1}{2}$ if $\nu > 1$.
\end{theorem}
\end{snugshade}

\noindent \textbf{Proof:} Suppose that the utility function is stationary. Let $\kappa$ denote the kernel and $\hat{\kappa}$ the corresponding power spectral density (PSD). 

Let $\alpha = [\alpha_1,\alpha_2,\hdots, \alpha_T]$ be a real valued vector with nonnegative entries and let $\mu = |\alpha| = \sum_{j} \alpha_j$. Then let $\partial_x^{\alpha} = \prod_{j} \partial_{x_j}^{\alpha_j}$ denote the mixed fractional partial derivatives. These partials exist in the mean-squared sense if and only if the corresponding fractional moment of the PSD exists, $\mathbb{E}_{\zeta \sim \hat{\kappa}}[\zeta^{2 \alpha}] = \mathbb{E}_{\zeta \sim \hat{\kappa}}[\prod_{j} \zeta_j^{2 \alpha_j}]$. 

If $\kappa$ is a stationary isotropic kernel, the $\hat{\kappa}$ can be expressed as a scale mixture of Gaussians with some mixture distribution $\pi_V$. Then:
\begin{equation}
    \mathbb{E}_{\zeta \sim \hat{\kappa}}[\zeta^{2 \alpha}] = \mathbb{E}_{V \sim \pi_V}\left[ \mathbb{E}_{\zeta|V}\left[ \zeta^{2 \alpha} \right]\right] = \mathbb{E}_{V \sim \pi_V}\left[\mathbb{E}_{\zeta|V}\left[ \prod_{j} \zeta_j^{2 \alpha_j} \right]\right] \propto \mathbb{E}_{V \sim \pi_V}\left[\prod_{j} V^{\alpha_j} \right] = \mathbb{E}_{V \sim \pi_V}[V^{\mu}].
\end{equation}
where the middle proportionality follows from the fact that the standard deviation is the scale parameter for a Gaussian distribution, and where the constants of proportionality are products of the $\alpha_j$ fractional moments of a Gaussian. For all integer $\alpha$ (standard derivatives) these are even moments, so are nonzero. Indeed, for all $\alpha$ that are not half-integers, the constants are nonzero. 

Therefore, the process is mean-square differentiable to order $\mu < \nu$ if and only if the mixture distribution $\pi_V$ admits finite raw moments of order $\mu < \nu$, but not beyond. Since $V \in [0,\infty)$, this requires that there exists constants $C_+ \geq C_- > 0$, and some threshold $v_{**} \geq 0$ such that:
\begin{equation}
    \pi_V(v) \in [C_-,C_+] v^{-(\nu + 1)} \text{ if } v \geq v_{**} .
\end{equation}

The corresponding shared-endpoint correlation is given by:
\begin{equation} \label{appeqn: generalized frequency mixture}
        \rho(r) = \frac{\mathbb{E}_{V \sim \pi_V}\left[ (1 + 2 r^2 V)^{-\frac{T}{2}} \right] - \mathbb{E}_{V \sim \pi_V}\left[ (1 + r^2 V)^{-\frac{T}{2}} (1 + 3 r^2 V)^{-\frac{T}{2}}\right]}{1 - \mathbb{E}_{V \sim \pi_V} \left[ (1 + 4 r^2 V)^{-\frac{T}{2}} \right]}.
\end{equation}

Let:
\begin{equation}
    J(r,c) = 1 - \mathbb{E}_{V \sim \pi_V}\left[(1 + r^2 V + c r^4 V^2)^{-T/2} \right].
\end{equation}

Then:
\begin{equation}
    \rho(r) = \frac{J(2r,\tfrac{3}{16}) - J(\sqrt{2} r,0)}{J(2r,0)}.
\end{equation}

We aim to show that, for any $c \geq 0$:
\begin{equation} \label{appeqn: integral function orders}
    J(r,c) = \begin{cases}
        \mathcal{O}(r^2) + \mathcal{O}(r^{2\nu}) + \mathcal{O}(r^4) & \text{ if } \nu \neq 1 \\
        \mathcal{O}(r^2 \log(r)) + \mathcal{O}(r^2) & \text{ if } \nu = 1.
    \end{cases} 
\end{equation}

If equation \eqref{appeqn: integral function orders} is true, then the theorem claim holds in all four $\nu$ regimes:

\begin{enumerate}
    \item If $\nu < 1$: $2 \nu < 2 < 4$ so the $2 \nu$ term dominates and the $4^{th}$ order terms are negligible. Then $\rho(r) = \tfrac{\mathcal{O}(r^{2 \nu}) + \mathcal{O}(r^2)}{\mathcal{O}(r^{2 \nu})} = \mathcal{O}(1) + \mathcal{O}(r^{2(1 - \nu)})$. 
    \item If $\nu = 1$: $r^2 \log(r)$ dominates $r^2$. Then $\rho(r) = \tfrac{\mathcal{O}(r^2 \log(r)) + \mathcal{O}(r^2)}{\mathcal{O}(r^{2} \log(r))} = \mathcal{O}(1) + \mathcal{O}(\log(r)^{-1})$.
    \item If $\nu \in (1,2)$: $2 < 2 \nu < 4$ so the second order term dominates and the $4^{th}$ order terms are negligible. Then $\rho(r) = \tfrac{\mathcal{O}(r^2) + \mathcal{O}(r^{2 \nu})}{\mathcal{O}(r^{2})} = \mathcal{O}(1) + \mathcal{O}(r^{2(\nu - 1)})$. 
    \item If $\nu \geq 2$: $2 < 4 \leq 2 \nu$ so the second order term dominates and the $2 \nu$ order termsa re negligible. Then $\rho(r) = \tfrac{\mathcal{O}(r^2) + \mathcal{O}(r^4)}{\mathcal{O}(r^{2})} = \mathcal{O}(1) + \mathcal{O}(r^{2})$. 
\end{enumerate}

So, to establish the Theorem 3, it is sufficient to show equation \eqref{appeqn: integral function orders}. Note, in each case, when we use $\mathcal{O}$ we assume that the unstated constants are nonzero.

As when deriving the asymptotic behavior of $\rho(r)$ for Mat\'ern processes, we proceed by dividing each integral into a bulk and tail component:
\begin{equation}
    J(r,c) = J_b(r,c) + J_t(r,c)
\end{equation}
where:
\begin{equation}
    \begin{aligned}
        & \textbf{Bulk: }  J_b(r,c) = \int_{v < v_{*}(r)}\left(1 -  (1 + r^2 v + c r^4 v^2)^{-T/2} \right) \pi_v(v) dv\\
        & \textbf{Tail: } J_t(r,c) = \int_{v \geq v_{*}(r)}\left( 1 -  (1 + r^2 v + c r^4 v^2)^{-T/2} \right) \pi_v(v) dv
    \end{aligned}
\end{equation}
and where $v_*(r) = r^{-(2 - \epsilon)}$ for some $\epsilon > 0$. 

Then, for sufficiently small $r$, $v_*(r) > v_{**}$, so the tail term can be bounded using the tail bound on the mixture distribution:
$$
J_b(r,c) \in [C_-,C_+] \int_{v \geq v_*(r)} \left( 1 -  (1 + r^2 v + c r^4 v^2)^{-T/2} \right) v^{-(\nu + 1)} dv.
$$

To produce a lower bound, notice that setting $c = 0$ makes $(1 + r^2 v + c r^4 v^2)$ smaller, so increases $(1 + r^2 v + c r^4 v^2)^{-T/2}$ so decreases $1 - (1 + \hdots)^{-T/2}$. Therefore, $J(r,c) > J(r,0)$. Thus:
$$
J_b(r,c) > C_- \int_{v \geq v_*(r)} \left(1 - (1 + r^2 v)^{-T/2} \right) v^{-(\nu + 1)} dv.
$$

We compute this integral in equation \eqref{appeqn: tail asymptotics}. The result is $\mathcal{O}(r^{(2 - \epsilon) \nu}) + \mathcal{O}(r^{(4 - 2\epsilon) \nu})$ so is $\mathcal{O}(r^{2\nu}) + \mathcal{O}(r^{4 \nu})$ after taking $\epsilon$ to zero from above.

To compute an upper bound, notice that $1 - (1 + \hdots)^{-T/2} < 1$ for all $r$, $c$, and $v$. Then:
$$
J_t(r,c) < C_+ \int_{v \geq v_*(r)} v^{-(\nu + 1)} dv = \frac{C_+}{\nu} v_*(t)^{-\nu} = \frac{C_+}{\nu} r^{(2 - \epsilon) \nu}
$$

So, after taking $\epsilon$ to zero from above, $J_t(r,c)$ is bounded above and below by functions of order $r^{2 \nu}$:
\begin{equation}
    J_t(r) = \mathcal{O}(r^{2 \nu}).
\end{equation}

Next, we need bounds on the bulk term. The argument of the expectation is analytic, so its Taylor series converges on an open interval containing zero. Since $v_*(r) r^2 = r^{\epsilon}$ and $\epsilon > 0$, the maximum input to the argument inside the integral defining the bulk term converges to zero as $r$ decreases. Therefore, for sufficiently small $r$, the Taylor expansion of the argumnet must converge over the entire region of integration. 

Taylor expanding the argument gives:
$$
1 - (1 + r^2 v + c r^4 v^2)^{-T/2} = \sum_{n=0}^{\infty} a_n r^{2n} v^n
$$
where:
$$
a_0 = 0, \quad a_1 = \frac{T}{2}, \quad a_2 = \frac{T}{2} \frac{T + 2 - 4 c}{4} 
$$

Then, for sufficiently small $r$:
\begin{equation}
    J_b(r,c) = \int_{v < v_*(r)} \left[\sum_{n=0}^{\infty} a_n r^{2n} v^n \right] \pi_V(v) dv = \sum_{n=0}^{\infty} a_n r^{2n} \int_{v=0}^{v_*(r)} v^n \pi_V(v) dv
\end{equation}
where the integral and infinite sum can be exchanged since each integral runs over a finite domain for finite $r$ and the infinite sum converges for sufficiently small $r$.

So, to evaluate the tail term, evaluate the conditional moments:
$$
M_n(r,\epsilon) = r^2 \int_{v=0}^{v_*(r)} v^n \pi_V(v) dv
$$

Split the integral defining the moment at $v_**$. Then:
$$
M_n(r,\epsilon) \in r^{2n} \int_{v=0}^{v_{**}} v^n \pi_v(v) dv + [C_-,C_+] r^{2n} \int_{v_{**}}^{v_*(r)} v^{n - (\nu + 1)} dv
$$

Since $v_{**}$ is a constant, the first term is $\mathcal{O}(r^{2n})$. If $\nu \neq 1$, then the second term is contained in the interval:
$$
\frac{[C_-,C_+]}{n - \nu} r^{2n} v^{n - \nu}|_{v_{**}}^{v_*(r)} = \frac{[C_-,C_+]}{n - \nu} r^{2n}(r^{-(2 - \epsilon)(n - \nu)} - v_{**}^{n - \nu} ) = \mathcal{O}(r^{2 \nu + \epsilon (n - \nu)}) + \mathcal{O}(r^{2n}) = \mathcal{O}(r^{2 \nu}) + \mathcal{O}(r^{2n})
$$
where the last equality assumes $\epsilon$ arbitrarily small. 

Therefore, when $\nu \neq 1$:
$$
J_b(r,c) = \sum_{n = 1}^{\infty} \mathcal{O}(r^{2n}) + \mathcal{O}(r^{2\nu}) = \mathcal{O}(r^2) + \mathcal{O}(r^{2 \nu}) + \mathcal{O}(r^4)
$$

In particular, the second order term is:
$$
\frac{T}{2} \left[ \int_{v < v_{**}} v \pi_V(v) dv \right] r^2 = b r^2
$$

When $\nu > 1$, the second order term dominates so:
\begin{equation}
    \rho_s = \lim_{r \rightarrow 0} \rho(r) = \lim_{r \rightarrow 0} \frac{J(2 r,\tfrac{3}{16}) - J(\sqrt{2} r,0)}{J(2 r,0)} =  \lim_{r \rightarrow 0} \frac{4 b r^2 - 2 b r^2}{4 b r^2} = \frac{1}{2}.
\end{equation}

This establishes the limiting claim in Theorem 3.

If $\nu = 1$, then the second term in the integral defining $M_n(r,\epsilon)$ is:
$$
[C_-,C_+] r^{2n} \log(v)|_{v_{**}}^{v_*(r)} = [C_-,C_+] r^{2n} \left(\log(r^{-(2 - \epsilon)}) - \log(v_{**}) \right) = \mathcal{O}(r^{2n} \log(r)) + \mathcal{O}(r^{2n}). 
$$

Therefore $J_b(r,c)$ is $\mathcal{O}(r^2 \log(r)) + \mathcal{O}(r^2)$ when $\nu = 1$. It follows that $J(r,c)$ is $\mathcal{O}(r^2 \log(r)) + \mathcal{O}(r^2)$ if $\nu = 1$ and is $\mathcal{O}(r^2) + \mathcal{O}(r^{2 \nu}) + \mathcal{O}(r^4)$ otherwise. Then equation \eqref{appeqn: integral function orders} holds, so the order statements in the Theorem hold. $\blacksquare$

\subsubsection{Proof of Corollary 2.1} \label{app: cor 2.1}

\setcounter{counter}{2}
\begin{snugshade}
\begin{corollary} \label{app cor 1: kernel sign}
    \textbf{[Sign of Kernel]} If $\mathcal{F} \sim \textnormal{GP}(\mu_f,k_f)$ for $\mu_f,k_f$ such that $\mathcal{F}(x,y) = -\mathcal{F}(y,x)$, and $\mathcal{F}(x,y) \neq \mu(x,y)$ then $k_f$ must take on both positive and negative values, and $k_f([x,y],[y,x])$ is negative where $\mathcal{F}(x,y)$ is nondeterministic.
\end{corollary}
\end{snugshade}

\textbf{Proof: } If $\mathcal{F}$ is random (non-deterministic), then $\mathbb{V}_{\mathcal{F}}[\mathcal{F}(x,y)] \neq 0$ for some $x,y$. Then $k_f([x,y],[x,y]) = \mathbb{V}_{\mathcal{F}}[\mathcal{F}(x,y)] \neq 0$ for some $x,y$. Moreover, variances are nonnegative, where the variance is nonzero, $k_f([x,y],[x,y]) > 0$. Then, by Theorem 2, $k_f([x,y],[y,x]) = - k_f([x,y],[x,y]) < 0$ for all $x,y$ such that $\mathcal{F}$ is non-deterministic. $\blacksquare$
\vspace{0.1 in}

\subsection{Supplementary Derivations}

\subsubsection{$\sigma^2$ and $\rho$ under Assumptions (1 - 2)} \label{app: 1 and 2}

Assume 1 and 2. Then:
$$
\sigma^2 = \mathbb{E}_{\mathcal{F}}[\mathbb{V}_{X,Y}[\mathcal{F}(X,Y)]].
$$

Since $\mathbb{E}_{X,Y \sim \pi_X}[f(X,Y)] = 0$ for any skew-symmetric $f$, we can rewrite the variance, then reverse the order of the expectations:
$$
\sigma^2 = \mathbb{E}_{\mathcal{F}}[\mathbb{E}_{X,Y}[\mathcal{F}(X,Y)^2]] = \mathbb{E}_{\mathcal{F},X,Y}[\mathcal{F}(X,Y)^2] = \mathbb{E}_{X,Y}[\mathbb{E}_{\mathcal{F}}[\mathcal{F}(X,Y)^2]]
$$

Finally, since we assumed that $\mathcal{F}$ is mean zero and drawn from a GP:
$$
\sigma^2 = \mathbb{E}_{X,Y}[\mathbb{V}[\mathcal{F}(X,Y)]] = \mathbb{E}_{X,Y}[\kappa_{F}([X,Y],[X,Y])].
$$

So, the intensity of competition (average squared payout) is:
\begin{equation}
    \sigma^2 = \mathbb{E}_{X,Y}[\kappa_{F}([X,Y],[X,Y])]
\end{equation}
where $X$ and $Y$ are drawn i.i.d.~from $\pi_x$. 

The shared-endpoint correlation follows via similar analysis:
$$
\sigma^2 \rho = \mathbb{E}_{\mathcal{F}}[\mathbb{C}_{X,Y,W}[\mathcal{F}(X,Y),\mathcal{F}(X,W)]].
$$

Since $f(X,Y)$ is mean zero for all skew-symmetric $f$ over choices of $X$ and $Y$, the covariance equals an expected product. Then we can reverse the order of the expectations:
$$
\sigma^2 \rho = \mathbb{E}_{\mathcal{F}}[\mathbb{E}_{X,Y,W}[\mathcal{F}(X,Y) \mathcal{F}(X,W)]] = \mathbb{E}_{\mathcal{F},X,Y,W}[\mathcal{F}(X,Y) \mathcal{F}(X,W)] = \mathbb{E}_{X,Y,W}[\mathbb{E}_{\mathbb{F}}[\mathcal{F}(X,Y) \mathcal{F}(X,W)]].
$$

Finally, since we assumed that $\mathcal{F}$ is mean zero, and drawn from a GP:
$$
\sigma^2 \rho = \mathbb{E}_{X,Y,W}[\mathbb{C}_{\mathbb{F}}[\mathcal{F}(X,Y), \mathcal{F}(X,W)]] = \mathbb{E}_{X,Y,W}[\kappa_f([X,Y],[X,W])].
$$

\subsubsection{$\rho$ when Squared Exponential} \label{app: rho SE}

\noindent \textbf{The General (Anisotropic) Case}

We will derive the general result first, then the isotropic result as a special case. 

For an anisotropic SE kernel set:
\begin{equation}
    h_{SE}(\phi|\Sigma_u) = \exp\left( - \frac{1}{2} \phi^{\intercal} \Sigma_u^{-1} s\phi\right).
\end{equation}

Now, notice that all of the expectations required for computing $\sigma$ and $\rho$ take the form:
\begin{equation}
    \mathbb{E}_{Z\sim \mathcal{N}(0,C)}\left[\exp\left( - \frac{1}{2} Z^{\intercal} D^{-1} Z \right) \right] = \det\left(I + D^{-1} C \right)^{-1/2}
\end{equation}
for some pair of covariance matrices $C$ and $D$. 

Here, the integral is computed by expanding the argument of the integral as a product of squared exponential functions, completing the square to express the argument as a un-normalized Gaussian density, then computing the corresponding normalization factor using the familiar form for the normalizing factor of a multivariate Gaussian. Then, by substitution into result 5.1:
\begin{equation} \label{appeqn: rho anisotropic}
    \rho(R) = \frac{\det\left(I + 2R \right)^{-1/2} - \det\left(I + R  \right)^{-1/2} \det\left(I + 3 R  \right)^{-1/2}}{1 - \det\left(I + 4 R \right)^{-1/2}}
\end{equation}

The second form provided in result 8 follows by expanding a determinant into a product of eigenvalues, and recalling that $\lambda(I + k A) = 1 + k \lambda(A)$ for any square matrix $A$ and scalar $k$. This yields:
\begin{equation} \label{appeqn: rho anisotropic roughness values}
        \rho(r) = \frac{\prod_{j=1}^T (1 + 2 r_j^2)^{-1/2} - \prod_{j=1}^T (1 + 2 r_j^2)^{-1/2} (1 + 3 r_j^2)^{-1/2}}{1 - \prod_{j=1}^T (1 + 4 r_j^2)^{-1/2}}
    \end{equation}
when $r_j = \sqrt{\lambda_j(R)}$. $\blacksquare$

\vspace{0.075 in}
\noindent \textbf{The Isotropic Case} \label{app: rho iso SE}
\vspace{0.05 in}

If isotropic, then $r_j = r$ for all $j$, so:
\begin{equation} \label{appeqn: rho isotropic SE gaussian traits}
\rho = \frac{\left(1 + 2 r^2 \right)^{-T/2} - (1 + r^2)^{-T/2} (1 + 3 r^2)^{-T/2}}{1 - (1 + 4 r^2)^{-T/2}}.
\end{equation}

Thus, the isotropic case follows trivially from the anisotropic case. $\blacksquare$

\vspace{0.1 in}
\noindent \textbf{Smoothness and Roughness Limits When Isotropic Squared Exponential} \label{app: smoothness and roughness iso SE}

\setcounter{result}{7}
\begin{snugshade}
\begin{limit_result} \label{app result 7.1: smoothness}
\textbf{[Smoothness Limit of $\rho$ under (\hyperref[table: Results + Assumptions]{1 - 7})]}  When the attributes are drawn independently of $\mathcal{F}$, the flow is mean zero, is drawn from a difference of stationary, isotropic, GP utilities with SE kernel with length scales $\Sigma_u$, and the attributes are normally distributed with covariance $\Sigma_x = \sigma_x^2 I$, in the limit as $r \rightarrow 0$:
\begin{equation} \label{appeqn: isotropic SE smoothness}
    \rho(r) \simeq \frac{1}{2} \left[1 - \frac{(T-1)}{2} r^2 + \mathcal{O}(r^4) \right] \xrightarrow{r \rightarrow 0} \frac{1}{2}.
\end{equation}
\end{limit_result}
\end{snugshade}

\textbf{Derivation:} All of the terms defining $\rho$ take the form:
$$(1 + k r^2)^{-T/2}$$
for some $k \in \{1,2,3,4\}$. When $\rho$ is small:
$$
(1 + k r^2)^{-T/2} \simeq 1 - \frac{k T}{2} r + \frac{k^2 T(T+2)}{8} r^2 + \mathcal{O}(r^3). 
$$

Then, substitution into $\rho$, yields:
$$
\begin{aligned}
    \rho(r) = \frac{( 1+ 2 r^2)^{-T/2} - (1 + r^2)^{-T/2} (1 + 3 r^2)^{-T/2}}{1 - (1 + 4 r^2)^{-T/2}} \simeq \frac{1}{2} - \frac{(T-1)}{4} r^2 + \mathcal{O}(r^4).
\end{aligned}
$$

The error in the approximation is $
\mathcal{O}(r^4)$ not $\mathcal{O}(r^3)$ since $\rho(r)$ is a composition of a function with $r^2$, so must be an even function of $r$ about zero. $\blacksquare$

\begin{snugshade}
\begin{limit_result} \label{app result 7: roughness}
\textbf{[Roughness Limit of $\rho$ under (\hyperref[table: Results + Assumptions]{1 - 7})]}  Under the same assumptions as stated in \ref{app result 7.1: smoothness}, in the limit as $r \rightarrow \infty$:
\begin{equation}  \label{appeqn: rho roughness limit isotropic SE}
    \rho(r) \simeq (\sqrt{2} r)^{-T} + \mathcal{O}(r^{-2T}) \xrightarrow{r \rightarrow \infty} 0.
\end{equation}

\end{limit_result}
\end{snugshade}

\textbf{Derivation:} In the rough limit $r \rightarrow \infty$. Then:
$$
\begin{aligned}
\lim_{r \rightarrow \infty} r^T \rho & = \lim_{r \rightarrow \infty} r^T \left(1 + 2 r^2 \right)^{-T/2} \left[ \frac{1 - (1 + 2 r^2)^{T/2}(1 + r^2)^{-T/2} (1 + 3 r^2)^{-T/2}}{1 - (1 + 4 r^2)^{-T/2}} \right] \\
& = 2^{-T/2} \lim_{r \rightarrow \infty} \left[ \frac{1 - \left(\frac{1 + 4 r^2 + 3 r^4}{1 + 2 r^2} \right)^{-T/2}}{1 - (1 + 4 r^2)^{-T/2}} \right] = 2^{-T/2} \lim_{r \rightarrow \infty} \left[ \frac{1 - \left(\frac{3}{2}\right)^{-T/2} r^{-T}}{1 - (4)^{-T/2} r^{-T}} \right] \\
& = 2^{-T/2}.
\end{aligned}
$$
Therefore:
\begin{equation}
    \rho(r) \simeq (\sqrt{2} r)^{-T} \left[1 - \left((3/2)^{-T/2} - 4^{-T/2} \right) r^{-T} \right] = (\sqrt{2} r)^{-T} + \mathcal{O}(r^{-2T}). \quad \blacksquare
\end{equation}

\pagebreak
\noindent \textbf{Smoothness and Roughness Limits When Anisotropic Squared Exponential} \label{app: smoothness and roughness aniso SE}

\setcounter{result}{8}
\setcounter{limit_result}{0}

\begin{snugshade}
\begin{limit_result} \label{app result 8: smoothness}
    \textbf{[Smoothness Limit of $\rho$, under (\hyperref[table: Results + Assumptions]{1 - 4, 6 - 7})]}  When the attributes are drawn independently of $\mathcal{F}$, the flow is mean zero, is drawn from a difference of stationary, GP utilities with SE kernel with length scales $\Sigma_u$, and the attributes are normally distributed with covariance $\Sigma_x$, in the limit as $r \rightarrow 0$:
    \begin{equation} \label{appeqn: half less rho smooth}
        \frac{1}{2} - \rho(R) \simeq \frac{1}{4} \left( \text{tr}(R) - \frac{\text{tr}(R^2)}{\text{tr}(R)}  \right) + \mathcal{O}(R^2) = \frac{1}{4}((T - 1) \text{av}(r^2) - \text{dis}(r^2)) + \mathcal{O}(r^4)
    \end{equation}
    where $T$ is the dimension of the trait space, $\text{av}$ denotes an arithmetic average and $\text{dis}$ denotes the index of dispersion (ratio of variance to average).
\end{limit_result}
\end{snugshade}

\textbf{Derivation:} Suppose that all of the roughness coefficients $r_j$ are small ($r_j \ll 1/4$ for all $j \in \{1,2,\hdots,T\}$). Then, expand the determinants in $\rho$:
$$
\prod_{j=1}^T (1 + k r_j^2) = 1 + k \sum_{j=1}^T r_j^2 + k^2 \sum_{i < j} r_i^2 r_j^2 + \hdots
$$

Formally the product expands into one plus a sum over each roughness value squared, then over all possible pairwise products squared, then all triplewise products, and so on. Since we've assumed that $k r_j^2 \ll 1$ for all $j$, we will truncate after the sum over pairs. This will allow us to formally work out the remainder term controlling convergence to hierarchy in the smoothness limit. Note that we only derived this term up to rate of convergence before.

Now:
$$
(1 + x)^{-1/2} = 1 - \frac{1}{2} x + \frac{3}{8} x^2 + \mathcal{O}(x^3).
$$

Therefore:
$$
\begin{aligned}
    \prod_{j=1}^T (1 + k r_j^2)^{-1/2} & = 1 - \frac{1}{2} \left(k \sum_{j=1}^T r_j^2 + k^2 \sum_{i < j} r_i^2 r_j^2 \right) + \frac{3}{8} k^2 \left(\sum_{j=1}^T r_j^2\right)^2 + \mathcal{O}(k^3 r^6) \\
    & = 1 - \frac{1}{2} k \sum_{j=1}^T r_j^2 + k^2 \left( \frac{3}{8} \left(\sum_{j=1}^T r_j^2\right)^2 - \frac{1}{2}  \sum_{i < j} r_i^2 r_j^2 \right) + \mathcal{O}(k^3 r^6) \\
    \end{aligned}
$$

To simplify note that $\sum_{i < j} r_i^2 r_j^2  = \frac{1}{2} \sum_{i \neq j} r_i^2 r_j^2 = \frac{1}{2} \sum_{i,j} r_i^2 r_j^2 - \frac{1}{2}\sum_{j} r_j^4  = \frac{1}{2} \left(\sum_{j} r_j^2 \right)^2 -  \frac{1}{2} \sum_{j} r_j^4$. Therefore:   

$$
\begin{aligned}
    \prod_{j=1}^T (1 + k r_j^2)^{-1/2} & = 1 - \frac{1}{2} k \sum_{j=1}^T r_j^2 + k^2 \left( \frac{3}{8} \left(\sum_{j=1}^T r_j^2\right)^2 - \frac{1}{2}  \sum_{i < j} r_i^2 r_j^2 \right) + \mathcal{O}(k^3 r^6)\\
    & 1 - \frac{1}{2} k \sum_{j=1}^T r_j^2 + k^2 \left( \frac{3}{8} \left(\sum_{j=1}^T r_j^2\right)^2 - \frac{1}{4} \left(  \sum_{j} r_j^2 \right) ^2 + \frac{1}{4} \sum_{j} r_j^4\right) + \mathcal{O}(k^3 r^6) \\
    & = 1 - \frac{1}{2} k \sum_{j=1}^T r_j^2 + k^2 \left( \frac{1}{8} \left(  \sum_{j} r_j^2 \right) ^2 + \frac{2}{8}  \sum_{j = 1}^T r_j^4 \right) + \mathcal{O}(k^3 r^6) 
    \end{aligned}
$$

Then:
$$
    \prod_{j=1}^T (1 + k r_j^2)^{-1/2} = 1 - \frac{1}{2} k \sum_{j=1}^T r_j^2 + \frac{1}{8} k^2 \left( \left(\sum_{j} r_j^2\right)^2 + 2  \sum_{j = 1}^T r_j^4 \right) + \mathcal{O}(k^3 r^6) 
$$

Equivalently, since $R$ has eigenvalues $r_j^2$:
\begin{equation}
    \prod_{j=1}^T (1 + k r_j^2)^{-1/2} = 1 - \frac{1}{2} k \text{tr}(R) + \frac{1}{8} k^2 \left( \text{tr}(R)^2 + 2  \text{tr}(R^2) \right) + \mathcal{O}(k^3 r^6) 
\end{equation}

Then, working a term at a time through $\rho$:
$$
1 - \prod_{j=1}^T (1 + 4 r_j^2)^{-1/2} = 2 \sum_{j=1}^T r_j^2 - 2 \left( \sum_{i,j} r_i^2 r_j^2 + 2  \sum_{j = 1}^T r_j^4 \right) + \mathcal{O}(r^6) 
$$
and:
$$
\begin{aligned}
& \textbf{(a) }\prod_{j=1}^T(1 + 2 r_j^2)^{-1/2}  = 1 - \sum_{j=1}^T r_j^2 + \frac{1}{2} \left( \sum_{i,j} r_i^2 r_j^2 + 2  \sum_{j = 1}^T r_j^4 \right) \\
& \textbf{(b) }\prod_{j=1}^T(1 + r_j^2)^{-1/2} (1 + 3 r_j^2)^{-1/2} = \left( 1 - \frac{1}{2} \sum_{j=1}^T r_j^2 + \frac{1}{8}\left( \sum_{i,j} r_i^2 r_j^2 + 2  \sum_{j = 1}^T r_j^4 \right) \right) \left( 1 - \frac{3}{2} \sum_{j=1}^T r_j^2 + \frac{9}{8} \left( \sum_{i,j} r_i^2 r_j^2 + 2  \sum_{j = 1}^T r_j^4 \right) \right) \hdots 
\end{aligned}
$$
where the errors in each approximation is $\mathcal{O}(r^6)$.

Simplifying the product:
$$
\begin{aligned}
\prod_{j=1}^T(1 + r_j^2)^{-1/2} (1 + 3 r_j^2)^{-1/2} & = 1 - 2 \sum_{j=1}^T r_j^2 + \frac{3}{4} \sum_{i,j} r_i^2 r_j^2 + \frac{10}{8} \left(\sum_{i,j} r_i^2 r_j^2 + 2  \sum_{j = 1}^T r_j^4 \right) + \mathcal{O}(r^6) \\
& = 1 - 2 \sum_{j=1}^T r_j^2 + 2 \sum_{i,j} r_i^2 r_j^2 + \frac{5}{2} \sum_{j=1}^T r_j^4 + \mathcal{O}(r^6)
\end{aligned}
$$

Then, the numerator is:
$$
\textbf{(a)} - \textbf{(b)} = \sum_{j=1}^T r_j^2 - \frac{3}{2} \sum_{i,j} r_i^2 r_j^2 - \frac{3}{2} \sum_{j} r_j^4 + \mathcal{O}(r^6).
$$

Therefore:
$$
\rho = \frac{\sum_{j=1}^T r_j^2 - \frac{3}{2} \sum_{i,j} r_i^2 r_j^2 - \frac{3}{2} \sum_{j} r_j^4 + \mathcal{O}(r^6)}{2 \sum_{j=1}^T r_j^2 - 2 \left( \sum_{i,j} r_i^2 r_j^2 + 2  \sum_{j = 1}^T r_j^4 \right) + \mathcal{O}(r^6) }
$$

Dividing the top and bottom by $\text{tr}(R)$ gives:
$$
\rho = \frac{1 - \frac{3}{2} \frac{\sum_{i,j} r_i^2 r_j^2}{\sum_{j=1}^T r_j^2} - \frac{3}{2} \frac{\sum_{j} r_j^4}{\sum_{j=1}^T r_j^2} + \mathcal{O}(r^6)}{2  - 2  \frac{\sum_{i,j} r_i^2 r_j^2}{\sum_{j=1}^T r_j^2} - 4  \frac{\sum_{j = 1}^T r_j^4}{\sum_{j=1}^T r_j^2} + \mathcal{O}(r^6) }
$$

Let $\epsilon_1(r) = (\sum_{i,j} r_i^2 r_j^2)/(\sum_{j=1}^T r_j^2)$ and $\epsilon_2(r) = \frac{\sum_{j=1}^T r_j^4}{\sum_{j=1}^T r_j^2}$. Then $\epsilon_1 = \sum_{j=1}^T r_j^2$. Both $\epsilon_1$ and $\epsilon_2$ are $\mathcal{O}(r^4)$ so will dominate the $\mathcal{O}(r^6)$ terms. Then:
$$
\rho \simeq \frac{1 - \frac{3}{2} \epsilon_1- \frac{3}{2} \epsilon_2}{2 - 2 \epsilon_1 - 4 \epsilon_2} \simeq \frac{1}{2} \left( {1 - \frac{3}{2} \epsilon_1} - \frac{3}{2} \epsilon_2 \right)\left(1 + \epsilon_1 + 2 \epsilon_2 \right) \simeq \frac{1}{2} \left( 1 - \frac{1}{2} \epsilon_1 + \frac{1}{2} \epsilon_2 \right) + \mathcal{O}(r^4)
$$

So, substituting back in for $\epsilon_1$ and $\epsilon_2$, when the roughness values are all small:
\begin{equation}
    \rho(R) = \frac{1}{2} - \frac{1}{4} \left( \text{tr}(R) -  \frac{\text{tr}(R^2)}{\text{tr}(R)}  \right) + \mathcal{O}(R^2).
\end{equation}

Alternately, let $\text{av}(r^2)$ be the arithmetic average of the roughness values (squared), and $\text{dis}(r^2)$ be the coefficient of dispersion in the roughness values squared. The coefficient of dispersion is the ratio of the variance to the mean so:
$$
\begin{aligned}
& \text{tr}(R) = \sum_{j=1}^T r_j^2 = T \text{av}(r^2) \\
& \frac{\text{tr}(R^2)}{\text{tr}(R)} = \frac{\sum_{j=1}^T r_j^4}{\sum_{j=1}^T r_j^2} = \frac{\sum_{j=1}^T r_j^2}{T \text{av}(r^2)} = \frac{\text{var}(r^2) + \text{av}(r^2)^2}{\text{av}(r^2)} = \text{av}(r^2) + \text{dis}(r^2).
\end{aligned}
$$

Therefore:
\begin{equation}
 \rho(R) = \frac{1}{2} - \frac{1}{4} \left( \text{tr}(R) -  \frac{\text{tr}(R^2)}{\text{tr}(R)}  \right) + \mathcal{O}(R^2) = \frac{1}{2} - \frac{1}{4}((T - 1) \text{av}(r^2) - \text{dis}(r^2)) + \mathcal{O}(r^4) \quad \blacksquare
\end{equation}

\begin{snugshade}
\begin{limit_result} \label{app result 8: roughness}
    \textbf{[Roughness Limit of $\rho$, under (\hyperref[table: Results + Assumptions]{1 - 4, 6 - 7})]}  Under the same assumptions as stated in \ref{app result 8: smoothness}, in the limit as $r \rightarrow \infty$:
    \begin{equation} \label{appeqn: rho rough smooth approx}
    \rho(R) \approx 2^{-|\mathcal{R}|/2} \bar{r}_{\mathcal{R}}^{-|\mathcal{R}|} \left[1 - \left((3/2)^{-|\mathcal{R}|/2} - 4^{-|\mathcal{R}|/2} \right) \bar{r}_{\mathcal{R}}^{-|\mathcal{R}|}\right] = \mathcal{O}(\bar{r}_{\mathcal{R}}^{-|\mathcal{R}|})\xrightarrow{\bar{r}_{\mathcal{R}} \rightarrow \infty} 0.
\end{equation}
\end{limit_result}
\end{snugshade}

\textbf{Derivation:} First, notice that all the terms in $\rho$ are soft-thresholded versions of the roughness parameters, bounded below by:
\begin{equation}
    (1 + k r^2) \geq \min\{k r^2, 1\}.
\end{equation}

When $r_j \ll 1$, the corresponding terms are all approximately 1. When $r_j^2 \gg 1$, then all terms are approximately $k r^2$. This suggests an initial approximation. If the roughness values span many orders of magnitude, and $T$ is small relative to the orders of magnitude spanned, then the range of possible values between $\max_j\{r_j\}$ and $\min_j\{r_j\}$ may be sparsely sampled, so there may not be any roughness values such that $k r_j^2$ is on the order of 1 for $k = 1,2,3,4$. Then, we may be able to approximate the shared-endpoint correlation by separating the roughness values into a rough set and smooth set:
\begin{equation}
\begin{aligned}
    & \textbf{Rough Set: } \mathcal{R} = \{j \mid k r_j^2 \gg 1\ \forall k \in \{1,2,3,4\}\}\\
    & \textbf{Smooth Set: } \mathcal{S} = \{j \mid k r_j^2 \ll 1\ \forall k \in \{1,2,3,4\}\}\\
\end{aligned}
\end{equation}

If $\mathcal{R} \cup \mathcal{S}$ covers all $j \in [1,T]$, then the set of roughness values can be partitioned into the rough and smooth sets. Let $|\mathcal{S}| = |\mathcal{S}|$ and $|\mathcal{R}| = |\mathcal{R}|$ denote the cardinality of these sets. Then $|\mathcal{S}| + |\mathcal{R}| = T$. 

Now, assuming such a partition exists, it is reasonable to adopt the approximation:
\begin{equation}
    (1 + k r_j^2) \simeq \begin{cases} & 1 \text{ if } j \in \mathcal{S} \\
    & k r_j^2 \text{ if } j \in \mathcal{R}
    \end{cases}
\end{equation}

Under this approximation, the shared-endpoint correlation simplifies dramatically since the smooth set drops out of the calculation, and the powers can be distributed.
$$
\begin{aligned}
\rho & \simeq \frac{\prod_{j\in \mathcal{R}} (2 r_j^2)^{-1/2} - \prod_{j\in \mathcal{R}}^T (3 r_j^4)^{-1/2}}{1 - \prod_{j\in \mathcal{R}} (4 r_j^2)^{-1/2}} = \frac{\prod_{j\in \mathcal{R}} (2 r_j^2)^{-1/2} - \prod_{j\in \mathcal{R}}^T (3 r_j^4)^{-1/2}}{1 - \prod_{j\in \mathcal{R}} (4 r_j^2)^{-1/2}}  = \frac{\prod_{j\in \mathcal{R}} (\sqrt{2} r_j)^{-1} - \prod_{j\in \mathcal{R}}^T (\sqrt{3} r_j^2)^{-1}}{1 - \prod_{j\in \mathcal{R}} (2 r_j)^{-1}} \\
& = \frac{1}{2^{|\mathcal{R}|/2} \prod_{j\in \mathcal{R}} r_j} \left[ \frac{1 - 3^{-|\mathcal{R}|/2} \prod_{j \in \mathcal{R}} r_j^{-2} 2^{|\mathcal{R}|/2} \prod_{j \in \mathcal{R}} r_j}{1 - (2^|\mathcal{R}| \prod_{j \in \mathcal{R}} r_j)^{-1}} \right] = \frac{1}{2^{|\mathcal{R}|/2} \prod_{j\in \mathcal{R}} r_j} \left[ \frac{1 - ((3/2)^{|\mathcal{R}|/2} \prod_{j \in \mathcal{R}} r_j)^{-1} }{1 - (2^|\mathcal{R}| \prod_{j \in \mathcal{R}} r_j)^{-1}} \right]
\end{aligned}
$$

Now, let:
\begin{equation}
    \bar{r}_{\mathcal{R}} = \left(\prod_{j \in \mathcal{R}} r_j \right)^{1/|\mathcal{R}|}
\end{equation}
denote the geometric average of the roughness values in the rough set (from here on out, the rough values). Then, provided we can partition the roughness values into a rough and smooth set:
\begin{equation}
    \rho \simeq 2^{-|\mathcal{R}|/2} \bar{r}_{\mathcal{R}}^{-|\mathcal{R}|} \left[ \frac{1 - (3/2)^{-|\mathcal{R}|/2} \bar{r}_{\mathcal{R}}^{-|\mathcal{R}|} }{1 - (4)^{-|\mathcal{R}|/2} \bar{r}_{\mathcal{R}}^{-|\mathcal{R}|}} \right] \longrightarrow \lim_{\bar{r}_{\mathcal{R}} \rightarrow \infty} \rho \simeq (\sqrt{2} \bar{r}_{\mathcal{R}})^{-|\mathcal{R}|} \quad \blacksquare
\end{equation}

\subsubsection{Smoothness Limits of $\rho$ for Mixture Models} \label{app:  smoothness limit mixture}

Suppose that the distribution of roughness coefficients concentrates at small $r$. Then, provided the tails of the $r^2$ distribution vanish sufficiently fast:

\setcounter{result}{9}
\setcounter{limit_result}{0}
\begin{snugshade}
\begin{limit_result}
    \textbf{[Smoothness Limit of $\rho$ under (\hyperref[table: Results + Assumptions]{1 - 4, 8})]} When the attributes are drawn independently of $\mathcal{F}$, the flow is mean zero, and is drawn from a difference of stationary, GP utilities with a mixture of SE kernels, in the limit as $r \rightarrow 0$:
    \begin{equation} \label{appeqn: smoothness limit of mixture }
        \frac{1}{2} - \rho \approx \frac{1}{4} \frac{\mathbb{E}_{r \sim \pi_r}[\text{tr}(R)^2 - \text{tr}(R^2)]}{\mathbb{E}_{r \sim \pi_r}[\text{tr}(R)]} = \begin{cases} & \frac{1}{4} \frac{\sum_{i \neq j} \mathbb{E}_{r \sim \pi_r}\left[ r_i^2 r_j^2 \right]}{\sum_{j} 
 \mathbb{E}_{r \sim \pi_r}\left[r_j^2 \right]}\\
        & \frac{1}{4}\left((T-1) \text{av}(\hat{r}^2) - \text{dis}(\hat{r}^2) + (T-1)\frac{\frac{1}{T (T-1)} \sum_{i \neq j} \mathbb{C}_{r \sim \pi_r}[r_i^2,r_j^2]}{\text{av}(\hat{r}^2) } \right) \\
        & \frac{1}{4} \left((T-1) \frac{\mathbb{E}_{r \sim \pi_r}[\mathbb{E}_{J}[r_J^2]^2]}{\mathbb{E}_{J,r \sim \pi_r}[r_J^2]} - \frac{\mathbb{E}_{r \sim \pi_r}[\mathbb{V}_{J}[r_J^2]]}{\mathbb{E}_{J,r \sim \pi_r}[r_J^2]}\right)
        \end{cases}
    \end{equation}
    where the three forms on the right hand side of Equation \eqref{appeqn: smoothness limit of mixture } are equivalent, $\text{tr}$ denotes the trace, $\hat{r}^2 = \mathbb{E}_{r \sim \pi_r}[r^2]$, $\text{av}$ and $\text{dis}$ denote the arithmetic average and coefficient of dispersion over the index $j \in \{1,2,\hdots, T\}$, $\mathbb{C}$ denotes covariance, $J$ denotes a uniformly drawn index, and $\mathbb{V}$ denotes variance. 
\end{limit_result}
\end{snugshade}

Consider a mixture model with $R \sim \pi_r$. Then:

    \begin{equation} \label{appeqn: rho mixed SE kernel}
        \rho = \frac{\mathbb{E}_{R \sim \pi_R}\left[\det\left(I + 2R \right)^{-1/2}\right] - \mathbb{E}_{R \sim \pi_R}\left[\det\left(I + R  \right)^{-1/2} \det\left(I + 3 R  \right)^{-1/2} \right]}{1 - \mathbb{E}_{R \sim \pi_R}\left[\det\left(I + 4 R \right)^{-1/2} \right]}
    \end{equation}

Notice the placement of the expectations in equation \eqref{appeqn: rho mixed SE kernel}. These appear outside all terms that involved an average of the kernel over samples $X$, $Y$, and $W$, but cannot be taken onto individual terms (as in the product in the numerator), or to the outside of the entire expression. 

To simplify in the smoothness limit, expand the each term inside the expectations assuming small $R$:
$$
\begin{aligned}
& \textbf{Generic: } \text{det}(I + k R)^{-1/2} = 1 - \frac{1}{2} k \text{tr}(R) + \frac{1}{8} k^2 \left( \text{tr}(R)^2 + 2  \text{tr}(R^2) \right) + \mathcal{O}(R^3) \\
& \textbf{(a) } \text{det}(I + 2 R)^{-1/2} = 1 - \text{tr}(R) + \frac{1}{2} \left( \text{tr}(R)^2 + 2  \text{tr}(R^2) \right) + \mathcal{O}(R^3) \\
& \textbf{(b) } \det\left(I + R  \right)^{-1/2} \det\left(I + 3 R  \right)^{-1/2} = 1 - 2 \text{tr}(R) + 2 \text{tr}(R)^2 + \frac{5}{2} \text{tr}(R^2) + \mathcal{O}(R^3)  \\
& \textbf{(c) } \text{det}(I + 4 R)^{-1/2} = 1 - 2 \text{tr}(R) + 2 \left( \text{tr}(R)^2 + 2  \text{tr}(R^2) \right) + \mathcal{O}(R^3) \\
\end{aligned}
$$

Then:
$$
\begin{aligned}
\rho & = \frac{\mathbb{E}_{R\sim\pi_R}[\textbf{(a)}] - \mathbb{E}_{R\sim\pi_R}[\textbf{(b)}]}{1 - \mathbb{E}_{R\sim\pi_R}[\textbf{(c)}]} = \frac{ \mathbb{E}_{R\sim\pi_R}[\text{tr}(R)] - \frac{3}{2}\mathbb{E}_{R\sim\pi_R}[\text{tr}(R)^2] - \frac{3}{2} \mathbb{E}_{R\sim\pi_R}[\text{tr}(R^2)] + \mathcal{O}(R^3)}{2 \mathbb{E}_{R\sim\pi_R}[\text{tr}(R)] - 2 \mathbb{E}_{R\sim\pi_R}[\text{tr}(R)^2] - 4 \mathbb{E}_{R\sim\pi_R}[\text{tr}(R^2)] + \mathcal{O}(R^3)} \\
& \simeq \frac{1}{2} \frac{1 - \frac{3}{2} \epsilon_1 - \frac{3}{2} \epsilon_1}{1 - 2 \epsilon_1 - 4 \epsilon_2} \simeq \frac{1}{2} \left(1 - \frac{3}{2} \epsilon_1 - \frac{3}{2} \epsilon_2 \right)\left(1 + \epsilon_1 + 2 \epsilon_2 \right) 
\simeq \frac{1}{2} \left(1 - \frac{1}{2} \epsilon_1 + \frac{1}{2} \epsilon_2\right) \\ & = \frac{1}{2} - \frac{1}{4} \frac{\mathbb{E}_{R \sim \pi_R}[\text{tr}(R)^2] - \mathbb{E}_{R \sim \pi_R}[\text{tr}(R^2)]}{\mathbb{E}_{R \sim \pi_R}[\text{tr}(R)]}
\end{aligned}
$$

Thus, when $\pi_R$ is chosen to put most of its weight to small $R$ ($k R \ll I$), then:
\begin{equation} \label{appeqn: rho smooth mixture}
\begin{aligned}
    \rho & \approx \frac{1}{2} - \frac{1}{4} \frac{\mathbb{E}_{R \sim \pi_R}[\text{tr}(R)^2] - \mathbb{E}_{R \sim \pi_R}[\text{tr}(R^2)]}{\mathbb{E}_{R \sim \pi_R}[\text{tr}(R)]} = \frac{1}{2} - \frac{1}{4} \frac{\mathbb{E}_{r \sim \pi_R}\left[\sum_{i \neq j} r_i^2 r_j^2 \right]}{\mathbb{E}_{r \sim \pi_R}\left[\sum_j r_j^2 \right]} = \frac{1}{2} \left(1 - \frac{\mathbb{E}_{r \sim \pi_R}\left[ \sum_{i < j} r_i^2 r_j^2 \right]}{\mathbb{E}_{r \sim \pi_R}\left[ \sum_{j} r_j^2 \right]} \right)
\end{aligned}
\end{equation}

To simplify in terms of the roughness values, note that:
$$
\begin{aligned}
\mathbb{E}_{r \sim \pi_R}\left[\sum_{i \neq j} r_i^2 r_j^2 \right] & = \sum_{i,j} \mathbb{E}_{r \sim \pi_R}\left[ r_i^2 r_j^2 \right] - \sum_j \mathbb{E}_{r \sim \pi_R}\left[ r_j^2 \right] \\
& = \sum_{i,j} \mathbb{E}_{r \sim \pi_R}\left[ r_i^2 \right] \mathbb{E}_{r \sim \pi_R}\left[ r_j^2 \right] + \sum_{i,j} \mathbb{C}_{r \sim \pi_R}[r_i^2,r_j^2] - \sum_{j} \mathbb{E}_{r \sim \pi_R}[r_j^2]^2 - \sum_{j} \mathbb{V}_{r \sim \pi_R}[r_j^2] \\
& \left( \mathbb{E}_{r \sim \pi_R} \left[\sum_j r_j^2 \right] \right)^2 - \sum_j \mathbb{E}_{r \sim \pi_R} \left[ r_j^2 \right]^2 + \sum_{i,j} \mathbb{C}_{r \sim \pi_R}[r_i^2,r_j^2] - \sum_{j} \mathbb{C}_{r \sim \pi_R}[r_j^2,r_j^2] \\
& = \left( \mathbb{E}_{r \sim \pi_R} \left[\sum_j r_j^2 \right] \right)^2 - \sum_j \mathbb{E}_{r \sim \pi_R} \left[ r_j^2 \right]^2 + \sum_{i \neq j} \mathbb{C}_{r \sim \pi_R}[r_i^2,r_j^2]
\end{aligned}
$$

Then:
$$
\frac{\mathbb{E}_{r \sim \pi_R}\left[ \sum_{i < j} r_i^2 r_j^2 \right]}{\mathbb{E}_{r \sim \pi_R}\left[ \sum_{j} r_j^2 \right]} = \mathbb{E}_{r \sim \pi_R}\left[ \sum_j r_j^2 \right] - \frac{\sum_j \mathbb{E}_{r \sim \pi_R} \left[ r_j^2 \right]^2}{\mathbb{E}_{r \sim \pi_R}\left[ \sum_j r_j^2 \right]} + \frac{\sum_{i \neq j} \mathbb{C}_{r \sim \pi_R}[r_i^2,r_j^2]}{\mathbb{E}_{r \sim \pi_R}\left[ \sum_j r_j^2 \right]}
$$

Let $\hat{r}_j^2 = \mathbb{E}_{r \sim \pi_R}[r_j^2]$, let $\text{av}$ denote an average over dimensions, and $\text{var}$ denote variance over dimension. Then, the second term can be written:
$$
\sum_j \mathbb{E}_{r \sim \pi_R} \left[ r_j^2 \right]^2 = T (\text{av}(\hat{r}^2) + \text{var}(\hat{r}^2))
$$

Then:
$$
\mathbb{E}_{r \sim \pi_R}\left[ \sum_j r_j^2 \right] - \frac{\sum_j \mathbb{E}_{r \sim \pi_R} \left[ r_j^2 \right]^2}{\mathbb{E}_{r \sim \pi_R}\left[ \sum_j r_j^2 \right]} = T\text{av}(\hat{r}^2) - \frac{T (\text{av}(\hat{r}^2) + \text{var}(\hat{r}^2))}{T \text{av}(\hat{r}^2)} = (T - 1) \text{av}(\hat{r}^2) - \frac{\text{var}(\hat{r}^2)}{\text{av}(\hat{r})} = (T - 1) \text{av}(\hat{r}^2) - \text{dis}(\hat{r}^2).
$$

Therefore, the first two terms match the terms derived for a single SE exponential kernel, with roughness values set to the average roughness values over the mixture $\pi_R$. Then:
\begin{equation} \label{appeqn: rho total variance 1}
\rho \approx \frac{1}{2} - \frac{1}{4}\left((T-1) \text{av}(\hat{r}^2) - \text{dis}(\hat{r}^2) + \frac{1}{T}\frac{\sum_{i \neq j} \mathbb{C}_{R \sim \pi_R}[r_i^2,r_j^2]}{\text{av}(\hat{r}^2) } \right).
\end{equation}

Thus, when $\pi_R$ only assigns weight to small $R$, replacing the SE kernel with a mixture of SE kernels only changes the expansion of $\rho$ by an added codispersion term associated with the the covariance in the roughness values over different $R \sim \pi_R$. Notice that equation \eqref{appeqn: rho total variance 1} is, in essence, the law of total variance applied to $\rho$. The first term accounts for the variance over dimensions, in expectation over $\pi_R$, while the second accounts for the variance over $R$. Let $J \sim \text{Unif}[1,T]$ represent a randomly drawn trait. Then:
$$
\begin{aligned}
    & \text{av}(\hat{r}^2) = \mathbb{E}_{J,R}[r_j^2]\\
    & \text{var}(\hat{r}^2) = \mathbb{V}_{J}[\mathbb{E}_R[r_j^2]]\\
\end{aligned}
$$

Moreover,
$$
\frac{1}{T} \sum_{i \neq j} \mathbb{C}_{R\sim \pi_R}[r_i^2,r_j^2] = \frac{T}{T^2} \sum_{i,j} \mathbb{C}_{R\sim\pi_R}[r_i^2,r_j^2] - \frac{1}{T} \sum_{j} \mathbb{V}_{R \sim \pi_R}[r_j^2] = T \mathbb{V}_{R}[\mathbb{E}_{J}[r_j^2]] - \mathbb{E}_{J}[\mathbb{V}_{R}[r_j^2]].
$$

Then:
$$
\begin{aligned}
    \rho & \approx \frac{1}{2} - \frac{1}{4} \left((T-1) \mathbb{E}_{J,R}[r_j^2] +  T \frac{\mathbb{V}_{R}[\mathbb{E}_{J}[r_j^2]]}{\mathbb{E}_{J,R}[r_j^2]} - \frac{\mathbb{E}_{J}[\mathbb{V}_{R}[r_j^2]] + \mathbb{V}_{J}[\mathbb{E}_{R}[r_j^2]]}{\mathbb{E}_{J,R}[r_j^2]}\right) \\
    & = \frac{1}{2} - \frac{1}{4} \left((T-1) \mathbb{E}_{J,R}[r_j^2] +  \frac{T \mathbb{V}_{R}[\mathbb{E}_{J}[r_j^2]] -\mathbb{V}_{J,R}[r_j^2]}{\mathbb{E}_{J,R}[r_j^2]}\right) \\
    & = \frac{1}{2} - \frac{1}{4} \left((T-1) \mathbb{E}_{J,R}[r_j^2] +  \frac{T \mathbb{V}_{R}[\mathbb{E}_{J}[r_j^2]] -\mathbb{V}_{J,R}[r_j^2]}{\mathbb{E}_{J,R}[r_j^2]}\right) \\
    & = \frac{1}{2} - \frac{1}{4} \left((T-1) \mathbb{E}_{J,R}[r_j^2] +  \frac{T \mathbb{V}_{R}[\mathbb{E}_{J}[r_j^2]] -\mathbb{V}_{R}[\mathbb{E}_{J}[r_j^2]] - \mathbb{E}_{R}[\mathbb{V}_{J}[r_j^2]]}{\mathbb{E}_{J,R}[r_j^2]}\right) \\
    & = \frac{1}{2} - \frac{1}{4} \left((T-1) \mathbb{E}_{J,R}[r_j^2] +  (T-1) \frac{\mathbb{V}_{R}[\mathbb{E}_{J}[r_j^2]]}{\mathbb{E}_{J,R}[r_j^2]} - \frac{\mathbb{E}_{R}[\mathbb{V}_{J}[r_j^2]]}{\mathbb{E}_{J,R}[r_j^2]}\right)
\end{aligned}
$$

In short:
    \begin{equation} \label{appeqn: rho total variance 2}
    \begin{aligned}
        \rho & \approx \frac{1}{2} - \frac{1}{4} \left((T-1) \mathbb{E}_{J,R}[r_j^2] +  (T-1) \frac{\mathbb{V}_{R}[\mathbb{E}_{J}[r_j^2]]}{\mathbb{E}_{J,R}[r_j^2]} - \frac{\mathbb{E}_{R}[\mathbb{V}_{J}[r_j^2]]}{\mathbb{E}_{J,R}[r_j^2]}\right) \\
        & = \frac{1}{2} - \frac{1}{4} \left((T-1) \frac{\mathbb{E}_{R}[\mathbb{E}_{J}[r_j^2]^2]}{\mathbb{E}_{J,R}[r_j^2]} - \frac{\mathbb{E}_{R}[\mathbb{V}_{J}[r_j^2]]}{\mathbb{E}_{J,R}[r_j^2]}\right) \quad \blacksquare
        \end{aligned}
    \end{equation}

\noindent \textbf{Remark:}  Note the usual dependence on $T$. If $T = 1$, then all of the correction terms of order $r^2$ vanish. 

\noindent \textbf{Remark:} The three equivalent forms of equation \eqref{appeqn: smoothness limit of mixture } are useful. The first is the simplest. The second is useful for comparison to the smoothness limit when $\pi_x$ is Gaussian and $\kappa_u$ is SE. Comparing equation \eqref{appeqn: smoothness limit of mixture } to the smoothness limit in the SE case, we see that, varying the trait and kernel variances over mixture distributions only change the lowest order terms in the smoothness limit by the added term:
$$
\frac{(T-1)}{4}\frac{\frac{1}{T (T-1)} \sum_{i \neq j} \mathbb{C}_{r \sim \pi_r}[r_i^2,r_j^2]}{\text{av}(\hat{r}^2) }.
$$
The numerator is the average over all distinct pairs of roughness coefficients, of the covariance in those coefficients squared. This represents the variation in the roughness coefficients associated with the mixture distributions. Dividing by $\text{av}(\hat{r}^2)$ converts from a measure of variance, to a measure of dispersion. Then, the second term can be interpreted as the dispersion in the squared roughness coefficient contributed by the mixture distribution, as opposed to $\text{dis}(\hat{r}^2)$ which is the dispersion across indices of the expected squared roughness coefficient. 

All three forms presented in equation \eqref{appeqn: smoothness limit of mixture } are all related by a variance decomposition which separates dispersion among the indices (directions of principal roughness), and over the mixture distribution. This decomposition is made explicit in the third form. 

In the isotropic case, the dispersion over indices vanishes, so:

\setcounter{result}{10}
\setcounter{limit_result}{0}
\begin{snugshade}
\begin{limit_result} \label{app result 10: smoothness}
    \textbf{[Smoothness Limit of $\rho$ under (\hyperref[table: Results + Assumptions]{1 - 5, 8})]}  When the attributes are drawn independently of $\mathcal{F}$, the flow is mean zero, and is drawn from a difference of stationary, isotropic, GP utilities with a mixture of SE kernels, in the limit as $r \rightarrow 0$:
    \begin{equation} \label{appeqn: smoothness limit of iso mixture}
        \frac{1}{2} - \rho \approx \frac{(T-1)}{4} \frac{\mathbb{E}_{r \sim \pi_r}\left[ r^4 \right]}{ 
 \mathbb{E}_{r \sim \pi_r}\left[r^2 \right]} =  \frac{(T-1)}{4}\left(\hat{r}^2 + \frac{ \mathbb{V}_{r \sim \pi_r}[r^2]}{\hat{r}^2 } \right) 
    \end{equation}
    where $\hat{r}^2 = \mathbb{E}_{r \sim \pi_r}[r^2]$, and $\mathbb{V}$ denotes variance.
\end{limit_result}
\end{snugshade} 

Notice, that mixing over $r$, decreases $\rho$ by the coefficient of dispersion of $r^2$. Thus, when smooth ($\pi_r$ concentrated at small $r$), mixture models produce smaller $\rho$ than the corresponding SE model with Gaussian traits, if the roughness of the SE model is set to $\hat{r} = (\mathbb{E}_{r \sim \pi_r}[r^2])^{1/2}$. The degree to which mixing reduces the correlation is proportional to the coefficient of dispersion in $r^2$.

\subsubsection{$\rho$ when Mat\'ern} 

\noindent \textbf{Dimensionality 
Reduction}\label{app: dimensionality reduction}

\setcounter{equation_env}{38}
\begin{snugshade}
\begin{equation_env}
    \textbf{[Mat\'ern Kernel Dimensionality Reduction]} The Mat\'ern kernel in $d$ dimensions can be related to the Mat\'ern kernel in $d'$ dimensions by the correct change of parameters.
    \begin{equation} \label{appeqn: necessary expectations pre frequency space}
    \begin{aligned}
        & \textbf{a\hspace{0.035 in})}\indent  \mathbb{E}_{Z \sim \mathcal{N}(0,2 \sigma_x^2 I)}\left[k_{2T}([Z,0]; \nu_{2T}, l_{2T}) \right] = \mathbb{E}_{Z \sim \mathcal{N}(0,2 \sigma_x^2 I)}\left[k_{T}(Z; \nu_{T}, l_{T}) \right] \\
        & \textbf{a')}\indent  \mathbb{E}_{Z \sim \mathcal{N}(0,2 \sigma_x^2 I)}\left[k_{2T}([\sqrt{2}Z,0]; \nu_{2T}, l_{2T}) \right] = \mathbb{E}_{Z \sim \mathcal{N}(0,4 \sigma_x^2 I)}\left[k_{2T}([Z,0]; \nu_{2T}, l_{2T}) \right]  = \mathbb{E}_{Z \sim \mathcal{N}(0,4 \sigma_x^2 I)}\left[k_{T}(Z; \nu_{T}, l_{T}) \right] \\
    \end{aligned}
\end{equation}
Therefore, we can replace the expectations \textbf{a)} and \textbf{a')} to an average over a $T$ dimensional space, provided we derive the correct parameter change such that the two kernels are equivalent when evaluated at $[z,0]$ and at $z$. In fact, the Mat\'ern kernel is parameterized so that this exchange is trivial - all parameters remain the same:
    \begin{equation} \label{appeqn: kernel dimension relationship}
        \kappa_{2T}([Z,0]; \nu, l) = \kappa_{T}(Z; \nu, l) = \mathcal{F}^{-1}[\hat{\kappa}_{T}(\zeta ;\nu, l)](z)
    \end{equation}
\end{equation_env}

\end{snugshade}

To begin, we consider:
\begin{equation} \label{appeqn: 2T to T conversion}
    \kappa_{2T}([Z,0]; \nu_{2T}, l_{2T}) = \alpha\cdot\kappa_{T}(Z; \nu_T, l_T).
\end{equation}
Thus, we are tasked with finding $\alpha$, as well as the relationship between, $\nu_{2T} \text{ and } \nu_T, \text{ along with } l_{2T} \text{ and } l_T$. 

Mat\'ern kernels can be equivalently expressed in frequency space through a Fourier transform. Therefore, the technical task is reduced to matching the Fourier transforms of the $2T$ dimensional kernel with the $T$ dimensional kernel, at least up to evaluation at $[Z,0]$ for all $Z$:
    \begin{equation} \label{appeqn: 2T Fourier Transform}
        \kappa_{2T}([Z,0];\nu_{2T}, l_{2T}) = \mathcal{F}^{-1}[\hat{\kappa}_{2T}({\zeta} ; \nu_{2T}, l_{2T})]([z,0]) 
    \end{equation}
    \begin{equation} \label{appeqn: T Fourier Transform}
        \alpha\cdot\kappa_{T}(Z;\nu_{T}, l_{T}) = \alpha\cdot\mathcal{F}^{-1}[\hat{\kappa}_{T}(\zeta ; \nu_{T}, l_{T})](z) 
    \end{equation}
    \begin{equation} \label{appeqn: 2T to T Fourier Transform}
        \mathcal{F}^{-1}[\hat{\kappa}_{2T}({\zeta} ; s_{2T}, l_{2T})]([z,0]) = \alpha\cdot\mathcal{F}^{-1}[\hat{\kappa}_{T}(\zeta ; s_{T}, l_{T})](z) 
    \end{equation}
    for some scaling constant alpha.
   
The Matern kernel is simple in frequency space and has a simple power spectral density:
    \begin{equation} \label{appeqn: Matern kernel in frequency space}
        \kappa_{d}(\tau; \nu, l) \leftrightarrow \hat{\kappa}_{d}(\zeta; \nu, l) = \frac{\Gamma(\nu+\tfrac{d}{2})}{\Gamma(\nu)} 2^d \pi^{d/2} \left(\frac{l^2}{2\nu}\right)^{d/2} \left(1 + \frac{l^{2} \|2 \pi \zeta\|^2}{2\nu} \right)^{-\left(\nu+\tfrac{d}{2}\right)}.
    \end{equation}
Thus, the kernels in question are:
    \begin{equation} \label{appeqn: Kernel in 2T}
        \hat{\kappa}_{2T}({\zeta} ; \nu_{2T}, l_{2T}) = \frac{\Gamma(\nu_{2T}+\tfrac{2T}{2})}{\Gamma(\nu_{2T})}2^{2T}\pi^{2T/2} \left(\frac{l_{2T}^{2}}{2\nu_{2T}}\right)^{2T/2}\left(1+\frac{l_{2T}^2\cdot
        ||2\pi{\zeta}||^2}{2\nu_{2T}}\right)^{-\left(\nu_{2T}+\tfrac{2T}{2}\right)}
    \end{equation}
    and 
    \begin{equation}\label{appeqn: Kernel in T}
        \hat{\kappa}_{T}({\zeta} ; s_{T}, l_{T}) = \frac{\Gamma(\nu_T+\tfrac{T}{2})}{\Gamma(\nu_T)}2^{T}\pi^{T/2} \left(\frac{l_{T}^2}{2\nu_T}\right)^{T/2}\left(1+\frac{l_{T}^2\cdot 
        ||2\pi{\zeta}||^2}{2\nu_T}\right)^{-\left(\nu_{T}+\tfrac{T}{2}\right)}
    \end{equation}

If we evaluate the inverse Fourier transform of the power spectral density, $\hat{\kappa}_{2T}$ at $[z,0]$ we see:
    \begin{equation} \label{appeqn: Fourier transform of spectral density in 2T dimensions}
        \mathcal{F}^{-1}[\hat{\kappa}_{2T}({\zeta} ; \nu_{2T}, l_{2T})]([z,0]) = \int_{\zeta_1, \zeta_2 \in \mathbb{R}^{2T}} \hat{\kappa}_{2T}({\zeta}) e^{i2\pi({\zeta} \cdot [z,0])} d{\zeta} = \int_{\zeta_1}\left[\int_{\zeta_2}\hat{\kappa}_{2T}(\zeta_1, \zeta_2; \nu_{2T}, l_{2T}) d\zeta_2\right] e^{i 2\pi z\cdot\zeta_1} d\zeta_1.
    \end{equation}
Which is expressed more compactly as,
    \begin{equation}
       \mathcal{F}_{\zeta_1}^{-1} \left[\int_{\zeta_2}\hat{\kappa}_{2T}(\zeta_1, \zeta_2; \nu_{2T}, l_{2T}) d\zeta_2 \right] (z).
    \end{equation}
    
Since Fourier transforms are invertible linear operators over the space of functions containing the Mat\'ern kernels, we simply need the kernels (and the coefficients) within the transforms in equation \eqref{appeqn: 2T to T Fourier Transform} to match in order to ensure that the entire equivalency remains true. This is because two functions are the same if and only if their Fourier transforms match. This leads us to the following, desired, equivalency:
    \begin{equation}\label{appeqn: Kernel equivalency by integration}
        \int_{\zeta_2}\hat{\kappa}_{2T}(\zeta_1, \zeta_2; \nu_{2T}, l_{2T}) d\zeta_2 = \alpha\cdot\hat{\kappa}_{T}(\zeta_1; \nu_T, l_T)
    \end{equation}.

Substituting the power spectral density for the kernel leads to
\begin{equation} \label{appeqn: expanded kernel equivalency}
    \begin{aligned} 
        \frac{\Gamma(\nu_{2T}+\tfrac{2T}{2})}{\Gamma(\nu_{2T})}2^{2T}\pi^{2T/2} \left(\frac{l_{2T}^{2}}{2\nu_{2T}}\right)^{2T/2}\int_{\zeta_2} & \left(1+\frac{l_{2T}^2\cdot
        ||2\pi{\zeta}||^2}{2\nu_{2T}}\right)^{-\left(\nu_{2T}+\tfrac{2T}{2}\right)}d\zeta_2 =  
        \\
        &\alpha \cdot \frac{\Gamma(\nu_T+\tfrac{T}{2})}{\Gamma(\nu_T)}2^{T}\pi^{T/2} \left(\frac{l_{T}^2}{2\nu_T}\right)^{T/2}\left(1+\frac{l_{T}^2\cdot 
        ||2\pi{\zeta}||^2}{2\nu_T}\right)^{-\left(\nu_{T}+\tfrac{T}{2}\right)}.
    \end{aligned}
\end{equation}

Dividing across by matching terms produces:
\begin{equation}\label{appeqn: isolating integral}
\begin{aligned}
    \int_{\zeta_2} \left(1+\frac{l_{2T}^2\cdot
    ||2\pi{\zeta}||^2}{2\nu_{2T}}\right)&^{-\left(\nu_{2T}+\tfrac{2T}{2}\right)}d\zeta_2 
    \\= &\alpha \cdot (2\sqrt{\pi})^{-T}\left(\frac{\Gamma(\nu_{2T})}{\Gamma(\nu_T)}\frac{\Gamma(\nu_T+\tfrac{T}{2})}{\Gamma(\nu_{2T}+\tfrac{2T}{2})}\right)\left(\frac{l_T}{l_{2T}^2}\right)^{T}\left(\frac{(2\nu_{2T})^2}{2\nu_T}\right)^{\frac{T}{2}}\left(1+\frac{l_{T}^2\cdot 
    ||2\pi{\zeta}||^2}{2\nu_T}\right)^{-\left(\nu_{T}+\tfrac{T}{2}\right)}.
\end{aligned}
\end{equation}

Now attempting to simplify the integral on the left hand side:
    \begin{equation} \label{appeqn: simplifying the integral}
    \begin{aligned}
        \int_{\zeta_2 \in \mathbb{R}^{T}} \left(1+\left(\frac{l_{2T}^2}{2\nu_{2T}}\right)(2\pi)^2(||\zeta_1||^2 + ||\zeta_2||^2)\right)^{-\left(\nu_{2T}+\tfrac{2T}{2}\right)}d\zeta_2 & = \int_{\zeta_2 \in \mathbb{R}^{T}}(a+b||\zeta_2||^2)^{-\left(\nu_{2T}+\tfrac{2T}{2}\right)}d\zeta_2  \\
        & = \int_{\zeta_2 \in \mathbb{R}^{T}} \frac{1}{a^{\left(\nu_{2T}+\tfrac{2T}{2}\right)}(1+\frac{b}{a}||\zeta_2||^2)^{\left(\nu_{2T}+\tfrac{2T}{2}\right)}} d\zeta_2 
        \\
        & = a^{-\left(\nu_{2T}+\tfrac{2T}{2}\right)}\int_{\zeta_2 \in \mathbb{R}^{T}}(1+c^2||\zeta_2||^2)^{-\left(\nu_{2T}+\tfrac{2T}{2}\right)}d\zeta_2.
    \end{aligned}
    \end{equation} 
where $a = (1+\left(\frac{l_{2T}^2}{2\nu_{2T}}\right)(2\pi)^2||\zeta_1||^2)$, $b = \left(\frac{l_{2T}^2}{2\nu_{2T}}\right)(2\pi)^2$, and $c^2 = \frac{b}{a}$.

Converting into Polar Coordinates yields the following, where $A_T$ is the surface area of a sphere in $T$-dimensions with a unit radius. If we let $u = cr \rightarrow du = cdr$:
    \begin{equation} \label{appeqn: converting to polar coordinates}
    \begin{aligned}
        \left(\frac{1}{2}A_T \cdot a^{-\left(\nu_{2T}+\tfrac{2T}{2}\right)}\right) \int_{r = -\infty}^{\infty}(1+(cr)^2)^{-\left(\nu_{2T}+\tfrac{2T}{2}\right)}r^{T-1} &dr 
        \\ 
        &= \left(\frac{1}{2}A_T \cdot a^{-\left(\nu_{2T}+\tfrac{2T}{2}\right)} c^{-T} \right)\int_{u = -\infty}^{\infty} \frac{u^{T-1}}{(1+u^2)^{\left(\nu_{2T}+\tfrac{2T}{2}\right)}}du.
    \end{aligned}
    \end{equation} 

This integral is proportional to the $(T-1)^{st}$ moment of the Student's-t distribution. 
    \begin{equation}\label{appeqn: proportional to Student's-t distribution}
        \int_{u = -\infty}^{\infty} \frac{u^{T-1}}{(1+u^2)^{\left(\nu_{2T}+T\right)}}du \propto \mathbb{E}[u^{T-1}]
    \end{equation} 

Let $\nu_{2T} +T = \frac{\nu + 1}{2} \rightarrow \nu = 2\nu_{2T} +2T-1$ and call $u^2 = \frac{x^2}{\nu} \rightarrow u = \frac{x}{\sqrt{\nu}}$
    \begin{equation}\label{appeqn: Making our integral look like a Student's-t}
        \int_{-\infty}^{\infty} \left(\frac{x}{\sqrt{\nu}}\right)^{T-1}\frac{1}{(1+\frac{x^2}{\nu})^{\left(\nu_{2T}+T\right)}}\frac{dx}{\sqrt{\nu}} = \left(\frac{1}{\sqrt{\nu}}\right)^T\int_{-\infty}^{\infty} \frac{x^{T-1}}{(1+\frac{x^2}{\nu})^{\frac{\nu+1}{2}}}dx 
\end{equation}

Meanwhile, the $(T-1)^{st}$ moment of X, which follows the Student's-t distribution with parameter $\nu$, is as follows. The following holds if $(T-1)$ is even and if $(T-1) < \nu$: 
    \begin{equation}\label{Finding the moment of student t}
       \mathbb{E}[X^{T-1}] = \int_{-\infty}^{\infty} \frac{\Gamma(\frac{\nu+1}{2})}{\sqrt{\pi\nu}\cdot\Gamma(\frac{\nu}{2})}\cdot x^{T-1}\left(1+\frac{x^2}{\nu}\right)^{-\frac{\nu+1}{2}}dx = \frac{1}{\sqrt{\pi}\cdot\Gamma(\frac{\nu}{2})}\cdot \Gamma\left(\frac{T}{2}\right)\cdot\Gamma\left(\frac{\nu-(T-1)}{2}\right)\cdot \nu^{\frac{T-1}{2}}
    \end{equation}
Thus,
    \begin{equation} \label{appeqn: Solving the integral using moments of the student t distribution}
    \begin{aligned}
       \left(\frac{1}{\sqrt{\nu}}\right)^T\int_{-\infty}^{\infty} \frac{x^{T-1}}{(1+\frac{x^2}{\nu})^{\left(\nu_{2T}+T\right)}}dx &= \left(\frac{1}{\sqrt{\nu}}\right)^T\cdot\frac{\sqrt{\pi\nu}\Gamma(\frac{\nu}{2})}{\Gamma(\frac{\nu+1}{2})}\cdot\frac{1}{\sqrt{\pi}}\cdot\frac{1}{\Gamma(\frac{\nu}{2})}\cdot\Gamma\left(\frac{T}{2}\right)\cdot\Gamma\left(\frac{\nu-(T-1)}{2}\right)\cdot\nu^{\frac{T-1}{2}}\\
       \\
       & = \frac{\Gamma(\frac{T}{2})\Gamma(\frac{(\nu+1)-T}{2})}{\Gamma(\frac{\nu+1}{2})}
       = \frac{\Gamma(\frac{T}{2})\Gamma(\nu_{2T}+\frac{T}{2})}{\Gamma(\nu_{2T} +T)}
    \end{aligned}
    \end{equation}

More compactly, this reads:
\begin{equation} \label{appeqn: integral equivalency as a ratio of gammas}
     \int_{u = -\infty}^{\infty} \frac{u^{T-1}}{(1+u^2)^{\left(\nu_{2T}+T\right)}}du =\frac{\Gamma(\frac{T}{2})\Gamma(\nu_{2T}+\frac{T}{2})}{\Gamma(\nu_{2T} +T)}
\end{equation} 

Re-substituting $a$ and $c$ back into equation \eqref{appeqn: converting to polar coordinates} produces:
    \begin{equation}
        \left(\frac{1}{2}A_T \cdot \frac{\left(1+\left(\frac{l_{2T}^2}{2\nu_{2T}}\right)(2\pi)^2||\zeta_1||^2\right)^{-\left(\nu_{2T}+\frac{T}{2}\right)}}{\left(\frac{l_{2T}^2}{2\nu_{2T}}\right)^{\frac{T}{2}}(2\pi)^T}\right) \left(\frac{\Gamma(\frac{T}{2})\Gamma(\nu_{2T}+\frac{T}{2})}{\Gamma(\nu_{2T} +T)}\right).
    \end{equation}

We can check both sides of \eqref{appeqn: expanded kernel equivalency} to determine the decay rate to zero as a function of $\zeta_1$. The left hand side of \eqref{appeqn: expanded kernel equivalency}, or what is shown in \eqref{appeqn: converting to polar coordinates}, has rate $\mathcal{O}\left(||\zeta_1||^{-\left(\nu_{2T}+\frac{T}{2}\right)}\right)$, while the right hand side of \eqref{appeqn: expanded kernel equivalency} decays at the same rate. Clearly, these to go to zero at the same speed, and we get our first equivalence:
\begin{equation} \label{appeqn: nu_T equivalency}
    \nu_T = \nu_{2T}
\end{equation} 
    
To match the second order Taylor expansion about $||\zeta_1||^2 = 0$, we would need,
    \begin{equation} \label{appeqn: solving for l_T equivalency}
        \left(\frac{l_{2T}^2}{2\nu_{2T}}\right) = \left(\frac{l_{T}^2}{2\nu_{T}}\right). 
    \end{equation} 
Since $\nu_T = \nu_{2T}$, his leads to our second equivalence,
    \begin{equation} \label{appeqn: l_T equivalency}
        l_T = l_{2T}.
    \end{equation} 

Finally, with $A_T = \frac{2\pi^{\frac{T}{2}}}{\Gamma(\frac{T}{2})}$, we get
$$    
    \begin{aligned}
    \int_{\zeta_2} \left(1+\frac{l_{2T}^2\cdot
    ||2\pi{\zeta}||^2}{2\nu_{2T}}\right)&^{-\left(\nu_{2T}+\tfrac{2T}{2}\right)}d\zeta_2 \\&= \left(\frac{1}{2}A_T \cdot \frac{\left(1+\left(\frac{l_{2T}^2}{2\nu_{2T}}\right)(2\pi)^2||\zeta_1||^2\right)^{-\left(\nu_{2T}+\frac{T}{2}\right)}}{\left(\frac{l_{2T}^2}{2\nu_{2T}}\right)^{\frac{T}{2}}(2\pi)^T}\right) \left(\frac{\Gamma(\frac{T}{2})\Gamma(\nu_{2T}+\frac{T}{2})}{\Gamma(\nu_{2T} +T)}\right) \\
        &=\alpha \cdot (2\sqrt{\pi})^{-T}\left(\frac{\Gamma(\nu_{2T})}{\Gamma(\nu_T)}\frac{\Gamma(\nu_T+\tfrac{T}{2})}{\Gamma(\nu_{2T}+\tfrac{2T}{2})}\right)\left(\frac{l_T}{l_{2T}^2}\right)^{T}\left(\frac{(2\nu_{2T})^2}{2\nu_T}\right)^{\frac{T}{2}}\left(1+\frac{l_{T}^2\cdot 
    ||2\pi{\zeta}||^2}{2\nu_T}\right)^{-\left(\nu_{T}+\tfrac{T}{2}\right)}.
    \end{aligned}
$$
This reduces to
    \begin{equation}
    \begin{aligned}
        \left(\frac{1}{2}A_T \cdot \left(\frac{l^2}{2\nu}\right)^{\frac{-T}{2}}(2\pi)^{-T}\right) = \alpha\cdot2^{-T}\pi^{-\frac{T}{2}}\frac{1}{\Gamma{\left(\tfrac{T}{2}\right)}}l^{-T}(2\nu)^{\frac{T}{2}}.
    \end{aligned}
    \end{equation}
Which allows us to isolate and solve for $\alpha$:
    \begin{equation}
    \begin{aligned} \label{appeqn: solving for alpha equivalency}
        \alpha = \frac{\frac{1}{2}A_T\left(\frac{l^2}{2\nu}\right)^{\frac{-T}{2}}(2\pi)^{-T}}{2^{-T}\pi^{-\frac{T}{2}}\frac{1}{\Gamma{\left(\tfrac{T}{2}\right)}}l^{-T}(2\nu)^{\frac{T}{2}}} = \frac{\frac{1}{2}A_T\Gamma{\left(\tfrac{T}{2}\right)}}{\pi^{\frac{T}{2}}} = \frac{\Gamma{\left(\tfrac{T}{2}\right)}}{2\pi^{\frac{T}{2}}} \cdot \frac{2\pi^{\frac{T}{2}}}{\Gamma(\frac{T}{2})} =1.
    \end{aligned}
    \end{equation}
    
This ultimately ends our chain of equivalences.
    \begin{equation} \label{appeqn: dimension conversion factor equivalency}
    \begin{aligned}
        &\nu_T =  \nu_{2T}\\
        &l_T = l_{2T}\\
        &\alpha = 1.
    \end{aligned}
    \end{equation} 
Therefore, we have the appropriate relationship of kernels and can proceed with the expected value calculations. The Mat\'ern kernel is conveniently parameterized so that this exchange is trivial - all parameters remain the same:
    \begin{equation} \label{appeqn: kernel dimension relationship restatement}
        \kappa_{2T}([Z,0]; \nu_{2T}, l_{2T}) = \kappa_{T}(Z; \nu_{T}, l_{T}) = \mathcal{F}^{-1}[\hat{\kappa}_{T}(\zeta ; \nu, l)](z).
    \end{equation}

\vspace{0.1 in}
\noindent \textbf{Roughness Limit}\label{app: roughness limit calculations}
\vspace{0.05 in}

Here, we will compute the correlation $\rho$ in the roughness limit. Recall that $\rho$ is defined as:
\begin{equation} \label{appeqn: Rho equation for roughness}
    \begin{aligned}
       \rho = \frac{\mathbb{E}_{Z \sim \mathcal{N}(0,2 \sigma_x^2 I)}\left[k_{2T}([Z,0]; \nu, l) \right] - \mathbb{E}_{Z^{(1)},Z^{(2)}\sim \mathcal{N}(0, \sigma_x^2 S^{-1})}\left[ k_{2T}([Z^{(1)},Z^{(2)}]; \nu, l)  \right]}{1 - \mathbb{E}_{Z \sim \mathcal{N}(0,4 \sigma_x^2 I)}\left[k_{2T}([Z,0]; \nu, l)\right] }
    \end{aligned}
\end{equation}
where:
$$
S  = \frac{1}{3} \left[\begin{array}{cc} 2 & -1 \\ -1 & 2 \end{array} \right] \otimes I_{T \times T} = \left( \left[\begin{array}{cc} 2 & 1 \\ 1 & 2 \end{array} \right] \otimes I_{T \times T} \right)^{-1} = \left[\begin{array}{cc} 2 I_{T \times T} & I_{T \times T} \\ I_{T \times T} & 2 I_{T \times T} \end{array} \right]^{-1}
$$

We first convert each expectation to frequency space via a Fourier Transform. Recall that the correct conversion for a Gaussian expectation is given by:
\begin{equation} \label{appeqn: gaussian expectation in frequency space}
    \mathbb{E}_{Z \sim \mathcal{N}(0,\Sigma_z)}[g(Z)] = \frac{1}{(2 \pi)^{d/2}| \text{det}(\Sigma_z)|^{1/2}} \mathbb{E}_{\zeta \sim \mathcal{N}(0, (2 \pi)^{-2} \Sigma_z^{-1})}[\hat{g}(\zeta)].
\end{equation}
where $d$ reflects the dimension. 

Therefore, the expectations required are:
\begin{equation} \label{appeqn: necessary expectations in frequency space}
    \begin{aligned}
        & \mathbb{E}_{Z \sim \mathcal{N}(0,2 \sigma_x^2 I)}\left[k_T(Z;\nu;l) \right] = ((2 \pi) 2 \sigma_x^2)^{-T/2}  \mathbb{E}_{\zeta \sim \mathcal{N}(0, (2\pi)^{-2} (2 \sigma_x^{2})^{-1} I)}[\hat{k}_T(\zeta;\nu,l)] \\
        & \mathbb{E}_{Z \sim \mathcal{N}(0,2 \sigma_x^2 I)}\left[k_T(\sqrt{2} Z;\nu;l) \right] = \mathbb{E}_{Z \sim \mathcal{N}(0,4 \sigma_x^2 I)}\left[k_T(Z;\nu;l) \right] \\& \hspace{2.135 in}  = ((2 \pi) 4 \sigma_x^2)^{-T/2} \mathbb{E}_{\zeta \sim \mathcal{N}(0,(2\pi)^{-2} (4 \sigma_x^{2})^{-1} I)}[\hat{k}_{T}(\zeta;\nu,l)] \\
        & \mathbb{E}_{Z^{(1)},Z^{(2)}}\left[ k_u([Z^{(1)},Z^{(2)}]) \right] = ((2 \pi) \sqrt{3} \sigma_x^2)^{-T} \mathbb{E}_{\zeta \sim \mathcal{N}(0,(2\pi)^{-2} \sigma_x^{-2} S)}[\hat{k}_{2T}(\zeta;\nu,l)]  
    \end{aligned}
\end{equation}
where:
$$
S  = \frac{1}{3} \left[\begin{array}{cc} 2 & -1 \\ -1 & 2 \end{array} \right] \otimes I_{T \times T} = \left( \left[\begin{array}{cc} 2 & 1 \\ 1 & 2 \end{array} \right] \otimes I_{T \times T} \right)^{-1} = \left[\begin{array}{cc} 2 I_{T \times T} & I_{T \times T} \\ I_{T \times T} & 2 I_{T \times T} \end{array} \right]^{-1}
$$
is the inverse covariance for $[Z^{(1)},Z^{(2)}]$, and where the factor of $3$ inside the scaling comes from the determinant of the covariance of $[Z^{(1)},Z^{(2)}]$.

Thus, the expectations we need to compute are (sans the previously computed normalizing factor) below. The first two expectations follow the same pattern. We will use $n$ to represent the 2 and 4 evident in their respective covariance matrices:
\begin{equation} \label{appeqn: expectations in frequency space}
    \begin{aligned}
         &\mathbb{E}_{\zeta \sim \mathcal{N}(0, (2\pi)^{-2} (n \sigma_x^{2})^{-1} I)}[\hat{k}_T(\zeta;\nu,l)] \\ &\hspace{0.5 in}=  \mathbb{E}_{\zeta \sim \mathcal{N}(0, (2\pi)^{-2} (n \sigma_x^{2})^{-1} I)}\left[\frac{\Gamma(\nu)}{\Gamma(\nu - \tfrac{T}{2})}2^T \pi^{\frac{T}{2}} \left(\frac{l^2}{2\nu - T} \right)^{\frac{T}{2}}\left(1+ \frac{l^2||2\pi\zeta||^2}{2\nu-T}\right)^{-\nu}\right]\\
         &\mathbb{E}_{\zeta \sim \mathcal{N}(0,(2\pi)^{-2} \sigma_x^{-2} S)}[\hat{k}_{2T}(\zeta;\nu,l)] \\ &\hspace{0.5 in} =  \mathbb{E}_{\zeta \sim \mathcal{N}(0, (2\pi)^{-2} \sigma_x^{-2} S)}\left[\frac{\Gamma(\nu+\tfrac{T}{2})}{\Gamma(\nu - \tfrac{T}{2})}2^{2T} \pi^{T} \left(\frac{l^2}{2\nu-T} \right)^{T}\left(1+ \frac{l^2||2\pi\zeta||^2}{2\nu-T}\right)^{-\left(\nu + \tfrac{T}{2}\right)}\right].
    \end{aligned}
\end{equation}

\noindent Applying linearity:
\begin{equation} \label{appeqn: expectations in frequency space after linearity}
    \begin{aligned}
         &\mathbb{E}_{\zeta \sim \mathcal{N}(0, (2\pi)^{-2} (n \sigma_x^{2})^{-1} I)}[\hat{k}_T(\zeta;\nu,l)] \\ &\hspace{0.5 in} =  \left[\frac{\Gamma(\nu)}{\Gamma(\nu - \tfrac{T}{2})}2^T \pi^{\frac{T}{2}} \left(\frac{l^2}{2\nu - T} \right)^{\frac{T}{2}}\right]\mathbb{E}_{\zeta \sim \mathcal{N}(0, (2\pi)^{-2} (n \sigma_x^{2})^{-1} I)}\left[\left(1+ \frac{l^2||2\pi\zeta||^2}{2\nu-T}\right)^{-\nu}\right]\\
         &\mathbb{E}_{\zeta \sim \mathcal{N}(0,(2\pi)^{-2} \sigma_x^{-2} S)}[\hat{k}_{2T}(\zeta;\nu,l)] \\ &\hspace{0.5 in} =  \left[\frac{\Gamma(\nu+\tfrac{T}{2})}{\Gamma(\nu - \tfrac{T}{2})}2^{2T} \pi^{T} \left(\frac{l^2}{2\nu-T} \right)^{T}\right]\mathbb{E}_{\zeta \sim \mathcal{N}(0, (2\pi)^{-2} \sigma_x^{-2} S)}\left[\left(1+ \frac{l^2||2\pi\zeta||^2}{2\nu-T}\right)^{-\left(\nu + \tfrac{T}{2}
         \right)}\right].
    \end{aligned}
\end{equation}

We begin with the top expectation:
\begin{equation} \label{appeqn: Replacing sigma_x and l tilde}
    \begin{aligned}
        \mathbb{E}_{\zeta \sim \mathcal{N}(0, (2\pi)^{-2} (n \sigma_x^{2})^{-1} I)}\left[\left(1+ \frac{l^2||2\pi\zeta||^2}{2\nu-T}\right)^{-\nu}\right] = \mathbb{E}_{\zeta \sim \mathcal{N}(0, \sigma_\zeta^{2} I)}\left[\left(1+ ||2\pi\Tilde{l}\zeta||^2\right)^{-\nu}\right].
    \end{aligned}
\end{equation}
Here, we also replace $\sigma_x^2$ and its coefficients by the term $\sigma_\zeta^2$. We also let $\Tilde{l} = \frac{l}{\sqrt{2\nu-T}}$. If we let everything inside the norm, $(2\pi\Tilde{l}\zeta)$ equal $u$, we make another substitution:
\begin{equation} \label{appeqn: Converting to u}
    \begin{aligned}
        \mathbb{E}_{\zeta \sim \mathcal{N}(0, \sigma_\zeta^{2} I)}\left[\left(1+ ||2\pi\Tilde{l}\zeta||^2\right)^{-\nu}\right] = \mathbb{E}_{u \sim \mathcal{N}(0, \sigma_u^{2} I)}\left[\left(1+ ||u||^2\right)^{-\nu}\right].
    \end{aligned}
\end{equation}

Above, if $u = (2\pi\Tilde{l}\zeta)$ and $\zeta \sim \mathcal{N}(0, \sigma_\zeta^{2} I)$, then $\sigma_u^{2} = (2\pi\Tilde{l})^2\sigma_\zeta^2$. Retracing the iteration backwards leads to $\sigma_u^{2} = (2\pi\Tilde{l})^2\sigma_\zeta^2 = \frac{(2\pi\Tilde{l})^2}{(2\pi)^2 n \sigma_x^2} = \frac{1}{n}(\frac{\Tilde{l}}{\sigma_x})^2 = [n(\frac{\sigma_x}{\Tilde{l}})^2]^{-1}$. If we let $\Tilde{r}= \frac{\sigma_x}{\Tilde{l}}$ be a representation of roughness, this all is equal to $\sigma_u^{2} = [n\Tilde{r}^2]^{-1}$, where again, n = 2 or 4. 

We now turn to approximation of the expectation using Taylor series expansion to gain intuition on the limiting behavior of $\sigma_u^2$.
\begin{equation} \label{appeqn: Taylor expanding integral instead}
    \begin{aligned}
        \mathbb{E}_{u \sim \mathcal{N}(0, \sigma_u^{2} I)}\left[\left(1+ ||u||^2\right)^{-\nu}\right] &= \mathbb{E}_{W \sim \chi^2(T)}\left[\left(1+ \sigma_u^2W\right)^{-\nu}\right] \\
        &\approx \mathbb{E}_{W \sim \chi^2(T)}\left[\sum_{n=0}^\infty \frac{(-1)^n}{n!}\frac{\Gamma(\nu + n)}{\Gamma(\nu)}(\sigma_u^2 W)^n \right] \\
        &= \sum_{n=0}^m \binom{\nu+n-1}{n} (-\sigma_u^2)^n\mathbb{E}_{W \sim \chi^2(T)}[W^n] + \mathcal{O}(\sigma_u^{2(m+1)}\mathbb{E}_{W \sim \chi^2(T)}[W^{m+1}]).
    \end{aligned}
\end{equation}

Expanding $\mathbb{E}_{W \sim \chi^2(T)}[W^n]$ as moments of the chi-squared distribution with $T$ degrees of freedom produces $1 - Tb\sigma_u^2 + \mathcal{O}(\sigma_u^{4})$. Re-substituting $\sigma_u^2$ back into the expression, and taking the second-order Taylor expansion, we see: 
\begin{equation} 
    \begin{aligned}
        \mathbb{E}_{u \sim \mathcal{N}(0, \sigma_u^{2} I)}\left[\left(1+ ||u||^2\right)^{-\nu}\right] = 1 - T(\nu\sigma_u^2 + \mathcal{O}(\sigma_u^{4}) = 1 - \frac{T}{2}\frac{\nu}{\nu-\frac{T}{2}}n^{-1}\Tilde{r}^{-2} + \mathcal{O}(\Tilde{r}^{-4}).
    \end{aligned}
\end{equation}

Now, we continue with the second expectation. We once again make the substitution for $\Tilde{l} = \frac{l}{\sqrt{2\nu-T}}$ and the appropriate replacement for $u = 2\pi\Tilde{l}\zeta$:
\begin{equation} \label{appeqn: second integral initial substitutions}
    \begin{aligned}
        \mathbb{E}_{\zeta \sim \mathcal{N}(0, (2\pi)^{-2} \sigma_x^{-2} S)}\left[\left(1+ \frac{l^2||2\pi\zeta||^2}{2(\nu-\frac{T}{2})}\right)^{-\left(\nu+\frac{T}{2}\right)}\right] &= \mathbb{E}_{\zeta \sim \mathcal{N}(0, (2\pi)^{-2} \sigma_x^{-2} S)}\left[\left(1+ ||2\pi\Tilde{l}\zeta||^2\right)^{-\left(\nu+\frac{T}{2}\right)}\right] \\&= \mathbb{E}_{u \sim \mathcal{N}(0, \Tilde{l}^{2} \sigma_x^{-2} S)}\left[\left(1+ ||u||^2\right)^{-\left(\nu+\frac{T}{2}\right)}\right].
    \end{aligned}
\end{equation}

Exploiting the rotational symmetry of $u$ allows us to vastly simplify the covariance matrix of the Gaussian. We use the notation $\textit{diag}(1(\times T), \frac{1}{3}(\times T))$ to represent a diagonal matrix with $T$ diagonal entries equalling $1$ and $\frac{1}{3}$ respectively. This is equivalent to $S$ up to rotation. In future steps this will be further simplified to simply $D$.
\begin{equation} \label{appeqn: exploiting rotational symmetry}
    \begin{aligned}
        \mathbb{E}_{u \sim \mathcal{N}(0, \Tilde{l}^{2} \sigma_x^{-2} S)}\left[\left(1+ ||u||^2\right)^{-\left(\nu+\frac{T}{2}\right)}\right] = \mathbb{E}_{u \sim \mathcal{N}(0, \Tilde{l}^{2} \sigma_x^{-2} \textit{diag}(1(\times T), \frac{1}{3}(\times T))}\left[\left(1+ ||u||^2\right)^{-\left(\nu+\frac{T}{2}\right)}\right].
    \end{aligned}
\end{equation}

 As we did above, the next step is to approximate this expectation through a second order Taylor series expansion in $u$ about $0$. We do this as $(\frac{\Tilde{l}}{\sigma_x})^2$ gets small (i.e. we are computing a roughness limit).
\begin{equation} \label{appeqn: Taylor approximation for second integral}
    \begin{aligned}
        \mathbb{E}_{u \sim \mathcal{N}(0, \Tilde{l}^{2} \sigma_x^{-2} D)}\left[\left(1+ ||u||^2\right)^{-\left(\nu+\frac{T}{2}\right)}\right] \approx \mathbb{E}_{u \sim \mathcal{N}(0, \Tilde{l}^{2} \sigma_x^{-2} D)}\left[1 + 0 + \frac{1}{2}u^T H u + \mathcal{O}(u^3) \right].
    \end{aligned}
\end{equation}
Above, $H$ represents the Hessian matrix of $(1+ ||u||^2)^{-(\nu+\frac{T}{2})}$ at $u=0$. Since the covariance matrix of $u$ is diagonal, 
\begin{equation} \label{appeqn: Hessian inner product}
    \begin{aligned}
         \mathbb{E}_{u \sim \mathcal{N}(0, \Tilde{l}^{2} \sigma_x^{-2} D)}\left[1 + 0 + \frac{1}{2}u^T H u + \mathcal{O}(u^3) \right] \approx 1 + \frac{1}{2}\mathbb{E}_{u \sim \mathcal{N}(0, \Tilde{l}^{2} \sigma_x^{-2} D)}[u^T H u] = 1 + \frac{1}{2}\left(\langle H, \Tilde{l}^{2} \sigma_x^{-2} D\rangle\right).
    \end{aligned}
\end{equation}

In the following, $h_{jj}$ represents the $j^{th}$ diagonal entry of the Hessian matrix: 
\begin{equation} \label{Simplifying taylor approximation with Hessian}
    \begin{aligned}
          1 + \frac{1}{2}\left(\langle H, \Tilde{l}^{2} \sigma_x^{-2} D\rangle\right) &=  1 + \frac{1}{2}\left[\sum_{j=1}^{T} h_{jj}(\Tilde{l}^2\sigma_x^{-2}\cdot 1) + \sum_{j=T+1}^{2T} h_{jj}(\Tilde{l}^2\sigma_x^{-2}\cdot \frac{1}{3})\right] \\ &= 1 + \frac{1}{2}\sum_{j=1}^{T} (\partial_{u_j}^2[(1 + ||u||^2)^{-\left(\nu+\frac{T}{2}\right)}]|_{u=0})(\Tilde{l}^2\sigma_x^{-2}\cdot 1) \\ & \hspace{7mm}+ \sum_{j=T+1}^{2T} (\partial_{u_j}^2[(1 + ||u||^2)^{-(\nu+\frac{T}{2})}]|_{u=0})(\Tilde{l}^2\sigma_x^{-2}\cdot \frac{1}{3}) \\
          &= \left(1+ (\Tilde{l}^2\sigma_x^{-2})(\tfrac{2}{3}T)\cdot (\partial_{u_1}^2(1 + ||u||^2)^{-\left(\nu+\frac{T}{2}\right)}|_{u=0})\right) \\
          &= \left(1+ (\tfrac{2}{3}T)(\Tilde{l}\sigma_x^{-1})^{-2}\cdot (\partial_{u_1}^2(1 + u_1^2)^{-\left(\nu+\frac{T}{2}\right)}|_{u_1=0})\right) \\
          &= 1-\tfrac{2}{3}T (\Tilde{l}\sigma_x^{-1})^{-2} \left(2\nu+T\right).
    \end{aligned}
\end{equation}

Thus, we finally have an approximation for our original expectation. This is in the limit as $(\Tilde{l}^2\sigma_x^{-2}) = \Tilde{r}^{-2} \rightarrow 0$, which is the roughness limit:
\begin{equation} \label{appeqn: roughness limit of second expectation}
    \begin{aligned}
         \mathbb{E}_{\zeta \sim \mathcal{N}(0, (2\pi)^{-2} \sigma_x^{-2} S)}\left[\hat{k}_{2T}(\zeta;\nu,l)\right] &\approx \left[\frac{\Gamma(\nu+\tfrac{T}{2})}{\Gamma(\nu - \tfrac{T}{2})}2^{2T} \pi^{T} \left(\frac{l^2}{2\nu-T} \right)^{T}\right] \left(1-\tfrac{4}{3}(\nu+\frac{T}{2})T \left(\frac{\Tilde{l}}{\sigma_x}\right)^{2} + \mathcal{O}\left(\frac{\Tilde{l}}{\sigma_x}\right)^{4} \right) \\
         &= \frac{\Gamma(\nu+\tfrac{T}{2})}{\Gamma(\nu - \tfrac{T}{2})} (4\pi)^T\left(\frac{l^2}{2\nu-T}\right)^T\left(1- \frac{4}{3}\frac{T}{2}\frac{(\nu+\frac{T}{2})}{(\nu-\frac{T}{2})}\Tilde{r}^{-2} + \mathcal{O}(\Tilde{r}^{-4})\right).
    \end{aligned}
\end{equation}

Returning to \eqref{appeqn: Rho equation for roughness}, if we view $\rho$ as solely the ratio of the order terms of $\Tilde{r}$, we find that, in the limit as $\Tilde{r} \rightarrow \infty$, 

\begin{equation} \label{appeqn: Rho equation as ratio of order terms}
    \begin{aligned}
       \rho &\approx \frac{\mathcal{O}(\Tilde{r}^{-T}) - \mathcal{O}(\Tilde{r}^{-2T})}{1 - \mathcal{O}(\Tilde{r}^{-T})} \simeq (\mathcal{O}(\Tilde{r}^{-T}) - \mathcal{O}(\Tilde{r}^{-2T})) (1+ \mathcal{O}(\Tilde{r}^{-T})) \\ &\simeq \mathcal{O}(\Tilde{r}^{-T}) \simeq  2^{T/2}\frac{\Gamma(\nu)}{\Gamma(\nu - \tfrac{T}{2})} 2 ^{-T/2} \Tilde{r}^{-T} \left( 1 - \frac{T}{2}\frac{\nu}{\nu-\frac{T}{2}}\frac{1}{2}\Tilde{r}^{-2} + \mathcal{O}(\Tilde{r}^{-4}) \right) \\ &=
       \frac{\Gamma(\nu)}{\Gamma(\nu-\frac{T}{2})} \Tilde{r}^{-T}\mathcal{O}(\Tilde{r}^{-T-2}) = \frac{\Gamma(\nu)}{\Gamma(\nu-\frac{T}{2})} \left(\frac{\sigma_x}{\Tilde{l}}\right)^{-T} = 2^{-T/2} \frac{\Gamma(\nu)}{\Gamma(\nu-\frac{T}{2})(\nu-\frac{T}{2}) ^{T/2}}\left(\frac{\sigma_x}{l}\right)^{-T}.
    \end{aligned}
\end{equation}

\setcounter{result}{11}
\setcounter{limit_result}{0}

Therefore, letting $r = \frac{\sigma_x}{l}$, in the roughness limit we see,
\begin{snugshade}
\begin{limit_result} \label{app result 11: roughness}
\textbf{[Roughness Limit of $\rho$ under (\hyperref[table: Results + Assumptions]{1 - 5, 9})]} When the attributes are drawn independently of $\mathcal{F}$, the flow is mean zero, and is drawn from a difference of stationary, isotropic, GP utilities with a Mat\'ern kernels, in the limit as $r \rightarrow \infty$:
\begin{equation} \label{appeqn: roughness limit}
    \rho(r) \simeq \frac{\Gamma(\nu)}{\Gamma(\nu-\tfrac{T}{2})(\nu-\tfrac{T}{2})^{\tfrac{T}{2}}}(\sqrt{2}r)^{-T}+ \mathcal{O}\left(r^{-(T + 2)} \right)
\end{equation}
\end{limit_result}
\end{snugshade}

\noindent \textbf{Smoothness Limit: Naive Expansion} \label{app: naive expansion}

The smoothness limit is considerably more involved in the joint limit when both $\nu$ and $r$ are small. The smaller $\nu$, the less regular the ensemble of functions, so small $\nu$ and small $r$ correspond to functions that are simultaneously smooth and irregular. Formally, small $\nu$ presents a problem since, when $\nu$ is small, naive expansion methods fail. 

Consider, for example, the scale-mixture form of $\rho$: 
\begin{equation}
    \rho(r,\nu) = \frac{\mathbb{E}_{V \sim \text{ Gamma}^{-1}(\nu,\nu)}\left[ (1 + 2 r^2 V)^{-\frac{T}{2}} \right] - \mathbb{E}_{V \sim \text{ Gamma}^{-1}(\nu,\nu)}\left[ (1 + r^2 V)^{-\frac{T}{2}} (1 + 3 r^2 V)^{-\frac{T}{2}}\right]}{1 - \mathbb{E}_{V \sim \text{ Gamma}^{-1}(\nu,\nu)}\left[ (1 + 4 r^2 V)^{-\frac{T}{2}} \right]}.
\end{equation}

When $r$ is small, it is tempting to Taylor expand the argument in $r$, then expand each expectation into a linear combination of the raw moments of the inverse gamma distribution. This approach is invalid since the inverse gamma distribution only has finite moments up to degree $\nu$. Thus, while the Taylor expansion is valid inside the expectation, the integral and the limit implicit in the sequence of partial sums, cannot be exchanged. 

Naive Taylor expansion inside the expectation fails since the tails of inverse gamma distribution decay slowly, according to a power law of the form $v^{-(\nu + 1)}$. Similarly, the argument of the expectations decays according to a power law of the form $v^{-T/2}$ for large $v$. The slow decay of these tails contributes to the expectations even in the smoothness limit when $r$ vanishes. So, to accurately capture the asymptotics for small $r$, it is important to retain an approximation to both a bulk term, that accounts for the bulk of the mass of the inverse gamma distribution, and a tail term, that accounts for the slowly decaying tails.  

To illustrate the issue, consider the following approximation approaches:
\begin{enumerate}
    \item Do not adopt the scale mixture approach. Instead, express the original expectations as the expectation of a Gaussian density function evaluated at a Student's T random variable. When $r$ is small, the Student's T distribution is concentrated relative to the Gaussian. Taylor expand the Gaussian to second order about zero, then evaluate the expectations using the moments of the Student's T distribution.
    \item Adopt the scale mixture approach, and Taylor expand the argument of each expectation in the following equation \eqref{appeqn: Matern mixture simplified}  with respect to $r$. Equivalently, Taylor expand in the scale parameter $V$ about $V = 0$. 
    
    \begin{equation} \label{appeqn: Matern mixture simplified}
    \rho(r,\nu) = \frac{\mathbb{E}_{V \sim \text{ Gamma}^{-1}(\nu,\nu)}\left[ (1 + 2 r^2 V)^{-\frac{T}{2}} \right] - \mathbb{E}_{V \sim \text{ Gamma}^{-1}(\nu,\nu)}\left[ (1 + r^2 V)^{-\frac{T}{2}} (1 + 3 r^2 V)^{-\frac{T}{2}}\right]}{1 - \mathbb{E}_{V \sim \text{ Gamma}^{-1}(\nu,\nu)}\left[ (1 + 4 r^2 V)^{-\frac{T}{2}} \right]}.
\end{equation}

    \item Recall that any inverse gamma random variable $V$ can be expressed as $\Theta^{-1}$ for a gamma random variable $\Theta \sim \text{Gamma}(\nu,\nu)$. Rewrite the expectations as expectations over the inverse scale parameter $\Theta$, then Taylor expand about $\mathbb{E}[\Theta] = 1$. 
\end{enumerate}

These approaches do not produce consistent asymptotics. The ensuing approximations are, to second order in $r$:
\begin{equation} \label{appeqn: naive smoothness limits}
    \begin{aligned}
        & \text{1. }\rho(r) \simeq \frac{1}{2} \left[1 - \frac{T-1}{2} r^2 \right], 
        & \quad  \text{2. } \rho(r) \simeq \frac{1}{2} \left[1 - \frac{T-1}{2} \frac{\nu}{\nu - 2} r^2 \right], 
        & \quad  \text{3. } \rho(r) \simeq \frac{1}{2} \left[1 - \frac{T-1}{2} \frac{\nu + 3}{\nu + 1} r^2 \right]. \\
    \end{aligned}
\end{equation}
While these approximations share a generic form, they differ in their dependence on $\nu$. The first is equivalent to the squared exponential result, and is only valid in the regular limit, $\nu \rightarrow \infty$. The second and third both converge to the first in this limit, but differ for small $\nu$. The smaller $\nu$, the larger the rational $\nu$ dependent terms. This illustrates an intuitive effect. The less regular the ensemble of sample functions, the smaller the shared-endpoint correlation. However, the two approximations do not agree on how the correlation should decrease as $\nu$ decreases. The second diverges at $\nu = 2$, and is useless for $\nu < 2$. Thus, it can only give predictions for ensembles that contain second differentiable functions. It cannot explain the desired regimes, $\nu < 1$, $\nu = 1$, $\nu \in (1,2)$, and $\nu = 2$, where the degree of differentiability changes, and where the assumptions of \cite{cebra2023similarity} are relaxed. In numerical tests, the third is less accurate than the second and fails to provide accurate approximations for small $\nu$. 

These approximations all fail since they neglect the slowly decaying tails of either the PSD, or the scale mixture distribution. They all assume that, for small $r$, $\rho(r)$ approaches $1/2$ from below at rate $r^2$. This result is necessarily true when the sample functions are almost surely second differentiable \cite{cebra2023similarity}. However, when $\nu \leq 2$, sample functions are almost never second differentiable. In this case, the analysis in \cite{cebra2023similarity} does not apply, so $\rho(r)$ need not differ from $1/2$ by a decrement of $\mathcal{O}(r^2)$. We will see that, when $\nu \in (1,2)$, $\rho(r)$ approaches $1/2$ as $r$ vanishes, but does so slower than $\mathcal{O}(r^2)$. When $\nu = 1$, $\rho(r)$ approaches $1/2$ for small $r$, since sample functions can be approximated locally with linear functions and vanishing error, however, $\rho(r)$ approaches $1/2$ logarithmically. When $\nu < 1$, the sample draws do not admit accurate linear approximations, so $\rho(r)$ converges to a value smaller than $1/2$.

\vspace{0.1 in}
\noindent \textbf{Smoothness Limit: Bulk and Tail}\label{app: smoothess limit Matern}
\vspace{0.05 in}

Recall that we aim to compute a lower bound on the shared-endpoint correlation of the form:
\begin{equation} \label{appeqn: rho in integral functions}
    \rho(r,\nu) \gtrsim \tilde{\rho}(r,\nu) = \frac{J(\sqrt{2} r, \nu) - J(r, \nu)}{J(\sqrt{2} r, \nu)} = 1 - \frac{J(r,\nu)}{J(\sqrt{2} r, \nu)}
\end{equation}
where $J(r,\nu)$ is the integral function: $J(r,\nu) = 1 - \mathbb{E}_{V \sim \text{Gamma}^{-1}}(\nu,\nu)\left[(1 + r^2 V)^{-T/2} \right]$.

Equation \eqref{appeqn: rho in integral functions} is easier to analyze, since the asymptotics of $\tilde{\rho}(r,\nu)$ are fully determined by the asymptotics of a single integral function, $J(r,\nu)$. Specifically, we only need to derive accurate approximations to $J(r,\nu)$ in the limit as $r \rightarrow 0$.

To analyze the asymptotic behavior of the integral function, $J(r,\nu)$, we split the integral over $v > 0$ into two pieces, a bulk term: $v \in (0,v_*(r))$, and a tail term: $v \geq v_*(r)$ where $v_*(r)$ is an $r$-dependent split term. Next, we develop separate approximations for the bulk and the tail that are each asymptotically accurate as $r$ vanishes. 

The bulk and tail terms are:
\begin{equation}
J(r,\nu) = J_{\text{b}}(r,\nu) + J_{\text{t}}(r,\nu) \text{ where }
\begin{cases}
    & J_{\text{b}}(r,\nu) = \left(1 - \mathbb{E}[(1 + r^2 V)^{-T/2} \mid V < v_*(r)] \right) \text{Pr}(V < v_*(r)) \\
    & J_{\text{t}}(r,\nu) = \left(1 - \mathbb{E}[(1 + r^2 V)^{-T/2} \mid V \geq v_*(r)] \right) \text{Pr}(V \geq v_*(r)) \\
\end{cases}
\end{equation}
for $V \sim \text{Gamma}^{-1}(\nu,\nu)$.

To evaluate each term, we approximate the conditional expectations. The cumulative probabilities outside each expectation can be evaluated using the CDF of the inverse gamma. In particular:
\begin{equation}
    \text{Pr}(V < v_*(r)) = \frac{\Gamma(\nu,\nu v_*(r)^{-1})}{\Gamma(\nu)}
\end{equation}
where $\Gamma(\cdot,\cdot)$ is the incomplete gamma function. 

Consider the bulk term first. If $r$ is small, then $r^2 V$ is small for all sufficiently small $V$. In particular, if $v_*(r)$ is chosen so that $r^2 v_*(r) \rightarrow 0$ as $r \rightarrow 0$, then $r^2 V$ is small for all $V$ in the bulk term. So, to approximate the bulk, adopt $v_*(r)$ such that $v_*(r) r^2 \rightarrow 0$ as $r \rightarrow \infty$, then Taylor expand the argument of the expectation. In practice, we will choose $v_*(r) = r^{-(2 - \epsilon)}$ for $\epsilon > 0$, then will take $\epsilon$ to zero from above.

Taylor expanding:
\begin{equation}
    (1 + r^2 V)^{-T/2} = 1 - \frac{T}{2} r^2 V + \frac{T(T+2)}{4} r^4 V^2 + \mathcal{O}(r^6 V^3). 
\end{equation}

Therefore, the expectation inside the bulk term takes the form:
\begin{equation}
\begin{aligned}
    1 - \mathbb{E}[(1 + r^2 V)^{-T/2} \mid V < v_*(r)] & = \frac{T}{2} \mathbb{E}_{V}[V|V < v_*(r)] r^2 + \frac{T(T+2)}{4} \mathbb{E}_{V}[V^2|V < v_*(r)] r^4 + \mathcal{O}(r^6 V^3) \\
\end{aligned}
\end{equation}
and the bulk term is, asymptotically:
\begin{equation} \label{appeqn: bulk asymptotics}
    J_{\text{b}}(r,\nu) \simeq \frac{\Gamma(\nu,\nu v_*(r)^{-1})}{\Gamma(\nu)} \left( \frac{T}{2} \mathbb{E}_{V}[V|V < v_*(r)] r^2 + \frac{T(T+2)}{4} \mathbb{E}_{V}[V^2|V < v_*(r)] r^4 \right)
\end{equation}

Note that, unlike the naive expansion approaches discussed before, this expansion of the bulk into a linear combination of moments, is asymptotically exact, since the integral runs over a finite domain for all nonzero $r$, and since the argument expanded is analytic at zero, the Taylor expansion converges for all $v$ in the domain of integration provided $r^2 v_*(r)$ vanishes as $r \rightarrow 0$. Each \textit{conditional} moment remains finite, so the integral and the infinite sum implied by the Taylor series commute. Each conditional moment of the inverse gamma exists when $v_*(r)$ is finite, and can be evaluated exactly. Thus, equation \eqref{appeqn: bulk asymptotics} provides exact, analytic, asymptotics for the bulk term in the smoothness limit. For detailed analysis, see Appendix section \ref{app: smoothess limit Matern}.

To approximate the tail term we take the opposite approach. Instead of developing a Taylor series approximation to the argument of the expectation in small $r^2 v$, we develop a Laurent series approximation to the inverse gamma density in large $v$. This is valid since the density does not depend on $r$, so may be approximated assuming large $v$ when $v_*(r)$ diverges. We adopt this approach, since, when $r^2 v_*(r)$ converges to zero, the argument of the expectation is always evaluated at both small and large inputs when evaluating the tail term. To approximate the tail term, expand the inverse gamma density into its Laurent series in $v$:
\begin{equation}
    f_{V}(v) = \frac{\nu^{\nu}}{\Gamma(\nu)} v^{-(\nu + 1)} e^{-\nu v^{-1}} = \frac{\nu^{\nu}}{\Gamma(\nu)} v^{-(\nu + 1)}(1 - \nu v^{-1} + \mathcal{O}(v^{-2})) \simeq \frac{\nu^{\nu}}{\Gamma(\nu)} v^{-(\nu + 1)}.
\end{equation}

This approximation converges to the true density with errors order $\mathcal{O}(v_*(r)^{-(\nu + 2)})$. This approximation sets $e^{-\nu/v} \simeq 1$ for large $v$. More accurate approximations follow by including more terms in the Laurent series. 

Under this approximation, the conditional inverse gamma density is Pareto on the tail. So, the conditional distribution over $V > v_*(r)$ should be replaced with:
\begin{equation}
    f_{V|V > v_*(r)}(v) \simeq \nu v_*(r)^{\nu} v^{-(\nu + 1)}.
\end{equation}

Then, the tail term is, asymptotically:
\begin{equation} \label{appeqn: tail asymptotics}
    J_{\text{t}}(r,\nu) \simeq \left( 1 - \frac{ \Gamma(\nu,\nu v_*(r)^{-1})}{\Gamma(\nu)} \right)\left(1 - \mathbb{E}_{V \sim \text{Pareto}(v_*(r),\nu)} [(1 + r^2 V)^{-T/2}]\right)
\end{equation}
with error $\mathcal{O}(v_*(r)^{-(\nu + 2)})$ since $1 - (1 + r^2 v)^{-T/2}$ is bounded above and below by 1 and 0. If $v_*(r) = r^{-(2 - \epsilon)}$ then $\mathcal{O}(v_*(r)^{-(\nu + 2)}) = \mathcal{O}(r^{4 + 2 \nu - \delta})$ for some $\delta > 0$ that converges to zero as $\epsilon$ converges to zero. Therefore, the higher order terms in the Laurent series associated with the tail will be dominated by the fourth order terms in $r$ that appear in the bulk term for all $\nu > 0$ and sufficiently small $\epsilon$. Like the asymptotic approximation to the bulk, the asymptotic approximation to the tail can be evaluated in closed form.

Evaluating the expectations in equations \eqref{appeqn: bulk asymptotics} and \eqref{appeqn: tail asymptotics}, and setting $v_*(r) = r^{-(2 - \epsilon)}$ for some $\epsilon > 0$ gives the combined approximation:
\begin{equation} \label{appeqn: combining bulk and tail}
\begin{aligned}
    J(r,\nu) = & J_{\text{b}}(r.\nu) + J_{\text{t}}(r,\nu) \simeq \hdots \\
    & \textbf{bulk:  } \frac{1}{\Gamma(\nu)} \left[ \frac{T}{2} (\nu r^2) \Gamma(\nu - 1, \nu r^{2 - \epsilon}) - \frac{T(T+2)}{8} (\nu r^2)^2 \Gamma(\nu - 2, \nu r^{2 - \epsilon}) \right] + \hdots \\
    & \textbf{tail: } \frac{\gamma(\nu,\nu r^{2 - \epsilon})}{\Gamma(\nu)} \left[1 - \frac{\nu}{\nu + \tfrac{T}{2}} r^{-\epsilon \tfrac{T}{2}} {_2\text{F}_1} \left(\tfrac{T}{2}, \nu + \tfrac{T}{2}; \nu + \tfrac{T}{2} +1 ; - r^{-\epsilon} \right)  \right]
\end{aligned}
\end{equation}
where $\Gamma(\cdot,\cdot)$ is the upper incomplete gamma function, $\gamma(\cdot,\cdot)$ is the lower incomplete gamma function, and where ${_2\text{F}_1}(\cdot,\cdot;\cdot;\cdot)$ is the hypergeometric function. To understand the asymptotic behavior of $\rho$ in the smoothness limit, expand each special function in the limit of small $r$, then take $\epsilon$ to zero from above. Expanding yields:
\begin{equation} \label{appeqn: expanded tail and bulk}
\begin{aligned}
    & \textbf{bulk: } \frac{1}{\Gamma(\nu)} \left[ \frac{T}{2} (\nu r^2 ) \left(\Gamma(\nu - 1) - \frac{(\nu r^{2 - \epsilon})^{\nu - 1}}{\nu - 1} \right) - \frac{T(T+2)}{8} (\nu r^2)^2 \Gamma(\nu - 2) + \mathcal{O}(r^6)\right] \\
    & \textbf{tail: } \frac{\nu^{\nu - 1}}{\Gamma(\nu)} r^{(2 - \epsilon) \nu} \left[1 - r^{2 \nu + \epsilon(T - \nu)}\left(1 - \frac{T}{2} \frac{\nu}{\nu - 1} r^{\epsilon} + \frac{T(T+2)}{8} \frac{\nu}{\nu - 2} r^{2 \epsilon} + \mathcal{O}(r^{3 \epsilon}) \right) \right]\\
\end{aligned}
\end{equation}

Notice that, for arbitrarily small $\epsilon$ the bulk includes a term of order $r^2$, of order $r^{2(\nu - 1) + 2} = r^{2 \nu}$ and a term of order $r^4$. The middle term is of intermediate order. It's order depends on the shape parameter $\nu$. It matches the lowest order term in the tail, which is also of order $r^{2 \nu}$ for $\epsilon$ sufficiently small. Together, these terms account for the slow tail decay that was ignored by the naive analysis. 

If $\nu > 2$, then $2 \nu= 2 \nu > 4$, so the order second and fourth order terms dominate. In this case, the expansion of the bulk term matches the fourth order naive expansion based on moments of the inverse gamma. If $\nu = 2$, then the middle term is of order $r^4$. If $\nu \in (1,2)$, then it is of intermediate order between $r^2$ and $r^4$. If $\nu = 1$, then it is order $r^2$, and if $\nu < 1$ then it is of intermediate order between $r^0$ and $r^2$. In other words, the smaller $\nu$, the slower the intermediate term approaches zero. At $\nu = 2$ it decays as slowly as the fourth order term. At $\nu = 1$ it decays as slowly as the second order term. Thus, the dominating term, and the largest correction to the dominating term, depend on $\nu$. For $\nu > 2$, we can drop the intermediate term since it is dominated by the second and fourth order terms. When $\nu < 2$ the intermediate term supplants the forth order term. When $\nu < 1$ it supplants the second-order term. This produces distinct asymptotic regimes separated by the critical values $\nu = 1$, and $\nu = 2$, where the sample functions drawn from the corresponding ensemble gain almost sure differentiability, and almost sure second differentiability.

Taking $\epsilon$ to zero from above, and keeping only the terms in the bulk and tail that approach zero slowly enough to contribute to the limiting behavior of $\rho$ leaves:
\begin{equation} \label{appeqn: expanded tail and bulk by nu}
    \begin{aligned}
        & \textbf{bulk: } \begin{cases}
            \frac{1}{\Gamma(\nu)} \frac{T}{2}(\nu r^2) \left[ \Gamma(\nu - 1) - \frac{(\nu r^2)^{\nu - 1}}{\nu - 1}\right] = \mathcal{O}(r^2) + \mathcal{O}(r^{2 \nu}) & \text{ if } \nu \in (0,2) \setminus 1 \\
            T r^2 (\log(r^{-1}) - \gamma) = \mathcal{O}(r^2 \log(r^{-1})) + \mathcal{O}(r^2) & \text{ if } \nu = 1 \\
            \frac{1}{\Gamma(\nu)} \left[\frac{T}{2} \Gamma(\nu - 1) (\nu r^2) - \frac{T(T+2)}{8} \Gamma(\nu - 2) (\nu r^2)^2 \right] = \mathcal{O}(r^2) + \mathcal{O}(r^4) & \text{ if } \nu > 2
        \end{cases} \\
        & \textbf{tail: } \begin{cases}\frac{\nu^{\nu - 1}}{\Gamma(\nu)} r^{2 \nu} + \mathcal{O}(r^{4 \nu}) = \mathcal{O}(r^{2 \nu}) & \text{ if } \nu \leq 2 \\
        \text{negligible} & \text{ if } \nu > 2
        \end{cases}
    \end{aligned}
\end{equation}

\end{document}